\documentclass[3p,times]{elsarticle}

\usepackage{graphicx}
\usepackage{psfrag}
\usepackage[dvipsnames]{xcolor}

\usepackage{amsmath}
\usepackage{amsthm}
\usepackage{amsfonts}
\usepackage{dsfont}
\usepackage{rotating}
\usepackage{multirow}
\usepackage{epstopdf}
\usepackage{newtxmath}
\usepackage{mathtools}
\usepackage{bm}
\usepackage{bbm}
\usepackage{comment}
\usepackage{booktabs}
\usepackage{accents}
\usepackage{empheq}

\usepackage{tikz}
\usepackage{pgfplots}
\usepackage{pgffor}
\usepackage{booktabs}
\usepackage{rotating}
\usepackage{lscape}
\usepackage{graphics}
\usepackage{dutchcal}

\newtheorem*{remark}{Remark}

\DeclareMathAlphabet\mathbfcal{OMS}{cmsy}{b}{n}  

\newtheorem{theorem}{Theorem}

\newcommand{\A}{\mathbf{A}} 
\newcommand{\J}{\mathbf{J}}
\newcommand{\detA}{|\A|} 
\newcommand{\dt}{\Delta t}  
\newcommand{\xx}{\mathbf{x}}  
\renewcommand{\vv}{\mathbf{v}} 
\newcommand{\U}{\mathbf{U}}  
\newcommand{\nn}{\nonumber}
\newcommand{\q}{\mathbf{q}}
\newcommand{\p}{\mathbf{p}}
\renewcommand{\r}{\mathbf{r}}
\renewcommand{\u}{\mathbf{u}}
\newcommand{\w}{\mathbf{w}}
\renewcommand{\vv}{\mathbf{v}}  
\newcommand{\halb}{\frac{1}{2}}
\newcommand{\err}{\mathcal{r}}
\newcommand{\nv}{\mathbf{n}}
\newcommand{\rhoS}{\mathcal{S}}

\def\be{\begin{equation}}
	\def\ee{\end{equation}}
\newcommand{\beann}[0]{\begin{eqnarray*}}
	\newcommand{\eeann}[0]{\end{eqnarray*}}
\def\bea{\begin{eqnarray}}
	\def\eea{\end{eqnarray}}
\def\ba{\begin{array}{l}\displaystyle}
	\def\ea{\end{array}}

\newcommand{\Frac}{\displaystyle\frac}
\newcommand{\Int}{\displaystyle\int}

\begin{document}
	
	\begin{frontmatter}
		
		\journal{Journal of Computational Physics}
		
		

		\title{A geometrically and thermodynamically compatible finite volume scheme for continuum mechanics on unstructured polygonal meshes}
		
		        
        \author[usb,dmi]{Walter Boscheri$^*$}
        \ead{walter.boscheri@unife.it}
        \cortext[cor1]{Corresponding authors}
        
        \author[imb]{Rapha{\"e}l Loub{\`e}re}
        \ead{raphael.loubere@math.u-bordeaux.fr}
        
        \author[cea]{Jean-Philippe Braeunig}
        \ead{jean-philippe.braeunig@cea.fr}
        
        \author[cea]{Pierre-Henri Maire}
        \ead{pierre-henri.maire@cea.fr}
        
        \address[usb]{LAMA UMR-5127 CNRS, Universit{\'e} Savoie Mont Blanc, 73376 Le Bourget du Lac, France}
        
		\address[dmi]{Department of Mathematics and Computer Science, University of Ferrara, 44121 Ferrara, Italy}
		
		\address[imb]{Institut de Math\'{e}matiques de Bordeaux (IMB), Universit{\'e} de Bordeaux, CNRS, Bordeaux INP, 33400 Talence, France}
		
		\address[cea]{CEA CESTA, 33116 Le Barp, France}
%
\begin{abstract}
We present a novel Finite Volume (FV) scheme on unstructured polygonal meshes that is provably compliant with the Second Law of Thermodynamics and the Geometric Conservation Law (GCL) at the same time. The governing equations are provided by a subset of the class of symmetric and hyperbolic thermodynamically compatible (SHTC) models introduced by Godunov in 1961. Specifically, our numerical method discretizes the equations for the conservation of momentum, total energy, distortion tensor and thermal impulse vector, hence accounting in one single unified mathematical formalism for a wide range of physical phenomena in continuum mechanics, spanning from ideal and viscous fluids to hyperelastic solids. By means of two conservative corrections directly embedded in the definition of the numerical fluxes, the new schemes are proven to satisfy two extra conservation laws, namely an entropy balance law and a geometric equation that links the distortion tensor to the density evolution. As such, the classical mass conservation equation can be discarded. Firstly, the GCL is derived at the continuous level, and subsequently it is satisfied by introducing the new concepts of general potential and generalized Gibbs relation. The new potential is nothing but the determinant of the distortion tensor, and the associated Gibbs relation is derived by introducing a set of dual or thermodynamic variables such that the GCL is retrieved by dot multiplying the original system with the new dual variables. Once compatibility of the GCL is ensured, thermodynamic compatibility is tackled in the same manner, thus achieving the satisfaction of a local cell entropy inequality. The two corrections are orthogonal, meaning that they can coexist simultaneously without interfering with each other. The compatibility of the new FV schemes holds true at the semi-discrete level, and time integration of the governing PDE is carried out relying on Runge-Kutta schemes. A large suite of test cases demonstrates the structure preserving properties of the schemes at the discrete level as well. 
\end{abstract}
%
\end{frontmatter}

\section{Introduction} \label{sec:introduction}
%
Born in 1929, Godunov embarked on a remarkable career that spanned several decades until his death in July 2023. 
He made profound contributions to the development of numerical techniques for solving Partial Differential Equations (PDEs), particularly in the field of fluid dynamics. 
The celebrated Godunov theorem and Godunov scheme revolutionized Computational Fluid Dynamics (CFD), enabling scientists and engineers to simulate complex fluid flows.
However, Godunov research developed beyond fluid dynamics, and in \cite{God1961} he found a connection between symmetric hyperbolicity in the sense of Friedrichs \cite{FriedrichsSymm} and thermodynamic compatibility. In this work of Godunov, one learns that for hyperbolic systems having an underlying variational formulation,  the total energy conservation law can be derived as the dot product of the other equations with the so-called \textit{thermodynamic dual variables} that are given by the partial derivative of the total energy
potential with respect to the conservative variables of the system. These variables are also known as \textit{main field}, or even \textit{Godunov variables}, see for instance \cite{Freistuehler2019}. 
Later on, Godunov and Romenski  \cite{Rom1998,GodRom2003} extended the theory of symmetric hyperbolic and thermodynamically compatible (SHTC) systems to a wide class of mathematical models: magnetohydrodynamics \cite{God1972MHD}, nonlinear hyperelasticity \cite{GodRom1972}, compressible multi-phase flows \cite{RomenskiTwoPhase2007,RomenskiTwoPhase2010} as well as relativistic fluid and solid mechanics \cite{Godunov2012,GRGPR}.
Within this theory the total energy potential plays a crucial role which is coming from the variational principle from which the system is derived. Moreover, the entropy density equation is part of the master system, while the total energy conservation law is an extra conservation law, since it is obtained by a linear combination of the other conservation equations. Due to the very general formalism, SHTC systems have been derived for modeling a wide range of physical phenomena that cover magnetohydrodynamics \cite{godunov1972symmetric}, nonlinear hyperelasticity \cite{GodRom1972Nonstationary}, compressible multi-phase flows \cite{romenski2010conservative} and even relativistic fluid and solid mechanics \cite{godunov2012thermodynamic,romenski2020new}.

The SHTC models are therefore compliant with the Second Law of Thermodynamics by construction, and they are derived as first order hyperbolic systems, where the stress tensor is a function of the
inverse deformation gradient $\A$ rather than velocity gradients, even for fluids. Indeed, irreversible dissipative processes are accounted by the presence of source terms with one or more characteristic strain relaxation times $\tau$. The hyperbolicity of SHTC models implies finite wave speeds for all involved physical processes, even dissipative
ones, thus making their mathematical structure substantially different from those PDE systems which admit parabolic dissipation and diffusion terms. Indeed, in \cite{PeshRom2014}, heat conduction is derived in first order hyperbolic form proving consistency with the Fourier law in the asymptotic regime \cite{GPRmodel}. Likewise, the stress tensor in the SHTC model proposed in \cite{PeshRom2014} is asymptotically consistent with the Navier-Stokes model. 
The distortion tensor is defined as the inverse of the deformation gradient, hence it is defined by construction as the inverse of the Jacobian matrix associated to the Lagrange-Euler mapping between the Lagrangian (or material) to the Eulerian (or updated) configuration. Consequently, the distortion tensor accounts for the deformation and rotation of the matter subject to mechanical and thermal loads. A direct link exists between the scalar density and the distortion tensor at the continuous level, which is also known as Geometric Conservation Law in the Lagrangian formulation of the governing equations \cite{Despres2005,Maire2007}. As a consequence, the density equation in the original model is redundant at the continuous level \cite{Nikiforadis21}. This can be viewed as an internal consistency constraint. However, ensuring this compatibility at the discrete level is not obvious, and this is one goal of this work. 

As already mentioned, in the SHTC formalism the total energy equation plays the role of an extra conservation law that can be deduced from the other equations of the system at the continuous level. This means that the entropy balance law is part of the master system, and it becomes an equality in the absence of shock waves. Nevertheless, at the discrete level, the energy equation is typically solved, hence ensuring energy conservation and numerical stability in the energy norm, and a lot of research has been conducted in order to achieve thermodynamic compatibility, i.e. obtaining an entropy balance law as a consequence of the chosen discretization. This research line started from the pioneering work presented in \cite{Tadmor1}, with the aim of devising provably entropy preserving and entropy stable numerical schemes that has been further investigated in \cite{Gassner_entropy2021,Gallice2022,kuya2023kinetic,lin2023positivity,lei2022high,berthon2018simple,kopriva2022theoretical,kuzmin2022limiter}. Other important contributions to the design and implementation of entropy preserving and stable schemes can be found for instance in \cite{fjordholm2012arbitrarily,castro2013entropy,hiltebrand2014entropy,Kling2016_entropy,gaburro2023high,fjordholm2012accurate} and references therein. The numerical strategy proposed in \cite{Abgrall2018,AbgrallOeffnerRanocha} has been recently employed to construct a new family of thermodynamically compatible schemes in which the entropy inequality is solved instead of the energy \cite{SWETurbulence,DumbserHTC,HTCAbgrall,HTCMHD,ThomannDumbser2023}, hence strictly mimicking the SHTC framework at the discrete level. The numerical methods are provably energy preserving at the semi-discrete level thanks to a scalar correction factor that is directly embedded in the definition of the numerical fluxes, hence ensuring conservation.
In Lagrangian hydrodynamics, thermodynamically compatible schemes have been developed in order to obtain the total energy conservation and the satisfaction of an entropy inequality as a consequence of a compatible discretization of the equations of continuity, momentum and internal energy, see for instance \cite{CaramanaCompatible1,CaramanaCompatible2,Maire2020}. A recent attempt in directly solving the entropy inequality and obtaining conservation and stability in the energy has been forwarded in \cite{boscheri2023HTCLag}.

In this work, we make use of the general framework introduced in \cite{Abgrall2018} for the construction of thermodynamically compatible schemes. We choose to discretize the total energy conservation law and deduce the entropy equation as a consequence, hence implying that the entropy inequality is one extra conservation law that is verified at continuous level, but must also be fulfilled at discrete level.  
This choice is the classical one because it is simpler to monitor energy or temperature in experimental devices compared to measuring entropy variations. From the continuous point of view, choosing the total energy equation or the entropy one is totally equivalent.
%
Furthermore, the numerical scheme must feature a discrete compatibility with the entropy inequality, and a discrete internal consistency between the determinant of the inverse of the deformation gradient and the discrete mass equation. We propose to resort to the approach originally forwarded in \cite{Abgrall2018}, and subsequently used in \cite{HTCAbgrall} for achieving thermodynamic compatibility for the SHTC model presented in \cite{PeshRom2014}. Our novel idea is to define a new geometrical potential that plays the role of total energy in SHTC schemes, and consequently to derive the associated dual variables. In this way, another extra conservation law can be obtained which accounts for the geometric consistency, that is nothing but the Geometric Conservation Law written in the Eulerian framework. Up to the knowledge of the authors, no geometrically compatible schemes on fixed unstructured meshes are part of the state-of-the-art numerical schemes for continuum mechanics. In this work we will design a first order Finite Volume (FV) on unstructured two-dimensional polygonal grids that is compatible with the Second Law of Thermodynamics and with the GCL, meaning that two extra conservation laws are satisfied by the scheme at the same time. This will ultimately allow to discard the classical mass conservation equation since the density can be deduced by the geometric compatibility achieved by the numerical method.

The paper is organized in three main sections. In Section \ref{sec.pde} we introduce the governing equations, the extra conservation laws and the final reduced compatible continuous model that is derived. Section \ref{sec.numscheme} is devoted to the design of the numerical scheme, including two theorems that demonstrate the compatibility of the new methods at the semi-discrete level. \ref{app.gcl} contains all the details related to the compatibility of the reduced model with the geometric constraint in the framework of SHTC systems. The numerical experiments are gathered in the dedicated Section \ref{sec.results}, where we numerically prove that the structural properties of the continuous model are preserved at the discrete level. Finally, we draw some conclusions and an outlook to future developments in Section \ref{sec.concl}.

\section{Mathematical model} \label{sec.pde}

\subsection{Governing equations}
The governing equations are given by the unified first order hyperbolic model of continuum mechanics proposed in \cite{PeshRom2014} that belongs to the class of hyperbolic thermodynamically compatible (HTC) systems \cite{God1961,GodRom1972,Rom1998, GodRom2003}. Let us assume Einstein summation convention over repeated indices, and let us adopt bold symbols to label vectors and matrices. Following \cite{HTCAbgrall}, the mathematical model is written in three space dimensions with indices $1\leq i,k,m \leq 3$ as follows:
\begin{subequations}\label{eqn.GPR}
	\begin{align}
		\frac{\partial \rho}{\partial t} 
		+\frac{\partial (\rho v_k)}{\partial x_k} 
		- \frac{\partial}{\partial x_m} \!\! \left( \epsilon \frac{\partial \rho}{\partial x_m} \right) 
		&= 0, \label{eqn.conti} \\
		\frac{\partial \rho v_i}{\partial t}
		+\frac{\partial \left( \rho v_i v_k + p \, \delta_{ik}   
			+  \sigma_{ik} + \phi_{ik} \right)}{\partial x_k}  
		-\frac{\partial}{\partial x_m}  \left( \epsilon \frac{\partial \rho v_i}{\partial x_m} \right) 
		&= 0,  
		\label{eqn.momentum} \\ 
		\frac{\partial \rhoS}{\partial t}
		+\frac{\partial \left( \rhoS v_k  \, + \, \beta_k \right)}{\partial x_k} 
		- \frac{\partial}{\partial x_m}  \left( \epsilon \frac{\partial \rhoS}{\partial x_m} \right)  
		&= \Pi  + \dfrac{\alpha_{ik} \alpha_{ik} }{\theta_1(\tau_1) T}  + \dfrac{\beta_i \beta_i}{\theta_2(\tau_2) T}  \geq 0,   
		\label{eqn.entropy} \\	
		\frac{\partial A_{i k}}{\partial t}
		+\frac{\partial (A_{im} v_m)}{\partial x_k} 
		+ v_m \left(\frac{\partial A_{ik}}{\partial x_m}-\frac{\partial A_{im}}{\partial x_k}\right)   -\frac{\partial}{\partial x_m}  \left( \epsilon \frac{\partial A_{ik}}{\partial x_m} \right) 
		&=  -\dfrac{ \alpha_{ik} }{\theta_1(\tau_1)} ,\label{eqn.deformation} 
		\\
		\frac{\partial J_k}{\partial t}
		+\frac{\partial \left( J_m v_m + T \right)}{\partial x_k} 
		+ v_m \left(\frac{\partial J_{k}}{\partial x_m} -\frac{\partial J_{m}}{\partial x_k}\right)  
		-\frac{\partial}{\partial x_m} \left( \epsilon \frac{\partial J_{k}}{\partial x_m} \right)  
		&= -\dfrac{\beta_k}{\theta_2(\tau_2)}, 
		\label{eqn.heatflux} 
	\end{align}
\end{subequations}
with $t \in \mathds{R}^0_+$ being the time and $\xx=\{x_k\}$ denoting the spatial position vector. The state vector $\q=\{q_j\} = (\rho, \rho v_i, \rhoS, A_{ik}, J_k)^\top$ is composed of mass density $\rho>0$, momentum $\rho \vv = \{\rho v_i\}$, total entropy $\rhoS$, distortion tensor $\A=\{A_{ik}\}$ and thermal impulse $\J=\{J_k\}$. Furthermore, $p>0$ and $T>0$ denote the fluid pressure and temperature, respectively. The fluid is also characterized by a polytropic index $\gamma=c_p/c_v$, given as the ratio of specific heats at constant pressure and volume, namely $c_p$ and $c_v$, respectively. The above system also accounts for parabolic vanishing viscosity terms with the parameter $\epsilon > 0$, that yield a production contribution $\Pi$ in the entropy equation \eqref{eqn.entropy} given by
\begin{equation}
	\label{eqn.Pi}
	\Pi = \frac{\epsilon}{T} 
	\, (\partial_{x_m} q_i) \,
	\, (\partial^2_{q_i q_j} \mathcal{E}  ) 
	\, (\partial_{x_m} q_j) \geq 0.   
\end{equation}
The positivity of the production term is ensured by assuming that the Hessian of the total energy potential is at least positive semi-definite, i.e. $\partial^2_{q_i q_j} \mathcal{E} \geq 0$, therefore the physical entropy is increasing, in accordance with the Second Law of Thermodynamics. Here, $\mathcal{E} = \mathcal{E}_1+\mathcal{E}_2+\mathcal{E}_3+\mathcal{E}_4$ is the total energy of the system which is obtained as the sum of four terms \cite{GPRmodel}:
\begin{equation} 
	\mathcal{E}_1 = \frac{\rho^\gamma}{\gamma-1} e^{S/c_v}, \quad 
	\mathcal{E}_2 = \halb \rho v_i v_i, 	
	\quad 
	\mathcal{E}_3 = \frac{1}{4} \rho c_s^2 \mathring{G}_{ij} \mathring{G}_{ij}, 	
	\quad 
	\mathcal{E}_4 = \halb c_h^2 \rho J_i J_i,  	
	\label{eq:energies}
\end{equation}
where $\mathbf{G}=\{ G_{ik}\}:=\{A_{ij} A_{kj}\}$ represents the metric tensor and $\mathring{\mathbf{G}}=\{\mathring{G}_{ik}\} = \{{G}_{ik} - \frac{1}{3} \, G_{mm} \delta_{ik}\}$ denotes its trace-free part with $\delta_{ik}$ being the Kronecker delta. The first term $\mathcal{E}_1$ corresponds to the internal energy, for which we assume the ideal gas equation of state, then the kinetic energy is considered by $\mathcal{E}_2$, whereas $\mathcal{E}_3$ is the shear energy with the shear sound speed $c_s$, and the last term $\mathcal{E}_4$ takes into account the thermal energy with $c_h$ being the heat wave speed. 
Let us now introduce the set of thermodynamic dual variables $\p:=\partial_{\q} \mathcal{E}=\{p_j\} = \left( r, v_i, T, \alpha_{ik}, \beta_k \right)^T$ which are explicitly given by the derivative of the energy potential \eqref{eq:energies} with respect to the state vector $\q$, that is
\begin{equation}
	r = \partial_{\rho} \mathcal{E}, 
	\qquad 
	v_i = \partial_{\rho v_i} \mathcal{E}, 
	\qquad 
	T = \partial_{\rhoS} \mathcal{E}, 
	\qquad 
	\alpha_{ik} = \partial_{A_{ik}} \mathcal{E}, 
	\qquad 
	\beta_{k} = \partial_{J_{k}} \mathcal{E}. 
\end{equation} 
In the momentum equation \eqref{eqn.momentum}, the shear stress tensor  $\bm{\sigma}=\{ \sigma_{ik}\}$ and the thermal stress tensor $\bm{\phi}=\{ \phi_{ik}\}$ are defined in terms of the dual variables $\alpha_{ik}$ and $\beta_k$ as 
\begin{equation}
	\sigma_{ik} = A_{ji} \partial_{A_{jk}} \mathcal{E} = 
	A_{ji} \alpha_{jk} = \rho c_s^2 G_{ij} \mathring{G}_{jk}, 
	\qquad 
	\phi_{ik} = J_i \partial_{J_{k}} \mathcal{E} = J_i \beta_{k} = \rho c_h^2 J_i J_k.
\end{equation}
The work of the shear and thermal stress tensors $\chi_k$ as well as the heat flux $h_k$ are given by
\begin{equation}
	\chi_k = \partial_{\rho v_i} \mathcal{E} \, \left( A_{ji} \partial_{A_{jk}} \mathcal{E} + J_i \partial_{J_{k}} \mathcal{E} \right)= v_i \, ( \sigma_{ik} + \phi_{ik} ), \qquad
	h_k = \partial_{\rhoS} \mathcal{E} \, \partial_{J_k} \mathcal{E} = T \beta_k = \rho c_h^2 T J_k,
\end{equation}

Finally, the mathematical model \eqref{eqn.GPR} is also endowed with algebraic source terms which contain two positive functions $\theta_1(\tau_1)>0$ and $\theta_2(\tau_2)>0$ that depend on $\q$ and on the relaxation times $\tau_1>0$ and $\tau_2>0$ as follows:   
\begin{equation}
	\theta_1 = \frac{1}{3}  \rho z_1 \tau_1 \, c_s^2 \, \left| \mathbf{A} \right|^{-\frac{5}{3}},
	\qquad 
	\theta_2 = \rho z_2 \tau_2 \, c_h^2,
	\qquad
	z_1 = \frac{\rho_0}{\rho},
	\qquad 
	z_2 = \frac{\rho_0 T_0}{\rho \, T},
\end{equation}
with $ \rho_0 $ and $ T_0 $ being a reference density and a reference temperature, respectively, and $|\A|$ denoting the determinant of $\A$. The asymptotic limit of the model \eqref{eqn.GPR} has been analyzed in \cite{GPRmodel} at the continuous level and in \cite{LGPR} in the fully discrete setting, showing that for small relaxation times, i.e. when $\tau_1 \to 0$ and $\tau_2 \to 0$, the Navier-Stokes-Fourier limit is obtained. Indeed, the stress tensor $\sigma_{ik}$ and the heat flux $h_k$ tend to  
\begin{equation}
	\sigma_{ik} = -\frac{1}{6} \rho_0 c_s^2 \tau_1 \left( \partial_k v_i + \partial_i v_k
	- \frac{2}{3} \left( \partial_m v_m\right) \delta_{ik} \right), 
	\qquad 
	h_k = -  \rho_0 T_0 c_h^2 \tau_2 \partial_k T,
	\label{eqn.asymptoticlimit}
\end{equation}
with $\sigma_{ik}$ fulfilling Stokes hypothesis. In the asymptotic regime, the relaxation time $\tau_1$ is directly related to the viscosity of the fluid by $ \mu = \frac{1}{6} \rho_0 c_s^2 \tau_1 $. Analogously, there is a direct link between the relaxation time $\tau_2$ and the thermal conductivity coefficient which is explicitly given by $\kappa = \rho_0 T_0 c_h^2 \tau_2 $.

The eigenstructure of the system \eqref{eqn.GPR} has been studied in \cite{JACKSON2020}. Here, we are only interested in an estimate of the maximum eigenvalues that can be chosen according to \cite{LGPR} as
\begin{equation}
	\label{eqn.lambda}
	\lambda = \sqrt{\frac{\gamma \, p}{\rho} + \frac{4}{3} c_s^2 + c_h^2}.
\end{equation}

\subsection{Overdetermined systems: extra conservation laws}

By construction, see \cite{PeshRom2014}, the model \eqref{eqn.GPR} is an overdetermined hyperbolic system, thus implying the satisfaction of additional (or extra) conservation laws. Firstly, we obtain total energy conservation from the HTC framework, then we focus on the derivation of the Geometric Conservation Law (GCL) that imposes a geometric constraint on the determinant of the distortion tensor $|\A|$.

\subsubsection{Total energy conservation law} 
By dot multiplying equations \eqref{eqn.conti}-\eqref{eqn.heatflux} with the associated thermodynamic variables $\p$, one obtains the total energy equation
\begin{equation}
\frac{\partial \mathcal{E} }{\partial t}
+\frac{ \partial \left( \mathcal{E} v_k + v_i \, p \, \delta_{ik} + \chi_k + h_k \right) }{\partial x_k}  
-\frac{\partial}{\partial x_m} \left( \epsilon \frac{\mathcal{E}}{\partial x_m} \right)
= 0, 
\label{eqn.energy} 	
\end{equation}
meaning that the following Gibbs relation is satisfied:
\begin{equation}
	\label{eqn.Gibbs}
	\begin{array}{rllllllllllllllll}
		1 \cdot d \mathcal{E}  &= 
		r &\cdot d\rho &+& 
		v_i &\cdot d (\rho v_i) &+&
		T &\cdot d \rhoS &+ &
		\alpha_{ik} &\cdot d A_{ik} &+ &
		\beta_k &\cdot dJ_k &:=
		\p &\cdot d\q \\
		1 \cdot \eqref{eqn.energy} &=
		r &\cdot \eqref{eqn.conti} &+& 
		v_i &\cdot \eqref{eqn.momentum} &+&
		T &\cdot \eqref{eqn.entropy} &+ &
		\alpha_{ik} &\cdot \eqref{eqn.deformation} &+& 
		\beta_{k} &\cdot \eqref{eqn.heatflux} 
	\end{array}.
\end{equation}
This also implies that the entropy production term $\Pi$ in \eqref{eqn.entropy} must be compatible with the parabolic dissipation terms 
\begin{equation}
	\partial_{\rhoS} \mathcal{E} \cdot \Pi + \p \cdot\partial_m \left( \epsilon \partial_m \q \right) = \partial_m  \left( \epsilon \partial_m \mathcal{E} \right),   
	\label{eqn.diss.comp} 
\end{equation}
and that the dot product of $\p$ with the algebraic relaxation source terms must vanish
\begin{equation}
	\p \cdot \mathbf{S}(\q) = 0.   
	\label{eqn.src.comp} 
\end{equation}

Although the rigorous formalism and derivation of HTC systems implies the use of the entropy as state variable, let us remark that the energy equation \eqref{eqn.energy} could be solved instead, and the associated entropy balance can be retrieved again from the Gibbs relation \eqref{eqn.Gibbs} as
\begin{equation}
	\label{eqn.Gibbs2}
	T \, d \rhoS  = - r \cdot d\rho - v_i \cdot d (\rho v_i) + 1 \cdot d \mathcal{E} - \alpha_{ik} \cdot d A_{ik} -\beta_k \cdot dJ_k,
\end{equation}
with a set of dual variables
\begin{equation}
	\label{eqn.dualvar_r}
	\r = \{r_j \} = \frac{1}{T} (-r,-v_i,1,-\alpha_{ik},-\beta_k)^\top.
\end{equation}
This implies the assumption of a physical entropy potential $\mathcal{S}$ such that $\r=\partial_\q \mathcal{S}$ with an associated positive semi-negative Hessian matrix $\partial^2_{q_i q_j} \mathcal{S} \leq 0$.

\subsubsection{Geometric Conservation Law (GCL)} 
The governing equations \eqref{eqn.GPR} also involve a geometric constraint on the determinant of the distortion tensor $\A$, which corresponds to the inverse deformation gradient for reversible processes in the material. To properly derive this geometric constraint, let us consider the Lagrange-Euler mapping between the Lagrangian domain $\Omega \subset  \mathbb{R}^3$ and the Eulerian domain $\omega(t)\subset  \mathbb{R}^3$ at time $t>0$, that deforms in time through the movement of the material.
\begin{figure}[!htbp]
	\begin{center}
		\includegraphics[clip,width=0.45\textwidth]{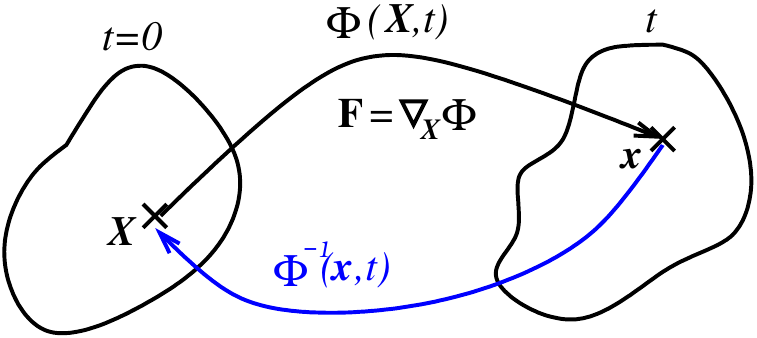} 
		\caption{Sketch of the Lagrangian-Eulerian transformation.}
		\label{fig:transformation}
	\end{center}
\end{figure}
At the aid of Figure \ref{fig:transformation}, let $\mathbf{X}=\{X_k\}$ and $\mathbf{x}=\{x_k\}$ represent the coordinate of any Lagrangian point in $\Omega$ and Eulerian point in $\omega(t)$, respectively. Then, the Lagrange-Euler mapping $\boldsymbol{\Phi}=\{\Phi_i\}$ is such that $\mathbf{x} = \boldsymbol{\Phi}( \mathbf{X},t) \in \omega(t)$, and the kinematic velocity of the material in the Eulerian frame is given by
\bea
\label{eqn.kinvel}
v_i( \mathbf{X}, t) = \frac{\partial \Phi_i( \mathbf{X},t )}{\partial t}.
\eea 
The deformation gradient tensor $\mathbf{F}=\{F_{ik}\}$ is nothing but the Jacobian matrix associated to the flow map $\Phi$ and verifies 
\bea \label{eq:F}
F_{ik}=\frac{\partial \Phi_i}{\partial X_k}.
\eea 
We assume that for all $t > 0$, the determinant of $\mathbf{F}$, called the Jacobian of the transformation, satisfies $\det(\mathbf{F}(\mathbf{X},t) ) := | \mathbf{F}(\mathbf{X},t)  | > 0$, so that the flow map is always invertible. The inverse of the transformation links the Eulerian coordinate to the Lagrangian one, i.e. $\mathbf{X}=\Phi^{-1}( \mathbf{x},t )$, and the distortion tensor in the mathematical model \eqref{eqn.GPR} is geometrically defined as $\mathbf{A}=\mathbf{F}^{-1}$ for reversible processes. The determinant of the deformation gradient represents the ratio of the Eulerian volume element to the Lagrangian volume element, that is
\begin{equation}
	\label{eqn.dv}
	dv = | \mathbf{F} | dV.
\end{equation}
Following \cite{Vilar1}, the mass conservation law with respect to the Lagrangian configuration is expressed for $t\geq 0$ by
\bea \label{eq:mass_cons}
\Frac{d}{dt} \Int_{\Omega} \rho (\mathbf{X},t) \,|\mathbf{F}(\mathbf{X},t)| \, dV = 0,
\eea
with the Lagrangian or material derivative given by
\begin{equation}
	\label{eqn.Lag_der}
	\Frac{d }{dt} = \Frac{\partial }{\partial t} + v_k \, \frac{\partial }{\partial x_k}.
\end{equation}
Since relation \eqref{eq:mass_cons} must hold for an arbitrary domain $\Omega$, it implies
\bea \label{eq:rhoJrho0}
\rho(\mathbf{X},t) |\mathbf{F}(\mathbf{X},t)| = \rho(\mathbf{X},0) \; \Longrightarrow \;
|\mathbf{F}| = \Frac{\rho_0}{\rho},
\eea
where we recall that $\rho_0=\rho(\mathbf{X},0)$ is the initial density of the material.
Consequently, thanks to the relationship $\mathbf{A}=\mathbf{F}^{-1}$, the determinant of the distortion matrix must obey the following constraint:
\bea \label{eq:detA}
|\mathbf{A}| = \Frac{\rho}{\rho_0}.
\eea 

This geometric constraint is extremely difficult to be respected at the discrete level, especially for Eulerian schemes. To the best knowledge of the authors, this has never been achieved so far on fixed grids. Therefore, our aim is to satisfy the constraint \eqref{eq:detA} by proposing a new approach, that requires the satisfaction of an extra conservation law for the quantity $|\A|$. Using the Lagrangian derivative \eqref{eqn.Lag_der} and neglecting viscous and source terms, the evolution equation \eqref{eqn.deformation} writes
\begin{equation}
	\label{eqn.deformationLag}
	\frac{d A_{ki}}{dt} + A_{kj} \, \frac{\partial v_j}{\partial x_i} = 0.
\end{equation}
Employing the Jacobi formula and the above relation, the time derivative of the determinant of the distortion tensor leads to
\begin{eqnarray}
	\label{eqn.detAJacobi}
	\frac{d |\A|}{dt} &=& \text{tr} \left( |\A| \, \mathbf{W} \, \frac{d \A}{dt} \right), \qquad \mathbf{W}= \A^{-1}, \nn \\
	&=& - |\A| \, W_{ik} \, A_{kj} \, \frac{\partial v_j}{\partial x_i} \nn \\
	&=& - |\A| \, \delta_{ij} \, \frac{\partial v_j}{\partial x_i}.
\end{eqnarray} 
By replacing the material derivative on the left hand side of \eqref{eqn.detAJacobi} with its Eulerian counterpart according to \eqref{eqn.Lag_der}, the Geometric Conservation Law is obtained as an extra conservation law satisfied by the governing equations \eqref{eqn.GPR}, which is explicitly given by
\begin{equation}
	\label{eqn.detApde}
	\frac{\partial |\A|}{\partial t} + \frac{\partial }{\partial x_k} (|\A| v_k)= 0.
\end{equation}

Therefore, satisfying the GCL \eqref{eqn.detApde} implies that the constraint \eqref{eq:detA} is also respected. To mimic the HTC approach, let us introduce a new set of pseudo-\textit{thermodynamic variables} $\w=\{w_{ik}\}$ that are dual with respect to a pseudo-\textit{potential} given by $|\A|$, thus obtaining $\w:=\partial_{\A} |\A|$. Then, by construction, one can verify that a pseudo-\textit{Gibbs relation} is satisfied, that is
\begin{equation}
	\label{eqn.GibbsA}
	\begin{array}{rlll}
  d(|\A|) &= w_{ik} &\cdot d(A_{ik}) & := \w \cdot d\A  \\
  \eqref{eqn.detApde} &= w_{ik} & \cdot \eqref{eqn.deformation} &
  \end{array} .
\end{equation}
More precisely, the source term of the distortion tensor equation \eqref{eqn.deformation}, referred to as $\mathbf{S}_{\A}$, has been designed in \cite{PeshRom2014} not to affect the mass conservation equation. Indeed, it is proportional to $\partial_{\A} \mathcal{E}$, namely 
\begin{equation}
	\mathbf{S}_{\A} = -\frac{\partial_{\A} \mathcal{E}}{\theta_1(\tau_1)} = -\frac{\boldsymbol{\alpha}}{\theta_1(\tau_1)} = - \frac{3}{\tau_1} |\A|^{\frac{5}{3}} \, \A \, \mathring{\mathbf{G}}.
\end{equation} 
On the other hand, the dual variables $\w$ are given by
\begin{equation}
	\w = |\A| \, \A^{-\top}.
\end{equation}
Therefore, the contraction $\w : \mathbf{S}_{\A}=\text{tr}(\w^\top \, \mathbf{S}_{\A})$ yields
\begin{equation}
	\label{eqn.zeroSrcA}
	\w : \mathbf{S}_{\A} = - \frac{3}{\tau_1} |\A|^{\frac{8}{3}} \text{tr} ( \A^{-1} \, \A \, \mathring{\mathbf{G}} ) = - \frac{3}{\tau_1} |\A|^{\frac{8}{3}} \text{tr} ( \mathring{\mathbf{G}} ) = 0,
\end{equation}
since $\mathring{\mathbf{G}}$ is the trace-free part of the metric tensor $\mathbf{G}=\mathbf{A}^\top \, \mathbf{A}$. The details concerning the derivation of the GCL in terms of the dual variables $\w$ can be found in \ref{app.gcl}.

\subsection{Reduced compatible model} \label{ssec:reduced_GPR}
The previous considerations incline us to consider a reduced model consisting of the following equations:
\begin{subequations}\label{eqn.GPR_red}
	\begin{align}
		\frac{\partial \rho v_i}{\partial t}
		+\frac{\partial \left( \rho v_i v_k + p \, \delta_{ik}   
			+  \sigma_{ik} + \phi_{ik} \right)}{\partial x_k}  
		-\frac{\partial}{\partial x_m}  \left( \epsilon \frac{\partial \rho v_i}{\partial x_m} \right) 
		&= 0,  
		\label{eqn.momentum_red} \\ 
		\frac{\partial \mathcal{E} }{\partial t}
		+\frac{ \partial \left( \mathcal{E} v_k + v_i \, p \, \delta_{ik} + \chi_k + h_k \right) }{\partial x_k}  
		-\frac{\partial}{\partial x_m} \left( \epsilon \frac{\mathcal{E}}{\partial x_m} \right)
		&= 0,   
		\label{eqn.energy_red} \\	
		\frac{\partial A_{i k}}{\partial t}
		+\frac{\partial (A_{im} v_m)}{\partial x_k} 
		+ v_m \left(\frac{\partial A_{ik}}{\partial x_m}-\frac{\partial A_{im}}{\partial x_k}\right)   -\frac{\partial}{\partial x_m}  \left( \epsilon \frac{\partial A_{ik}}{\partial x_m} \right) 
		&=  -\dfrac{ \alpha_{ik} }{\theta_1(\tau_1)} ,\label{eqn.deformation_red} 
		\\
		\frac{\partial J_k}{\partial t}
		+\frac{\partial \left( J_m v_m + T \right)}{\partial x_k} 
		+ v_m \left(\frac{\partial J_{k}}{\partial x_m} -\frac{\partial J_{m}}{\partial x_k}\right)  
		-\frac{\partial}{\partial x_m} \left( \epsilon \frac{\partial J_{k}}{\partial x_m} \right)  
		&= -\dfrac{\beta_k}{\theta_2(\tau_2)}. 
		\label{eqn.heatflux_red} 
	\end{align}
\end{subequations}
This system satisfies the entropy inequality \eqref{eqn.entropy} and the Geometric Conservation Law \eqref{eqn.detApde}. We underline that no evolution equation for the mass density is embedded in the model, since the material density can be easily computed from the determinant constraint \eqref{eq:detA} thanks to the GCL compatibility, that is $\rho = \rho_0 |\mathbf{A}|$. Likewise, the entropy balance is also satisfied by the reduced model \eqref{eqn.GPR_red} which is compliant with the Gibbs relation \eqref{eqn.Gibbs2}.

Here we consider the state variables $\u=( \rho v_i, \mathcal{E}, A_{ik}, J_k)$, and the governing equations can be written in a compact matrix-vector formulation as
\begin{equation}
	\frac{\partial \u}{\partial t}  + \frac{\partial \mathbf{f}_k(\u)}{\partial x_k}  +  \mathbf{B}_k(\u) \frac{\partial \u}{\partial x_k} - \frac{\partial }{\partial m} \left( \epsilon \frac{\partial \u}{\partial m}  \right) =  \mathbf{S}(\u),
	\label{eqn.pde:red} 
\end{equation}
where $\mathbf{f}_k(\q)$ is the nonlinear flux tensor and $\mathbf{B}_k(\u) \partial_k \u$ contains the non-conservative part of the system in block-matrix notation for $\A$ and $\J$. The algebraic sources are gathered in the term $\mathbf{S}(\u)$, while the regularizing viscous terms are given by $\partial_m \left( \epsilon \partial_m \u \right)$. 

The model \eqref{eqn.GPR_red} is solved with a finite volume method on general unstructured meshes that is proven to preserve both the geometric and the thermodynamic compatibility. Indeed, the novel numerical method only solves the reduced model \eqref{eqn.GPR_red} because it is compliant with \eqref{eqn.entropy} and \eqref{eqn.detApde}. All the details of the numerical scheme are provided in the next section.

\section{Numerical scheme} \label{sec.numscheme}

\subsection{Semi-discrete finite volume scheme on unstructured meshes}
\begin{figure}[!htbp]
	\begin{center}
		\begin{tabular}{cc} 
			\includegraphics[trim=0 0 0 0,clip,width=0.47\textwidth]{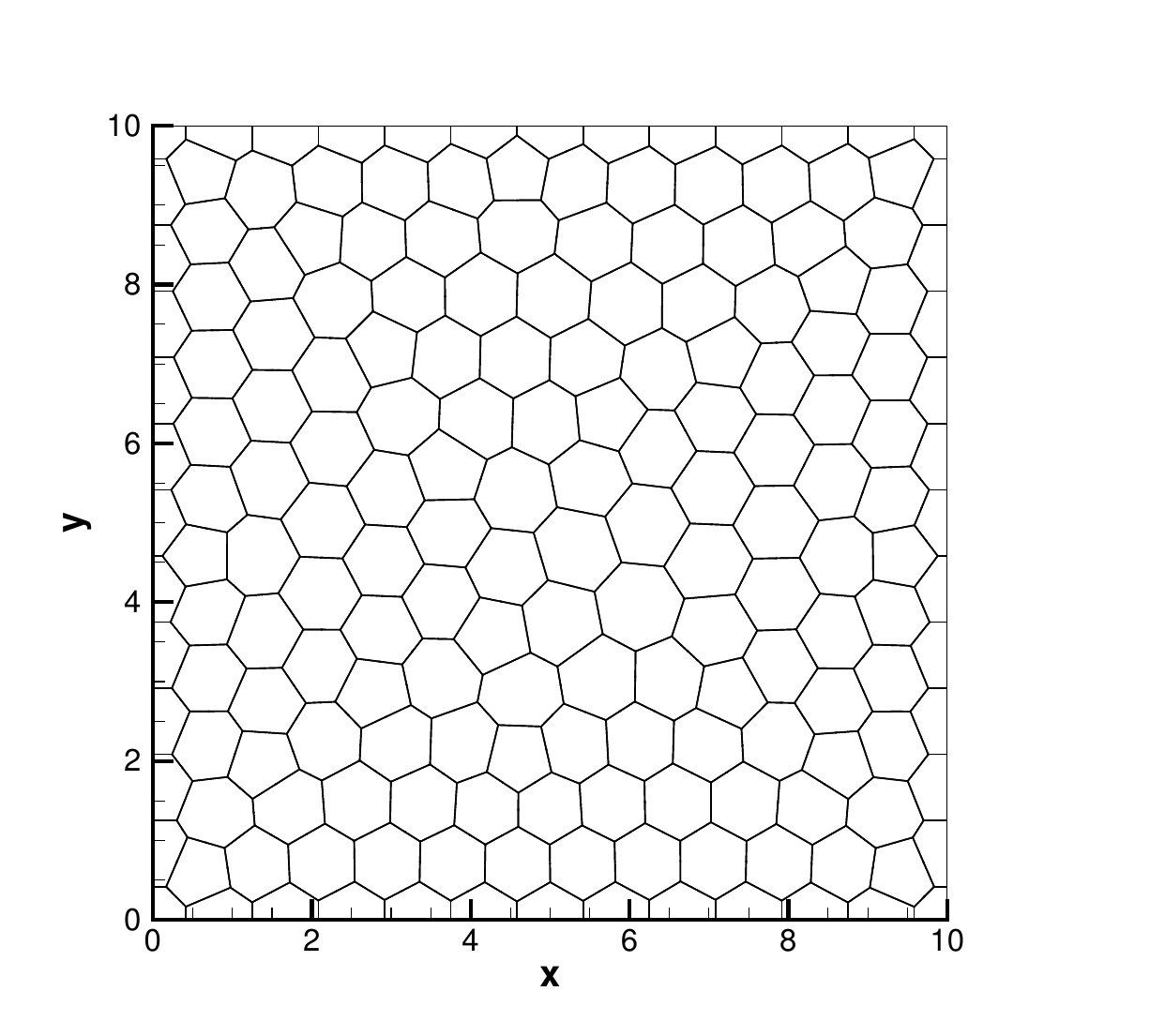} & 
			\includegraphics[trim=0 0 0 0,clip,width=0.47\textwidth]{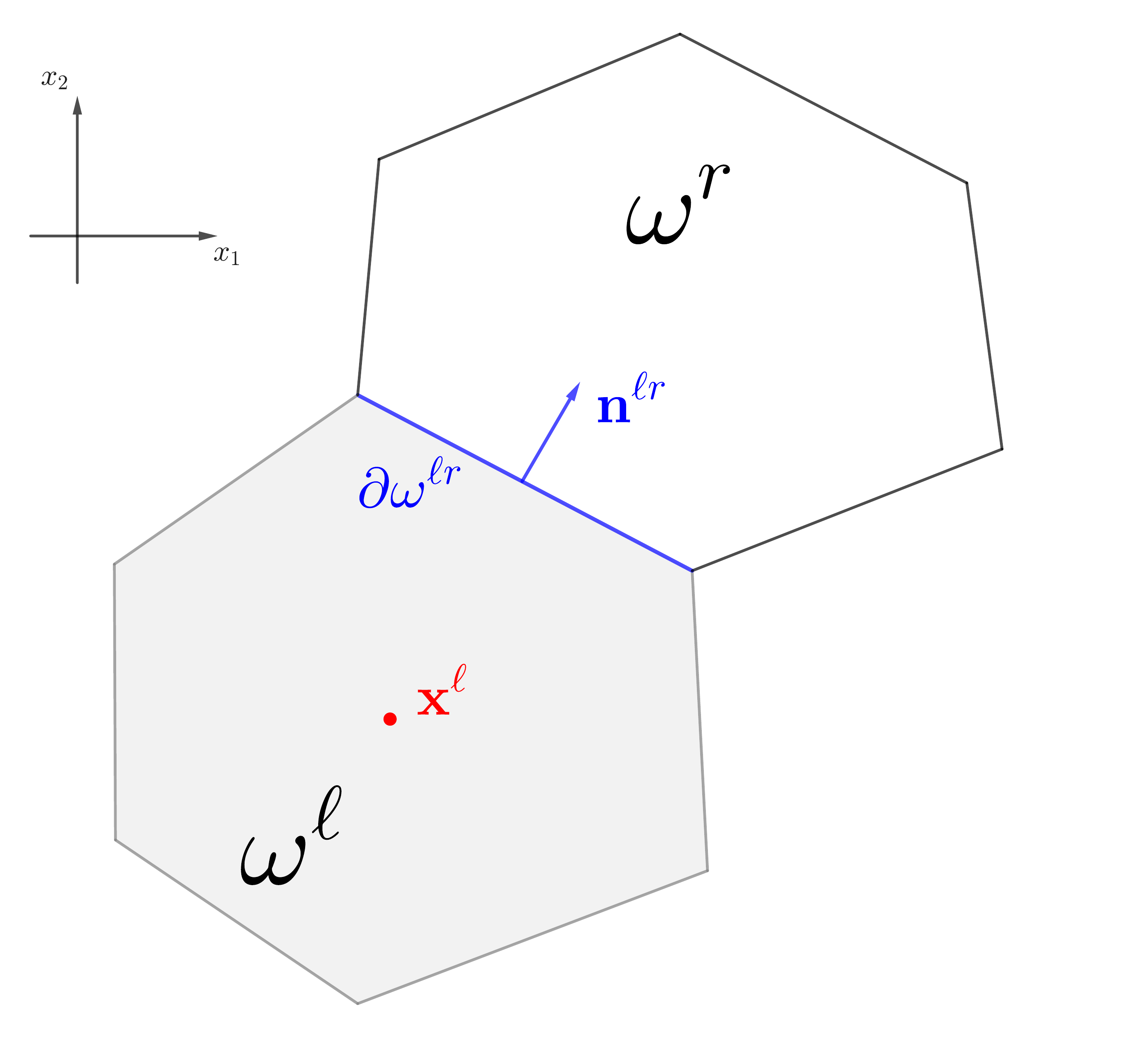} \\
		\end{tabular} 
		\caption{Left: example of unstructured Voronoi mesh. Right: notation used for cell $\omega^\ell$ and one direct neighbor cell $\omega^\err$.}
		\label{fig.notation_mesh}
	\end{center}
\end{figure}

To ease the notation and the readability, subscripts are used for tensor indices while superscripts denote the spatial discretization index. The two-dimensional computational domain $\Omega \in \mathds{R}^2$ is discretized with a total number $N^\ell$ of non-overlapping unstructured polygonal control volumes $\omega^{\ell}$ with border $\partial \omega^\ell$ and barycenter coordinates $\xx^\ell$. We underline that Voronoi meshes can also be employed as well as any other general polygonal elements. The surface of the element is denoted with $|\omega^\ell|$, whereas $|\partial \omega^\ell|$ refers to the cell perimeter. The set of neighbors of cell $\omega^\ell$ is labeled with $\mathcal{N}^\ell$, and $\partial \omega^{\ell \err}$ is the common edge shared by two adjacent elements $\omega^\ell$ and $\omega^\err$ with outward pointing unit normal vector $\nv^{\ell\err}$. Figure \ref{fig.notation_mesh} shows an example of an unstructured polygonal mesh and a sketch of the adopted notation.

A semi-discrete first order finite volume scheme for the governing equations \eqref{eqn.GPR_red} writes
\begin{equation}
	\label{eqn.semidisc_fv}
	\frac{\partial \u^\ell}{\partial t} = - \sum \limits_{\err \in \mathcal{N}^\ell} \frac{|\partial \omega^{\ell \err}|}{|\omega^\ell|} \left( \mathcal{F}(\u^\ell,\u^\err)^{\ell \err} + \mathcal{D}(\u^\ell,\u^\err)^{\ell \err} + \mathcal{G}(\u^\ell,\u^\err)^{\ell \err} \right) \cdot \nv^{\ell \err} + \mathbf{S}(\u^\ell).
\end{equation}
The numerical flux is given by a central approximation, that is
\begin{equation}
	\label{eqn.flux}
	\mathcal{F}(\u^\ell,\u^\err)^{\ell \err} \cdot \nv^{\ell \err} = \halb \left( \mathbf{f}_k^\ell + \mathbf{f}_k^\err \right) \, n_k^{\ell \err},
\end{equation}
and the following discretization is chosen for the non-conservative terms
\begin{equation}
	\label{eqn.ncprod}
	\mathcal{D}(\u^\ell,\u^\err)^{\ell \err} \cdot \nv^{\ell \err} = \halb \mathbf{B}_k(\bar{\u}^{\ell \err}) \, n_k^{\ell \err} \, ( \u^\err - \u^\ell ), \qquad \bar{\u}^{\ell \err} = \halb ( \u^\err + \u^\ell ).
\end{equation}
The dissipative numerical flux is computed by a Rusanov-type scheme with
\begin{equation}
	\label{eqn.diss}
	\mathcal{G}(\u^\ell,\u^\err)^{\ell \err} \cdot \nv^{\ell \err} = -\epsilon^{\ell \err} \, \|\xx^\err-\xx^\ell \| \, \frac{\u^\err - \u^\ell}{\|\xx^\err-\xx^\ell\|} = -\epsilon^{\ell \err} \, (\u^\err - \u^\ell ), \qquad \epsilon^{\ell \err} = \halb \, \max \left( |\lambda(\u^\ell)|, \, |\lambda(\u^\err)| \right),
\end{equation}
where we use the estimate of the maximum eigenvalue of the system according to \eqref{eqn.lambda} for the definition of the positive coefficient $\epsilon^{\ell \err}$. As it stands, the finite volume scheme \eqref{eqn.semidisc_fv} is not structure preserving: it is neither compliant with the extra conservation laws \eqref{eqn.entropy} nor with \eqref{eqn.detApde}. Therefore some \textit{ad hoc} modifications have to be designed in order to make the scheme \textit{structure preserving}.

\subsection{Geometrically compatible finite volume scheme}
We start by designing a modification of the finite volume scheme \eqref{eqn.semidisc_fv} such that geometric compatibility is ensured by satisfying the GCL equation \eqref{eqn.detApde} as an extra conservation law of the mathematical model. For the moment, we neglect the source terms and the dissipation fluxes, i.e we assume $\mathcal{G}(\u^\ell,\u^\err)^{\ell \err} \cdot \nv^{\ell \err} = \mathbf{0}$ and $\mathbf{S}(\u)=\mathbf{0}$. We rely on a very general method firstly proposed for achieving thermodynamic compatibility in \cite{Abgrall2018}, and more recently extended to hyperbolic systems of the type \eqref{eqn.GPR} to recover energy conservation from the direct discretization of the entropy inequality \cite{DumbserHTC,HTCAbgrall,HTCMHD}. Here, we apply this strategy for the first time to preserve a different structural property rather than thermodynamics.

Let us interpret the determinant of the distortion matrix $|\A|$ as a \textit{thermodynamic potential}, and the associated dual variables $\w=\{w_{ik}\} = \partial_{\A} |\A|$ as a set of \textit{thermodynamic variables}. Moreover, let $\mathcal{F}_{\A}(\u^\ell,\u^\err) \cdot \nv^{\ell \err}$ denote the central fluxes related to the distortion tensor equation \eqref{eqn.deformation_red} according to \eqref{eqn.flux}. Likewise, $\mathbf{f}_{\A,m}$ represents the physical flux of equation \eqref{eqn.deformation_red} and $\mathcal{D}_{\A}(\u^\ell,\u^\err) \cdot \nv^{\ell \err}$ are the corresponding fluctuations of the non-conservative terms restricted to \eqref{eqn.deformation_red}. According to \cite{Abgrall2018}, these fluxes are modified by a correction factor $\alpha_\A^{\ell \err}$, which is defined at the cell interface, hence obtaining the modified fluxes
\begin{equation}
 \label{eqn.modifiedfluxA}	
 \tilde{\mathcal{F}}_{\A}(\u^\ell,\u^\err) \cdot \nv^{\ell \err} = \mathcal{F}_{\A}(\u^\ell,\u^\err) \cdot \nv^{\ell \err} - \alpha_\A^{\ell \err} \, ( \w^\err - \w^\ell ) \cdot \nv^{\ell \err} = \halb \left( \mathbf{f}_{\A,k}^\ell + \mathbf{f}_{\A,k}^\err \right) n_k^{\ell \err} - \alpha_\A^{\ell \err} \, ( \w^\err - \w^\ell ) \cdot \nv^{\ell \err}.
\end{equation} 
The scalar correction factor $\alpha_\A^{\ell \err}$ has no sign, and it can add or subtract the total amount of the jump in the dual variables which is needed to reach geometric compatibility with the GCL \eqref{eqn.detApde}. To determine $\alpha_\A^{\ell \err}$, the conservation principle is invoked. Indeed, across each cell boundary $\partial \omega^{\ell \err}$, a consistent condition implies that the sum of the fluctuations must balance the sum of the fluxes which have to be preserved. This sum, namely $|\A|v_k$ in \eqref{eqn.detApde}, must be recovered as the dot product of equation \eqref{eqn.deformation_red} with the dual variables $\w^\ell$, that is
\begin{equation}
	\label{eqn.geom_eqn}
	\begin{aligned}
	\w^\ell \cdot \left( \tilde{\mathcal{F}}_{\A}(\u^\ell,\u^\err) \cdot \nv^{\ell \err} - \mathbf{f}_{\A,k}^\ell \cdot n_k^{\ell \err}\right) + \w^\err \cdot \left( \mathbf{f}_{\A,k}^\err \cdot n_k^{\ell \err} - \tilde{\mathcal{F}}_{\A}(\u^\ell,\u^\err) \cdot \nv^{\ell \err} \right) & +& \\ \w^\ell \cdot \mathcal{D}_{\A}(\u^\ell,\u^\err) \cdot \nv^{\ell \err} + \w^\err \cdot \mathcal{D}_{\A}(\u^\err,\u^\ell) \cdot \nv^{\err \ell}  & =& \left((|\A| v_k)^\err - (|\A| v_k)^\ell \right) \, n_k^{\ell \err}.
	\end{aligned}
\end{equation}
By inserting the definition of the modified fluxes \eqref{eqn.modifiedfluxA} in the condition \eqref{eqn.geom_eqn}, we obtain the sought correction factor defined at the cell interface $\partial \omega^{\ell \err}$:
\begin{equation}
	\label{eqn.alphaA}
	\begin{aligned}
	\alpha_\A^{\ell \err} & = \frac{\left((|\A| v_k)^\err - (|\A| v_k)^\ell \right) \, n_k^{\ell \err} + \left( \mathcal{F}_{\A}(\u^\ell,\u^\err) \cdot \nv^{\ell \err} \right) \cdot (\w^\err-\w^\ell) - \left( \w^\err \cdot \mathbf{f}_{\A,k}^\err - \w^\ell \cdot \mathbf{f}_{\A,k}^\ell\right)\, n_k^{\ell \err}}{(\w^\err - \w^\ell)^2} \\
	& - \frac{( \w^\err + \w^\ell ) \cdot \mathcal{D}_{\A}(\u^\ell,\u^\err) \cdot \nv^{\ell \err}}{(\w^\err - \w^\ell)^2}.
	\end{aligned}
\end{equation}
Obviously, if $\w^\ell=\w^\err$ then no correction occurs and we simply set $\alpha_\A^{\ell \err}=0$.

\paragraph{Geometric compatibility with dissipation fluxes and source terms} Even with smooth initial data, the solution of the PDE system \eqref{eqn.GPR_red} can exhibit shocks and other discontinuities, which require a stabilization of the numerical scheme that is carried out relying on parabolic vanishing viscosity terms. Also in this case, the compatibility with the geometric extra conservation law \eqref{eqn.detApde} must be respected. To that aim, let us add to the compatible fluxes \eqref{eqn.modifiedfluxA} also the dissipative fluxes $\mathcal{G}_{\A}(\u^\ell,\u^\err)^{\ell \err} \cdot \nv^{\ell \err}$ as well as the source terms $\mathbf{S}_{\A}(\u^\ell)=-\frac{\boldsymbol{\alpha}^\ell}{\theta (\tau_1)}$, so that the semi-discrete evolution equation for $\A$ \eqref{eqn.deformation_red} now becomes
\begin{equation}
	\label{eqn.Avisc_source}
	\frac{\partial \A^\ell}{\partial t} + \sum \limits_{\err \in \mathcal{N}^\ell} \frac{|\partial \omega^{\ell \err}|}{|\omega^\ell|} \left( \tilde{\mathcal{F}}_{\A}(\u^\ell,\u^\err)^{\ell \err} + \mathcal{D}_{\A}(\u^\ell,\u^\err)^{\ell \err} \right) \cdot \nv^{\ell \err} = -\sum \limits_{\err \in \mathcal{N}^\ell} \mathcal{G}_{\A}(\u^\ell,\u^\err) \cdot \nv^{\ell \err} + \mathbf{S}_{\A}(\u^\ell).
\end{equation}
The part on the left hand side of the above equation is already compatible with the GCL thanks to the modified fluxes \eqref{eqn.modifiedfluxA} with the scalar correction factor given by \eqref{eqn.alphaA}. Therefore, we focus on the compatibility of the right hand side of \eqref{eqn.Avisc_source}. Recalling the definition \eqref{eqn.diss} and following \cite{DumbserHTC}, after multiplication by the dual variables $\w^\ell$, for the viscous terms we obtain
\begin{equation}
	\label{eqn.dissA1}
	\begin{aligned}
		\w^\ell \cdot \mathcal{G}_{\A}(\u^\ell,\u^\err) \cdot \nv^{\ell \err} &= \halb \left( \w^\ell \cdot \mathcal{G}_{\A}(\u^\ell,\u^\err) \cdot \nv^{\ell \err} + \w^\err \cdot \mathcal{G}_{\A}(\u^\ell,\u^\err) \cdot \nv^{\ell \err} + \w^\ell \cdot \mathcal{G}_{\A}(\u^\ell,\u^\err) \cdot \nv^{\ell \err} - \w^\err \cdot \mathcal{G}_{\A}(\u^\ell,\u^\err) \cdot \nv^{\ell \err} \right) \\
		&= \underbrace{\halb ( \w^\err - \w^\ell ) \cdot \epsilon^{\ell \err} ( \A^\err - \A^\ell )}_{\mathcal{G}_1} - \underbrace{\halb ( \w^\err + \w^\ell ) \cdot \epsilon^{\ell \err} ( \A^\err - \A^\ell )}_{\mathcal{G}_2}.
	\end{aligned}
\end{equation}
The second term $\mathcal{G}_2$ is the approximation of the jump term related to the numerical dissipation in the GCL \eqref{eqn.detApde}, that is 
\begin{equation}
	\label{eqn.dissA2}
	-\halb ( \w^\err + \w^\ell ) \cdot \epsilon^{\ell \err} ( \A^\err - \A^\ell ) \approx -\epsilon^{\ell \err} ( |\A|^\err - |\A|^\ell ).
\end{equation}
Indeed,  applying path integration in the state variables $\A$, the following relation holds true by construction:
\begin{equation}
	\label{eqn.Roe_p}
	\int \limits_{\A^\ell}^{\A^\err} \w \cdot d\A = \int \limits_{\A^\ell}^{\A^\err} \partial_{\A} |\A| \cdot d\A = |\A|^\err - |\A|^\ell,
\end{equation} 
and the term $\halb ( \w^\err + \w^\ell ) ( \A^\err - \A^\ell )$ in \eqref{eqn.dissA2} can be seen as a numerical approximation of the path integral in \eqref{eqn.Roe_p} using a trapezoidal rule. Therefore, we still remain with an additional contribution given by the first term $\mathcal{G}_1$ in \eqref{eqn.dissA1}. To control its production of numerical dissipation, we reformulate the jump in the dual variables $\w$ as a jump in the state variables $\A$ through the Hessian matrix $\partial^2_{\A \A} |\A|^{\ell \err}$ which verifies the Roe property 
\begin{equation}
	\label{eqn.Hessian}
	\partial^2_{\A \A} |\A|^{\ell \err} \cdot ( \A^\err - \A^\ell ) = \w^\err - \w^\ell.
\end{equation}
The Hessian matrix at the cell interface is computed as
\begin{equation}
	\label{eqn.HessianA_face}
	\partial^2_{\A \A} |\A|^{\ell \err} = \int \limits_0^1  \partial^2_{\A \A} |\A|(\boldsymbol{\xi}(s)) \, ds, \qquad  \boldsymbol{\xi}(s) = \A^\ell + s \, ( \A^\err - \A^\ell ), \quad 0 \leq s \leq 1,
\end{equation}
where $\partial^2_{\A \A} |\A|$ explicitly writes  
\begin{equation}
	\label{eqn.HessianA}
	\partial^2_{\A \A} |\A| = \left[\begin{array}{ccccccccc}
		0 & 0 & 0 & 0 & \phantom{-}A_{33} & -A_{32} & 0 & -A_{23} & \phantom{-}A_{22} 
		\\
		0 & 0 & 0 & -A_{33} & 0 & \phantom{-}A_{31} & \phantom{-}A_{23} & 0 & -A_{21} 
		\\
		0 & 0 & 0 & \phantom{-}A_{32} & -A_{31} & 0 & -A_{22} & \phantom{-}A_{21} & 0 
		\\
		0 & -A_{33} & \phantom{-}A_{32} & 0 & 0 & 0 & 0 & \phantom{-}A_{13} & -A_{12} 
		\\
		\phantom{-}A_{33} & 0 & -A_{31} & 0 & 0 & 0 & -A_{13} & 0 & \phantom{-}A_{11} 
		\\
		-A_{32} & \phantom{-}A_{31} & 0 & 0 & 0 & 0 & A_{12} & -A_{11} & 0 
		\\
		0 & \phantom{-}A_{23} & -A_{22} & 0 & -A_{13} & \phantom{-}A_{12} & 0 & 0 & 0 
		\\
		-A_{23} & 0 & \phantom{-}A_{21} & \phantom{-}A_{13} & 0 & -A_{11} & 0 & 0 & 0 
		\\
		\phantom{-}A_{22} & -A_{21} & 0 & -A_{12} & \phantom{-}A_{11} & 0 & 0 & 0 & 0 
	\end{array}\right].
\end{equation}
These contributions, which come from all the faces $\mathcal{N}^\ell$ of the cell, must vanish in order to obtain compatibility with the GCL equation \eqref{eqn.detApde}. Consequently, a production term $\Pi_\A^\ell$ is introduced to balance these terms with opposite sign, that is given by
\begin{equation}
	\label{eqn.prodA}
	\Pi_\A^\ell = \sum \limits_{\err \in \mathcal{N}^\ell} \frac{|\partial \omega^{\ell \err}|}{|\omega^\ell|} \, \halb \epsilon^{\ell \err} ( \A^\err - \A^\ell ) \cdot  \partial^2_{\A \A} |\A|^{\ell \err} \, ( \A^\err - \A^\ell ),
\end{equation}
where only jumps in the state variables appear because of the use of the Roe-type Hessian matrix \eqref{eqn.HessianA_face}. The production term $\Pi_\A^\ell$ is a scalar, that now needs to be distributed among all the components $A_{ik}^\ell$ of the distortion tensor, hence obtaining new contributions $P_{ik}^\ell$. Here, we adopt a rescaling with respect to the trace of the distortion tensor, as proposed in \cite{SWETurbulence} for the redistribution of a production term associated to the Reynolds stress tensor, thus we define $P_{ik}^\ell$ as
\begin{equation}
	P_{ik}^\ell = \Pi_\A^\ell \, \frac{w_{ik}^\ell}{\text{tr}(\w^\ell \, \w^{\ell,\top})},
\end{equation}
with the positive trace $\text{tr}(\w^\ell \, \w^{\ell,\top}) = w_{ik}^\ell w_{ik}^\ell \geq 0$. 

At last, it remains to verify the compatibility with the source terms. Multiplication of $\mathbf{S}_{\A}(\u^\ell)$ by the dual variables yields
\begin{equation}
	\label{eqn.sourceAcomp}
	-\w^\ell \cdot \frac{\boldsymbol{\alpha}^\ell}{\theta (\tau_1)} = -w_{ik}^\ell \cdot \frac{\alpha_{ik}^\ell}{\theta(\tau_1)} = 0,
\end{equation}
therefore the compatibility is proven by construction of the dual variables $\w$, as demonstrated at the continuous level by \eqref{eqn.zeroSrcA}. All the related details are reported in \ref{app.gcl}.

\begin{theorem}[Geometric compatibility] \label{th1}
	The semi-discrete finite volume scheme for the equation of the distortion tensor $\A$ given by 
	\begin{equation}
		\label{eqn.dAdt}
		\frac{\partial \A^\ell}{\partial t}  + \sum \limits_{\err \in \mathcal{N}^\ell} \frac{|\partial \omega^{\ell \err}|}{|\omega^\ell|} \left( \tilde{\mathcal{F}}_{\A}(\u^\ell,\u^\err) + \mathcal{D}_{\A}(\u^\ell,\u^\err) + \mathcal{G}_{\A}(\u^\ell,\u^\err) \right) \cdot \nv^{\ell \err} = - \frac{\boldsymbol{\alpha}^\ell}{\theta(\tau_1)} + \Pi_\A^\ell \, \frac{\w^\ell}{\text{tr}(\w^\ell \, \w^{\ell,\top})},
	\end{equation} 
	with the geometrically compatible fluxes \eqref{eqn.modifiedfluxA}, the non-conservative products \eqref{eqn.ncprod}, the dissipation terms \eqref{eqn.diss} and the production term \eqref{eqn.prodA}, satisfies the extra conservation law \eqref{eqn.detApde} with the following conservative semi-discrete scheme:
	\begin{equation}
		\label{eqn.ddetAdt}
		\frac{\partial |\A|^\ell}{\partial t} + \sum \limits_{\err \in \mathcal{N}^\ell} \frac{|\partial \omega^{\ell \err}|}{|\omega^\ell|} \halb \left(F_{|\A|}^\ell + F_{|\A|}^\err \right) \cdot \nv^{\ell \err} = 0.
	\end{equation}
\end{theorem}

\begin{proof}
	To ease the notation, let us use the abbreviations
	\begin{equation}
		\label{eqn.th1_1}
	\tilde{\mathcal{F}}_{\A}^{\ell \err}:=\tilde{\mathcal{F}}_{\A}^{\ell \err}(\u^\ell,\u^\err), \qquad \mathcal{D}_{\A}^{\ell \err}:=\mathcal{D}_{\A}(\u^\ell,\u^\err), \qquad \mathcal{G}_{\A}^{\ell \err}:=\mathcal{G}_{\A}(\u^\ell,\u^\err).
	\end{equation} 
	Moreover, let us recall that the discrete Gauss theorem over a close surface yields the relation
	\begin{equation}
		\label{eqn.Gaussdiv}
		\sum \limits_{\err \in \mathcal{N}^\ell} {|\partial \omega^{\ell \err}|} \, \nv^{\ell \err} = \mathbf{0}. 
	\end{equation}
    Using the above shorthand notation and dot multiplying the distortion equation \eqref{eqn.dAdt} by the dual variables $\w^\ell$, we obtain
	\begin{equation}
		\label{eqn.th1_2}
		\w^\ell \cdot \frac{\partial \A^\ell}{\partial t} +  \sum \limits_{\err \in \mathcal{N}^\ell} \frac{|\partial \omega^{\ell \err}|}{|\omega^\ell|} \, \w^\ell \cdot \left( \tilde{\mathcal{F}}_{\A}^{\ell \err} + \mathcal{D}_{\A}^{\ell \err} + \mathcal{G}_{\A}^{\ell \err}  \right) \cdot \nv^{\ell \err} = \w^\ell \cdot \left( - \frac{\boldsymbol{\alpha}^\ell}{\theta(\tau_1)} + \Pi_\A^\ell \, \frac{\w^\ell}{\text{tr}(\w^\ell \, \w^{\ell,\top})} \right).
	\end{equation}
    On the right hand side, we have
    \begin{equation}
    	\label{eqn.th1_3}
    	\w^\ell \cdot \left( - \frac{\boldsymbol{\alpha}^\ell}{\theta(\tau_1)} + \Pi_\A^\ell \, \frac{\w^\ell}{\text{tr}(\w^\ell \, \w^{\ell,\top})} \right) = 0 + \Pi_\A^\ell,
    \end{equation}
    where the first term vanishes thanks to the compatibility condition \eqref{eqn.sourceAcomp} (see \ref{app.gcl}) and the second term verifies by construction the relation 
    \begin{equation}
    	\label{eqn.th1_4}
    	\w^\ell \cdot \Pi_\A^\ell \, \frac{\w^\ell}{\text{tr}(\w^\ell \, \w^{\ell,\top})} = w^\ell_{ik} \, \Pi_\A^\ell \, \frac{w^\ell_{ik}}{w^\ell_{ik} \, w^\ell_{ik}} = \Pi_\A^\ell.
    \end{equation}
    On the left hand side of \eqref{eqn.th1_2}, we add and subtract the terms $\halb \, \w^\err \cdot \tilde{\mathcal{F}}_{\A}^{\ell \err} \cdot \nv^{\ell \err}$, $\halb \, \w^\err \cdot {\mathcal{D}}_{\A}^{\err \ell} \cdot \nv^{\err \ell}$ and $\halb \, \w^\err \cdot {\mathcal{G}}_{\A}^{\ell \err} \cdot \nv^{\ell \err}$, hence obtaining
    \begin{equation}
    	\label{eqn.th1_5}
    \begin{aligned}
    	\frac{\partial |\A|^\ell}{\partial t} & + \halb \, \sum \limits_{\err \in \mathcal{N}^\ell} \frac{|\partial \omega^{\ell \err}|}{|\omega^\ell|} \, \left( ( \w^\ell + \w^\err ) \cdot \tilde{\mathcal{F}}_{\A}^{\ell \err} \cdot \nv^{\ell \err} + ( \w^\ell - \w^\err ) \cdot \tilde{\mathcal{F}}_{\A}^{\ell \err} \cdot \nv^{\ell \err} \right) \\
    	& + \halb \, \sum \limits_{\err \in \mathcal{N}^\ell} \frac{|\partial \omega^{\ell \err}|}{|\omega^\ell|} \,  \left( \w^\ell \cdot \mathcal{D}_{\A}^{\ell \err} \cdot \nv^{\ell \err} + \, \w^\err \cdot \mathcal{D}_{\A}^{\err \ell} \cdot \nv^{\err \ell} + \w^\ell \cdot \mathcal{D}_{\A}^{\ell \err} \cdot \nv^{\ell \err} - \, \w^\err \cdot \mathcal{D}_{\A}^{\err \ell} \cdot \nv^{\err \ell} \right) \\
    	& + \halb \, \sum \limits_{\err \in \mathcal{N}^\ell} \frac{|\partial \omega^{\ell \err}|}{|\omega^\ell|} \,  \left(  ( \w^\ell + \w^\err ) \cdot {\mathcal{G}}_{\A}^{\ell \err} \cdot \nv^{\ell \err} + ( \w^\ell - \w^\err ) \cdot {\mathcal{G}}_{\A}^{\ell \err} \cdot \nv^{\ell \err} \right) \\
    	& = \Pi_\A^\ell.
    \end{aligned}	
    \end{equation}
    Due to the continuity of the computational mesh, it holds that $\nv^{\ell \err}=-\nv^{\err \ell}$. Furthermore, the term $( \w^\ell - \w^\err ) \cdot \tilde{\mathcal{F}}_{\A}^{\ell \err} \cdot \nv^{\ell \err}$ can be rewritten by means of the compatibility condition \eqref{eqn.geom_eqn}, which leads to
    \begin{equation}
    	\label{eqn.th1_6}
    	\begin{aligned}
    		\frac{\partial |\A|^\ell}{\partial t} & + \halb \, \sum \limits_{\err \in \mathcal{N}^\ell} \frac{|\partial \omega^{\ell \err}|}{|\omega^\ell|} \, ( \w^\ell + \w^\err ) \cdot \tilde{\mathcal{F}}_{\A}^{\ell \err} \cdot \nv^{\ell \err} \\
    		& + \halb \, \sum \limits_{\err \in \mathcal{N}^\ell} \frac{|\partial \omega^{\ell \err}|}{|\omega^\ell|} \, \left( \left((|\A| v_k)^\err - (|\A| v_k)^\ell \right) + \left( \w^\ell \cdot \mathbf{f}_{\A,k}^\ell - \w^\err \cdot \mathbf{f}_{\A,k}^\err \right) \right)\, n_k^{\ell \err} \\
    		& + \halb \, \sum \limits_{\err \in \mathcal{N}^\ell} \frac{|\partial \omega^{\ell \err}|}{|\omega^\ell|} \,  \left(
    		\w^\ell \cdot \mathcal{D}_{\A}^{\ell \err} + \w^\err \cdot \mathcal{D}_{\A}^{\err \ell}
    		\right) \cdot \nv^{\ell \err} \\
    		& + \halb \, \sum \limits_{\err \in \mathcal{N}^\ell} \frac{|\partial \omega^{\ell \err}|}{|\omega^\ell|} \,  \left(  ( \w^\ell + \w^\err ) \cdot {\mathcal{G}}_{\A}^{\ell \err} \cdot \nv^{\ell \err} + ( \w^\ell - \w^\err ) \cdot {\mathcal{G}}_{\A}^{\ell \err} \cdot \nv^{\ell \err} \right) \\
    		& = \Pi_\A^\ell.
    	\end{aligned}
    \end{equation}  
    By virtue of the discrete Gauss theorem \eqref{eqn.Gaussdiv}, we can add to the left hand side of the above equation a zero term given by 
    $$ \sum \limits_{\err \in \mathcal{N}^\ell} \frac{|\partial \omega^{\ell \err}|}{|\omega^\ell|} \, \left( (|\A| v_k)^\ell - \w^\ell \cdot \mathbf{f}_{\A,k}^\ell \right) \, n_k^{\ell \err} = 0, $$ and we reformulate the numerical dissipation according to \eqref{eqn.dissA1}-\eqref{eqn.dissA2} with the Roe-type property \eqref{eqn.Roe_p}, hence obtaining
    \begin{equation}
    	\label{eqn.th1_7}
    	\begin{aligned}
    		\frac{\partial |\A|^\ell}{\partial t} & + \halb \, \sum \limits_{\err \in \mathcal{N}^\ell} \frac{|\partial \omega^{\ell \err}|}{|\omega^\ell|} \, ( \w^\ell + \w^\err ) \cdot \tilde{\mathcal{F}}_{\A}^{\ell \err} \cdot \nv^{\ell \err} \\
    		& + \halb \, \sum \limits_{\err \in \mathcal{N}^\ell} \frac{|\partial \omega^{\ell \err}|}{|\omega^\ell|} \, \left((|\A| v_k)^\err + (|\A| v_k)^\ell \right) \, n_k^{\ell \err} \\
    		& - \halb \, \sum \limits_{\err \in \mathcal{N}^\ell} \frac{|\partial \omega^{\ell \err}|}{|\omega^\ell|} \, \left( \w^\ell \cdot \mathbf{f}_{\A,k}^\ell + \w^\err \cdot \mathbf{f}_{\A,k}^\err \right)\, n_k^{\ell \err} \\
    		& + \halb \, \sum \limits_{\err \in \mathcal{N}^\ell} \frac{|\partial \omega^{\ell \err}|}{|\omega^\ell|} \,  \left(
    		\w^\ell \cdot \mathcal{D}_{\A}^{\ell \err} + \w^\err \cdot \mathcal{D}_{\A}^{\err \ell}
    		\right) \cdot \nv^{\ell \err} \\
    		& + \phantom{\halb} \, \sum \limits_{\err \in \mathcal{N}^\ell} \frac{|\partial \omega^{\ell \err}|}{|\omega^\ell|} \, \halb \, \epsilon^{\ell \err} ( \A^\err - \A^\ell ) \cdot  \partial^2_{\A \A} |\A|^{\ell \err} \, ( \A^\err - \A^\ell ) - \epsilon^{\ell \err} \, ( |\A|^\err - |\A|^\ell ) \\
    		& = \Pi_\A^\ell.
    	\end{aligned}
    \end{equation}  
    The last term on the left hand side partially cancels with the production term $\Pi_\A^\ell$, that follows by the definition \eqref{eqn.prodA}. Therefore, the extra conservation law \eqref{eqn.detApde} is satisfied by defining the following fluxes in the semi-discrete finite volume scheme \eqref{eqn.ddetAdt}:
    \begin{equation}
    	\label{eqn.comp_flux_detA}
    	\begin{aligned}
    	F_{|\A|}^\ell \cdot \nv^{\ell \err} &= ( (|\A| v_k)^\ell - \w^\ell \cdot \mathbf{f}_{\A,k}^\ell ) \, n_k^{\ell \err} + \w^\ell \cdot \left(  \tilde{\mathcal{F}}_{\A}^{\ell \err} + \mathcal{D}_{\A}^{\ell \err} \right) \cdot \nv^{\ell \err} + 2 \, \epsilon^{\ell \err} \, |\A|^\ell, \\
    	F_{|\A|}^\err \cdot \nv^{\ell \err} &= ( (|\A| v_k)^\err - \w^\err \cdot \mathbf{f}_{\A,k}^\err ) \, n_k^{\ell \err} + \w^\err \cdot \left(  \tilde{\mathcal{F}}_{\A}^{\ell \err} + \mathcal{D}_{\A}^{\err \ell} \right) \cdot \nv^{\ell \err} - 2 \, \epsilon^{\ell \err} \, |\A|^\err.
	    \end{aligned}
    \end{equation}
\end{proof}

\subsection{Thermodynamically compatible finite volume scheme}
After achieving compatibility with the extra conservation law \eqref{eqn.detApde}, the semi-discrete finite volume scheme must be modified again to be compliant with the Second Law of Thermodynamics, meaning that it must fulfill also the entropy balance \eqref{eqn.entropy}. This is equivalent to satisfy the Gibbs relation \eqref{eqn.Gibbs2}, implying that we need to work with all the state variables $\u$ plus the density. However, thanks to the geometrically compatible discretization achieved for the distortion tensor $\A$, we can deduce the density directly from the determinant of $\A$ as $\rho=\rho_0|\A|$, therefore the full vector of state variables is simply given by $\tilde{\u}=(\rho_0|\A|,\u)^\top=(\rho_0|\A|,\rho v_i,\mathcal{E},A_{ik},J_k)^\top$. The thermodynamic correction is carried out in analogy with the one employed for the geometric compatibility, hence we introduce a modified set of numerical fluxes of the form
\begin{equation}
	\label{eqn.modifiedfluxS}	
	\hat{\mathcal{F}}(\tilde{\u}^\ell,\tilde{\u}^\err) \cdot \nv^{\ell \err} = \tilde{\mathcal{F}}(\tilde{\u}^\ell,\tilde{\u}^\err) - \alpha_S^{\ell \err} \,  ( \tilde{\r}^\err - \tilde{\r}^\ell ),
\end{equation} 
where $\alpha_S^{\ell \err}$ is a scalar correction factor that must be determined to obtain thermodynamic compatibility. The fluxes $\tilde{\mathcal{F}}(\tilde{\u}^\ell,\tilde{\u}^\err) \cdot \nv^{\ell \err}$ coincide with the central fluxes for all equations of \eqref{eqn.GPR_red} except for the distortion tensor equation, for which they are given by \eqref{eqn.modifiedfluxA}. Furthermore, the flux in the continuity equation \eqref{eqn.conti} is computed from the compatible fluxes \eqref{eqn.comp_flux_detA} of the semi-discrete equation for $|\A|$ given by \eqref{eqn.ddetAdt} upon multiplication by $\rho_0^\ell$. Consequently, we have that
\begin{equation}
	\label{eqn.ftilde}
	\tilde{\mathcal{F}}(\tilde{\u}^\ell,\tilde{\u}^\err) \cdot \nv^{\ell \err} = \left\{ \begin{array}{ll}
		\halb \, \rho_0^\ell \, \left(F_{|\A|}^\ell + F_{|\A|}^\err \right) \cdot \nv^{\ell \err} & \text{ for \eqref{eqn.conti}}  \text{ with fluxes \eqref{eqn.comp_flux_detA}} \\ [2mm]
		\mathcal{F}(\u^\ell,\u^\err) \cdot \nv^{\ell \err} & \text{ for \eqref{eqn.momentum_red}-\eqref{eqn.energy_red}-\eqref{eqn.heatflux_red}} \\ [2mm]
		\tilde{\mathcal{F}}_{\A}(\u^\ell,\u^\err) \cdot \nv^{\ell \err} & \text{ for \eqref{eqn.deformation_red}}
	\end{array} \right. .
\end{equation}
Since the compliance with the GCL must not be destroyed by this new modification, we deliberately choose to add the correction factor $\alpha_S^{\ell \err}$ to only a subset of dual variables among $\r$ in \eqref{eqn.dualvar_r}, which is referred to as $\tilde{\r}$ in \eqref{eqn.modifiedfluxS}. Specifically, we allow the momentum and the thermal impulse equations to account for the thermodynamic compatibility, thus defining
\begin{equation}
	\tilde{\r} = \frac{1}{T} \, \left(0,\partial_{\rho \vv} \mathcal{E}, 0, \mathbf{0}, \partial_{\J} \mathcal{E}  \right)^\top = \{\tilde{r}_j\} = \frac{1}{T} \, (0,-v_i,0,0_{ik},-\beta_k)^\top.
\end{equation}
We remark that the equations for density and distortion tensor are not affected by the thermodynamic corrections since they already carry the geometric compatibility correction. Let us also note that the total energy equation is not modified in order to maintain stationary solutions of the governing PDE, as explained at the end of this section.
To compute the correction factor $\alpha_S^{\ell \err}$ we can now proceed along the lines of \cite{Abgrall2018,AbgrallOeffnerRanocha}, hence requiring that the sum of all the fluctuations across an element interface is equal to the flux difference of the entropy equation \eqref{eqn.entropy}, thus leading to
\begin{equation}
	\label{eqn.entropy_cond}
	\begin{aligned}
		\r^\ell \cdot \left( \hat{\mathcal{F}}(\tilde{\u}^\ell,\tilde{\u}^\err) \cdot \nv^{\ell \err} - \mathbf{f}_{k}^\ell \cdot n_k^{\ell \err}\right) + \r^\err \cdot \left( \mathbf{f}_{k}^\err \cdot n_k^{\ell \err} - \hat{\mathcal{F}}(\tilde{\u}^\ell,\tilde{\u}^\err) \cdot \nv^{\ell \err} \right) & +& \\ \r^\ell \cdot \mathcal{D}(\tilde{\u}^\ell,\tilde{\u}^\err) \cdot \nv^{\ell \err} + \r^\err \cdot \mathcal{D}(\tilde{\u}^\err,\tilde{\u}^\ell) \cdot \nv^{\ell \err}  & =& \left((\rhoS v_k + \beta_k)^\err - (\rhoS v_k + \beta_k)^\ell \right) \, n_k^{\ell \err}.
	\end{aligned}
\end{equation}
By employing the flux definition \eqref{eqn.modifiedfluxS} in the conservation condition \eqref{eqn.entropy_cond}, the thermodynamic correction scalar $\alpha_S^{\ell \err}$ is found to be given by
\begin{equation}
	\label{eqn.alphaS}
	\begin{aligned}
		\alpha_S^{\ell \err} & = \frac{\left((\rhoS v_k + \beta_k)^\err - (\rhoS v_k + \beta_k)^\ell \right) \, n_k^{\ell \err} + \left( \tilde{\mathcal{F}}(\tilde{\u}^\ell,\tilde{\u}^\err) \cdot \nv^{\ell \err} \right) \cdot \left(\r^\err-\r^\ell\right) - \left( \r^\err \cdot \mathbf{f}_{k}^\err - \r^\ell \cdot \mathbf{f}_{k}^\ell\right)\, n_k^{\ell \err}}{\left(\tilde{\r}^\err - \tilde{\r}^\ell\right)^2} \\
		& - \frac{\left( \r^\err + \r^\ell \right) \cdot \mathcal{D}(\tilde{\u}^\ell,\tilde{\u}^\err) \cdot \nv^{\ell \err}}{\left(\tilde{\r}^\err - \tilde{\r}^\ell\right)^2}.
	\end{aligned}
\end{equation}
Obviously, we set $\alpha_S^{\ell \err}=0$ if $\tilde{\r}^\ell=\tilde{\r}^\err$. The source terms in equations \eqref{eqn.deformation_red} and \eqref{eqn.heatflux_red} are compatible by construction also for thermodynamic compatibility. Indeed, multiplication of the sources by the dual variables $\partial_{A_{ik}^\ell} \mathcal{E} = \alpha_{ik}^\ell$ and  $\partial_{J_k^\ell} \mathcal{E} = \beta_{k}^\ell$, with negative sign and divided by the temperature according to \eqref{eqn.Gibbs2}, yields
\begin{equation}
	-\frac{\alpha_{ik}^\ell}{T^\ell} \, \left( -\dfrac{ \alpha_{ik}^\ell }{\theta_1(\tau_1)} \right) - \frac{\beta_k^\ell}{T^\ell} \, \left( -\dfrac{\beta_k^\ell}{\theta_2(\tau_2)} \right) = \dfrac{\alpha_{ik}^\ell \alpha_{ik}^\ell }{\theta_1(\tau_1) T^\ell}  + \dfrac{\beta_k^\ell \beta_k^\ell}{\theta_2(\tau_2) T^\ell}  \geq 0, 
\end{equation}
which is exactly the source term in the entropy equation \eqref{eqn.entropy}.

The geometrically and thermodynamically compatible semi-discrete finite volume scheme without numerical dissipation is then given by
\begin{equation}
	\label{eqn.hgtc_no_diss}
	\frac{\partial \u^\ell}{\partial t} + \sum \limits_{\err \in \mathcal{N}^\ell} \frac{|\partial \omega^{\ell \err}|}{|\omega^\ell|} \left( \hat{\mathcal{F}}(\u^\ell,\u^\err)^{\ell \err} + \mathcal{D}(\u^\ell,\u^\err)^{\ell \err} \right) \cdot \nv^{\ell \err} = \mathbf{S}(\u^\ell),
\end{equation}
with the definition of the compatible fluxes \eqref{eqn.modifiedfluxS}-\eqref{eqn.ftilde} and the thermodynamic correction factor \eqref{eqn.alphaS}.

\begin{remark}[On the energy equation in the thermodynamic correction]
	The subset of dual variables $\tilde{\r}$ in the flux correction \eqref{eqn.modifiedfluxS} does not take into account the dual variable in the energy equation, which is equal to $1/T$. Without loss of generality, let us consider a computational domain $\Omega=[-L;L]^2$ with periodic boundaries and $L\in\mathbb{R}$, and the following initial condition with only a discontinuity in the density field located at $||\xx||=R_0 \subset \Omega$:
	\begin{equation}
		\label{eqn.ic0}
		\rho(\xx,t=0) = \left\{ \begin{array}{ll}
			\rho^L_0 & \text{ for } ||\xx|| \leq R_0 \\ [1.5mm]
			\rho^R_0 & \text{ for } ||\xx|| > R_0
		\end{array} \right. , \quad \vv(\xx,t=0)=\mathbf{0}, \quad p(\xx,t=0)=p_0, \quad \A(\xx,t=0)=\mathbf{I}, \quad \mathbf{J}(\xx,t=0) = \mathbf{0}.
	\end{equation}
	In this case, the total energy is only given by the internal energy contribution which is constant, namely
	\begin{equation}
		\label{eqn.en0}
		\mathcal{E}(\xx,t=0) = \frac{p_0}{\gamma - 1},
	\end{equation}
	therefore the semi-discrete finite volume scheme \eqref{eqn.hgtc_no_diss} yields
	\begin{equation}
		\label{eqn.dudtconst}
		\frac{\partial \u^\ell}{\partial t} = \mathbf{0},
	\end{equation}
	and the initial condition \eqref{eqn.ic0} also represents the exact solution. The correction factor is \textit{a priori} nonzero, i.e. $\alpha_S^{\ell \err} \neq 0$, because of the pressure and temperature terms in the momentum and thermal impulse equations, respectively. If we add the thermodynamic correction in the energy equation, the associated flux $\hat{\mathcal{F}}_{\mathcal{E}}(\u^\ell,\u^\err) \cdot \nv^{\ell \err}$ should be corrected with the jump term in the dual variable $1/T$, leading to 
	\begin{equation}
		\label{eqn.artfluxE}
		\hat{\mathcal{F}}_{\mathcal{E}}(\u^\ell,\u^\err) \cdot \nv^{\ell \err} = {\mathcal{F}}_{\mathcal{E}}(\u^\ell,\u^\err) \cdot \nv^{\ell \err} -\alpha_S^{\ell \err} \left(\frac{1}{T^\err} - \frac{1}{T^\ell}\right).
	\end{equation}
	This would no longer preserve the constant energy density \eqref{eqn.en0} because the discontinuity in the density profile causes a jump in the temperature which is defined as
	\begin{equation}
		T(\xx,t=0) = \frac{p_0}{\rho(\xx,t=0) \, c_v \, (\gamma-1)}.
	\end{equation}
	As a consequence, the artificial flux \eqref{eqn.artfluxE} is not physical and thus the thermodynamic dual variable $1/T$ is not included in the dual vector $\tilde{\r}$ for computing the scalar factor $\alpha_S^{\ell \err}$ in \eqref{eqn.alphaS}, so that the geometrically and thermodynamically compatible scheme still maintains this physical equilibrium, i.e. we obtain again the correct stationary solution given by \eqref{eqn.dudtconst}.
\end{remark}

\paragraph{Thermodynamic compatibility with dissipation terms} As done for the geometric compatibility, to ensure the stability of the scheme in case of discontinuous solutions we also take into account the parabolic dissipation terms \eqref{eqn.diss}, thus we supplement the geometrically and thermodynamically compatible scheme \eqref{eqn.hgtc_no_diss} with the dissipative fluxes $\mathcal{G}(\u^\ell,\u^\err) \cdot \nv^{\ell \err}$, hence giving rise to the non-negative production term $\Pi^\ell$ in the entropy inequality \eqref{eqn.entropy}.

\begin{theorem}[Thermodynamic compatibility] \label{th2} 
    The semi-discrete finite volume scheme for the reduced model \eqref{eqn.GPR_red}
	\begin{equation}
		\label{eqn.hgtc}
		\frac{\partial \u^\ell}{\partial t} + \sum \limits_{\err \in \mathcal{N}^\ell} \frac{|\partial \omega^{\ell \err}|}{|\omega^\ell|} \, \left( \hat{\mathcal{F}}(\u^\ell,\u^\err)^{\ell \err} + \mathcal{D}(\u^\ell,\u^\err)^{\ell \err} + \mathcal{G}(\A^\ell,\A^\err) \right) \cdot \nv^{\ell \err} = \mathbf{S}(\u^\ell),
	\end{equation} 
	with the geometrically and thermodynamically compatible fluxes \eqref{eqn.modifiedfluxS}-\eqref{eqn.ftilde}, the non-conservative products \eqref{eqn.ncprod} and the dissipation terms \eqref{eqn.diss}, satisfies the extra conservation law \eqref{eqn.entropy} with the following conservative semi-discrete scheme:
	\begin{equation}
		\label{eqn.entropyfv}
		\frac{\partial \rhoS^\ell}{\partial t}
		+ \sum \limits_{\err \in \mathcal{N}^\ell} \frac{|\partial \omega^{\ell \err}|}{|\omega^\ell|} \, \halb \left(F_{\rhoS}^\ell + F_{\rhoS}^\err \right) \cdot \nv^{\ell \err} = \Pi^\ell  + \dfrac{\alpha_{ik}^\ell \alpha_{ik}^\ell }{\theta_1(\tau_1) T^\ell}  + \dfrac{\beta_k^\ell \beta_k^\ell}{\theta_2(\tau_2) T^\ell}.
	\end{equation}
    Furthermore, assuming $T^\ell>0$ and $\partial_{\u \u} \rhoS^{\ell \err} \leq 0$, the right hand side of the entropy balance is non-negative: 
    \begin{equation}
    	\label{eqn.posPi}
    	\Pi^\ell  + \dfrac{\alpha_{ik}^\ell \alpha_{ik}^\ell }{\theta_1(\tau_1) T^\ell}  + \dfrac{\beta_k^\ell \beta_k^\ell}{\theta_2(\tau_2) T^\ell}  \geq 0,
    \end{equation}
    therefore the scheme \eqref{eqn.hgtc} also satisfies a cell entropy inequality.
\end{theorem}

\begin{proof}
	The proof is similar to the one already carried out for Theorem \ref{th1}. We introduce the abbreviations 
	\begin{equation}
		\label{eqn.th2_1}
		\hat{\mathcal{F}}^{\ell \err}:=\hat{\mathcal{F}}^{\ell \err}(\u^\ell,\u^\err), \qquad \mathcal{D}^{\ell \err}:=\mathcal{D}(\u^\ell,\u^\err), \qquad \mathcal{G}^{\ell \err}:=\mathcal{G}(\u^\ell,\u^\err).
	\end{equation}
    We consider the dual variables $\tilde{\u}=(\rho_0|\A|,\u)^\top$ where the density is directly evaluated from the determinant of the distortion tensor, and the associated fluxes are computed using \eqref{eqn.comp_flux_detA} multiplied by $\rho_0^\ell$ in the GCL equation \eqref{eqn.ddetAdt}. With a little abuse of notation, let us omit the tilde symbol and assume that the additional fluxes related to the density equation are embedded in the flux tensor, thus we will simply use $\u^\ell$ (and thus the dual variables $\r^\ell$) and $\hat{\mathcal{F}}^{\ell \err}$. After dot multiplying the semi-discrete system \eqref{eqn.hgtc} by the dual variables $\r^\ell$, and adding and subtracting the terms $\halb \, \r^\err \cdot \hat{\mathcal{F}}^{\ell \err} \cdot \nv^{\ell \err}$, $\halb \, \r^\err \cdot {\mathcal{D}}^{\err \ell} \cdot \nv^{\err \ell}$ and $\halb \, \r^\err \cdot {\mathcal{G}}^{\ell \err} \cdot \nv^{\ell \err}$, we obtain
    \begin{equation}
    	\label{eqn.th2_2}
    	\begin{aligned}
    		\frac{\partial \rhoS^\ell}{\partial t} & + \halb \, \sum \limits_{\err \in \mathcal{N}^\ell} \frac{|\partial \omega^{\ell \err}|}{|\omega^\ell|} \, \left( ( \r^\ell + \r^\err ) \cdot \hat{\mathcal{F}}^{\ell \err} \cdot \nv^{\ell \err} + ( \r^\ell - \r^\err ) \cdot \hat{\mathcal{F}}^{\ell \err} \cdot \nv^{\ell \err} \right) \\
    		& + \halb \, \sum \limits_{\err \in \mathcal{N}^\ell} \frac{|\partial \omega^{\ell \err}|}{|\omega^\ell|} \,  \left( \r^\ell \cdot \mathcal{D}^{\ell \err} \cdot \nv^{\ell \err} + \, \r^\err \cdot \mathcal{D}^{\err \ell} \cdot \nv^{\err \ell} + \r^\ell \cdot \mathcal{D}^{\ell \err} \cdot \nv^{\ell \err} - \, \r^\err \cdot \mathcal{D}^{\err \ell} \cdot \nv^{\err \ell} \right) \\
    		& + \halb \, \sum \limits_{\err \in \mathcal{N}^\ell} \frac{|\partial \omega^{\ell \err}|}{|\omega^\ell|} \,  \left(  ( \r^\ell + \r^\err ) \cdot {\mathcal{G}}^{\ell \err} \cdot \nv^{\ell \err} + ( \r^\ell - \r^\err ) \cdot {\mathcal{G}}^{\ell \err} \cdot \nv^{\ell \err} \right) \\
    		& = \r^\ell \cdot \mathbf{S}(\u^\ell).
    	\end{aligned}	
    \end{equation}
    We analyze the compatibility of the source terms, which explicitly write
    \begin{equation}
    	\label{eqn.th2_3}
    		\r^\ell \cdot \mathbf{S}(\u^\ell) = -\frac{\alpha_{ik}^\ell}{T^\ell} \, \left( -\dfrac{ \alpha_{ik}^\ell }{\theta_1(\tau_1)} \right) - \frac{\beta_k^\ell}{T^\ell} \, \left( -\dfrac{\beta_k^\ell}{\theta_2(\tau_2)} \right) - \underbrace{\alpha_{ik}^\ell \Pi_\A^\ell \, \frac{w^\ell_{ik}}{w^\ell_{ik} \, w^\ell_{ik}}}_{=0} = \dfrac{\alpha_{ik}^\ell \alpha_{ik}^\ell }{\theta_1(\tau_1) T^\ell}  + \dfrac{\beta_k^\ell \beta_k^\ell}{\theta_2(\tau_2) T^\ell}, 
    \end{equation}
    where the production term vanishes thanks to the relation $\alpha_{ik}^\ell w_{ik}^\ell = 0$ as proven for the geometric compatibility (see \ref{app.gcl}). Thus we retrieve the source terms of the entropy balance law \eqref{eqn.entropyfv}. Back to equation \eqref{eqn.th2_2}, we use the compatibility condition \eqref{eqn.entropy_cond} to rewrite the term $( \r^\ell - \r^\err ) \cdot \hat{\mathcal{F}}^{\ell \err} \cdot \nv^{\ell \err}$, and we get
    \begin{equation}
    	\label{eqn.th2_4}
    	\begin{aligned}
    		\frac{\partial \rhoS^\ell}{\partial t} & + \halb \, \sum \limits_{\err \in \mathcal{N}^\ell} \frac{|\partial \omega^{\ell \err}|}{|\omega^\ell|} \, ( \r^\ell + \r^\err ) \cdot \hat{\mathcal{F}}^{\ell \err} \cdot \nv^{\ell \err} \\
    		& + \halb \, \sum \limits_{\err \in \mathcal{N}^\ell} \frac{|\partial \omega^{\ell \err}|}{|\omega^\ell|} \, \left( \left((\rhoS v_k + \beta_k)^\err - (\rhoS v_k + \beta_k)^\ell \right) + \left( \r^\ell \cdot \mathbf{f}_k^\ell - \r^\err \cdot \mathbf{f}_k^\err \right) \right)\, n_k^{\ell \err} \\
    		& + \halb \, \sum \limits_{\err \in \mathcal{N}^\ell} \frac{|\partial \omega^{\ell \err}|}{|\omega^\ell|} \,  \left(
    		\r^\ell \cdot \mathcal{D}^{\ell \err} + \r^\err \cdot \mathcal{D}^{\err \ell}
    		\right) \cdot \nv^{\ell \err} \\
    		& + \halb \, \sum \limits_{\err \in \mathcal{N}^\ell} \frac{|\partial \omega^{\ell \err}|}{|\omega^\ell|} \,  \left(  ( \r^\ell + \r^\err ) \cdot {\mathcal{G}}^{\ell \err} \cdot \nv^{\ell \err} + ( \r^\ell - \r^\err ) \cdot {\mathcal{G}}^{\ell \err} \cdot \nv^{\ell \err} \right) \\
    		& = \dfrac{\alpha_{ik}^\ell \alpha_{ik}^\ell }{\theta_1(\tau_1) T^\ell}  + \dfrac{\beta_k^\ell \beta_k^\ell}{\theta_2(\tau_2) T^\ell}.
    	\end{aligned}
    \end{equation} 
    Adding on the left hand side term $$ \halb \, \sum \limits_{\err \in \mathcal{N}^\ell} \frac{|\partial \omega^{\ell \err}|}{|\omega^\ell|} \, \left( (\rhoS v_k + \beta_k)^\ell  - \r^\ell \cdot \mathbf{f}_{k}^\ell \right) \, n_k^{\ell \err}=0,$$ which corresponds to a zero contribution thanks to the property \eqref{eqn.Gaussdiv}, leads to
    \begin{equation}
    	\label{eqn.th2_6}
    	\begin{aligned}
    		\frac{\partial \rhoS^\ell}{\partial t} & + \halb \, \sum \limits_{\err \in \mathcal{N}^\ell} \frac{|\partial \omega^{\ell \err}|}{|\omega^\ell|} \, ( \r^\ell + \r^\err ) \cdot \hat{\mathcal{F}}^{\ell \err} \cdot \nv^{\ell \err} \\
    		& + \halb \, \sum \limits_{\err \in \mathcal{N}^\ell} \frac{|\partial \omega^{\ell \err}|}{|\omega^\ell|} \, \left((\rhoS v_k + \beta_k)^\err + (\rhoS v_k + \beta_k)^\ell \right)\, n_k^{\ell \err} \\
    		& - \halb \, \sum \limits_{\err \in \mathcal{N}^\ell} \frac{|\partial \omega^{\ell \err}|}{|\omega^\ell|} \, \left( \r^\ell \cdot \mathbf{f}_k^\ell + \r^\err \cdot \mathbf{f}_k^\err \right)\, n_k^{\ell \err} \\
    		& + \halb \, \sum \limits_{\err \in \mathcal{N}^\ell} \frac{|\partial \omega^{\ell \err}|}{|\omega^\ell|} \,  \left(
    		\r^\ell \cdot \mathcal{D}^{\ell \err} + \r^\err \cdot \mathcal{D}^{\err \ell}
    		\right) \cdot \nv^{\ell \err} \\
    		& + \halb \, \sum \limits_{\err \in \mathcal{N}^\ell} \frac{|\partial \omega^{\ell \err}|}{|\omega^\ell|} \,  \left(  ( \r^\ell + \r^\err ) \cdot {\mathcal{G}}^{\ell \err} \cdot \nv^{\ell \err} + ( \r^\ell - \r^\err ) \cdot {\mathcal{G}}^{\ell \err} \cdot \nv^{\ell \err} \right) \\
    		& = \dfrac{\alpha_{ik}^\ell \alpha_{ik}^\ell }{\theta_1(\tau_1) T^\ell}  + \dfrac{\beta_k^\ell \beta_k^\ell}{\theta_2(\tau_2) T^\ell}.
    	\end{aligned}
    \end{equation} 
    Relying on the same reasoning applied for the geometric compatibility, the dissipation terms $\r^\ell \cdot \mathcal{G}^{\ell \err} \cdot \nv^{\ell \err}$ can be rearranged as in \eqref{eqn.dissA1}, that is
    \begin{equation}
    	\label{eqn.dissS1}
    	\begin{aligned}
    		\r^\ell \cdot \mathcal{G}^{\ell \err} \cdot \nv^{\ell \err} &= \halb ( \r^\err - \r^\ell ) \cdot \epsilon^{\ell \err} ( \u^\err - \u^\ell ) - \halb ( \r^\err + \r^\ell ) \cdot \epsilon^{\ell \err} ( \u^\err - \u^\ell ).
    	\end{aligned}
    \end{equation}
    Likewise in \eqref{eqn.dissA2}, due to the path integral
     \begin{equation}
    	\label{eqn.Roe_r}
    	\int \limits_{\u^\ell}^{\u^\err} \r \cdot d\u = \int \limits_{\u^\ell}^{\u^\err} \partial_{\u} \rhoS \cdot d\u = \rhoS^\err - \rhoS^\ell, 
    \end{equation} 
    we interpret the second term on the right hand side as an approximation of the jump term in the entropy variables, thus
    \begin{equation}
    	\label{eqn.dissS2}
    	-\halb ( \r^\err + \r^\ell ) \cdot \epsilon^{\ell \err} (\u^\err - \u^\ell ) \approx - \epsilon^{\ell \err} ( \rhoS^\err - \rhoS^\ell ).
    \end{equation}
   The jump in the dual variables present in the first term in \eqref{eqn.dissS1} is converted into a jump in the state variables by introducing the Hessian matrix $\partial^2_{\u \u} \rhoS^{\ell \err}$ which verifies the Roe property 
    \begin{equation}
    	\label{eqn.HessianS}
    	\partial^2_{\u \u} \rhoS^{\ell \err} \cdot ( \u^\err - \u^\ell ) = \r^\err - \r^\ell.
    \end{equation}
    Therefore, using \eqref{eqn.dissS2} and \eqref{eqn.HessianS} in \eqref{eqn.th2_6}, we arrive at
    \begin{equation}
    	\label{eqn.th2_7}
    	\begin{aligned}
    		\frac{\partial \rhoS^\ell}{\partial t} & + \halb \, \sum \limits_{\err \in \mathcal{N}^\ell} \frac{|\partial \omega^{\ell \err}|}{|\omega^\ell|} \, ( \r^\ell + \r^\err ) \cdot \hat{\mathcal{F}}^{\ell \err} \cdot \nv^{\ell \err} \\
    		& + \halb \, \sum \limits_{\err \in \mathcal{N}^\ell} \frac{|\partial \omega^{\ell \err}|}{|\omega^\ell|} \, \left((\rhoS v_k + \beta_k)^\err + (\rhoS v_k + \beta_k)^\ell \right)\, n_k^{\ell \err} \\
    		& - \halb \, \sum \limits_{\err \in \mathcal{N}^\ell} \frac{|\partial \omega^{\ell \err}|}{|\omega^\ell|} \, \left( \r^\ell \cdot \mathbf{f}_k^\ell + \r^\err \cdot \mathbf{f}_k^\err \right)\, n_k^{\ell \err} \\
    		& + \halb \, \sum \limits_{\err \in \mathcal{N}^\ell} \frac{|\partial \omega^{\ell \err}|}{|\omega^\ell|} \,  \left(
    		\r^\ell \cdot \mathcal{D}^{\ell \err} + \r^\err \cdot \mathcal{D}^{\err \ell}
    		\right) \cdot \nv^{\ell \err} \\
    		& - \phantom{\halb} \, \sum \limits_{\err \in \mathcal{N}^\ell} \frac{|\partial \omega^{\ell \err}|}{|\omega^\ell|} \epsilon^{\ell \err} \, ( \rhoS^\err - \rhoS^\ell ) \\
    		& = \Pi^\ell + \dfrac{\alpha_{ik}^\ell \alpha_{ik}^\ell }{\theta_1(\tau_1) T^\ell}  + \dfrac{\beta_k^\ell \beta_k^\ell}{\theta_2(\tau_2) T^\ell},
    	\end{aligned}
    \end{equation}
    with the production term given by
    \begin{equation}
    	\Pi^\ell = -\sum \limits_{\err \in \mathcal{N}^\ell} \frac{|\partial \omega^{\ell \err}|}{|\omega^\ell|} \, \halb \, \epsilon^{\ell \err} ( \u^\err - \u^\ell ) \cdot  \partial^2_{\u \u} \rhoS^{\ell \err} \, ( \u^\err - \u^\ell ).
    \end{equation}
     The compatibility with the extra conservation law \eqref{eqn.entropy} is then achieved by defining the following fluxes in the semi-discrete finite volume scheme \eqref{eqn.entropyfv}:
     \begin{equation}
     	\label{eqn.comp_flux_S}
     	\begin{aligned}
     		F_{\rhoS}^\ell \cdot \nv^{\ell \err} &= ( (\rhoS v_k + \beta_k)^\ell - \r^\ell \cdot \mathbf{f}_{k}^\ell ) \, n_k^{\ell \err} + \r^\ell \cdot \left(  \hat{\mathcal{F}}^{\ell \err} + \mathcal{D}^{\ell \err} \right) \cdot \nv^{\ell \err} + 2 \, \epsilon^{\ell \err} \, \rhoS^\ell, \\
     		F_{\rhoS}^\err \cdot \nv^{\ell \err} &= ( (\rhoS v_k + \beta_k)^\err - \r^\err \cdot \mathbf{f}_{k}^\err ) \, n_k^{\ell \err} + \r^\err \cdot \left(  \hat{\mathcal{F}}^{\ell \err} + \mathcal{D}^{\err \ell} \right) \cdot \nv^{\ell \err} - 2 \, \epsilon^{\ell \err} \, \rhoS^\err.
     	\end{aligned}
     \end{equation}
     Finally, in the presence of numerical viscosity, i.e. when $\epsilon^{\ell \err}>0$, the entropy inequality is retrieved since the resulting term on the right hand side of \eqref{eqn.th2_7} is non-negative, meaning that the positivity condition \eqref{eqn.posPi} is fulfilled  due to the assumptions $\theta_1>0$, $\theta_2>0$, $T^\ell >0$ and $\partial^2_{\u \u} \rhoS^{\ell \err} \leq 0$. The cell entropy inequality is thus satisfied at the semi-discrete level by the finite volume scheme \eqref{eqn.hgtc}.
\end{proof}

\subsection{Time discretization}
The explicit time marching algorithm is given by Runge-Kutta schemes that are listed in \ref{app.rk} for order one, two and four. The associated time step is computed according to a classical CFL-type stability condition based on the maximum hyperbolic eigenvalue estimate given by \eqref{eqn.lambda} and the maximum viscous eigenvalue related to the parabolic dissipative terms:
\begin{equation}
	\label{eqn.timestep}
	\Delta t \leq \text{CFL} \frac{\min \limits_{\ell \in N^\ell} h^\ell}{\max \limits_{\ell \in N^\ell} \left( |\lambda^\ell| + 2 \frac{\epsilon^\ell}{h^\ell} \right)},
\end{equation}
where $\text{CFL}$ is the Courant-Friedrichs-Lewy number and $h^\ell=\sqrt{|\omega^\ell|}$ is the characteristic cell size.

\section{Numerical results} \label{sec.results}
In this section, we propose a suite of test cases aiming at validating the accuracy and the robustness of the novel Hyperbolic Geometrically and Thermodynamically Compatible finite volume schemes \eqref{eqn.hgtc}, which will be labeled as HGTC. We demonstrate that the compatibility is preserved at the semi-discrete level up to the order of the time integrator, and we systematically measure the errors of mass conservation ($\varepsilon_{A}$) and total entropy balance ($\varepsilon_{S}$).
More precisely, we monitor over time the following quantities in $L_{\infty}$ norm over the entire computational domain $\Omega$:
\begin{equation}
	\delta_\A = \left\| \, \detA - \rho/\rho_0 \, \right\|_{\infty}, \qquad \delta_S = \left\| \, {\rhoS} - {\mathcal{S}}(\rho,p) \, \right\|_{\infty}.
	\label{eqn.epsA_S}
\end{equation}
For $\delta_A$, the quantity $|\A|$ is computed by evaluating the determinant of the distortion tensor $\A$ by using the components $A_{ik}$ that are evolved according to the semi-discrete scheme \eqref{eqn.dAdt}, whereas the quantity $\rho/\rho_0$ is obtained with $|\A|$ taken from the solution of the extra conservation law \eqref{eqn.detApde} discretized by the scheme \eqref{eqn.ddetAdt} with fluxes \eqref{eqn.comp_flux_detA}. For $\delta_S$, ${\rhoS}$ is the total entropy computed from the entropy equation \eqref{eqn.entropy} solved as an extra conservation law with the semi-discrete scheme \eqref{eqn.entropyfv}, while ${\mathcal{S}}(\rho,p)$ is evaluated from the equation of state given by $\mathcal{E}_1$ in \eqref{eq:energies}, namely
\begin{equation}
	{\mathcal{S}}(\rho,p) = \rho \, \log\left( \frac{p}{\rho^{\gamma}}  \, c_v \right), \qquad \rho=\rho_0 \, |\A|.
\end{equation}	
As such, the structure-preserving properties of the scheme are numerically investigated. If not stated otherwise, we set the CFL number to $\text{CFL}=0.5$ in \eqref{eqn.timestep} and the polytropic index of the gas is assumed to be $\gamma=7/5$, whereas the specific heat at constant volume is always chosen to be $c_v=2.5$. Whenever a viscosity coefficient $\mu$ is specified, the relaxation time $\tau_1$ is computed according to $\mu = \frac{1}{6} \rho_0 c_s^2 \tau_1$. Likewise, if a heat conduction coefficient $\kappa$ is set, the corresponding relaxation time $\tau_2$ is evaluated from the asymptotic relation $\kappa = \rho_0 T_0 c_h^2 \tau_2$. In the other cases, no source terms are considered, thus we set $\tau_1=\tau_2=10^{20}$ hence retrieving the behavior of elastic solids without heat conduction. The distortion matrix is always initialized as $\A=\mathbf{I}$, and the thermal impulse is initially given $\J=\mathbf{0}$. The reference density and temperature are set to $\rho_0=\rho(\xx,t=0)$ and $T_0=1$, if not specified. We depict the absolute values of the correction factors $\alpha_\A$ and $\alpha_S$ in \eqref{eqn.alphaA} and \eqref{eqn.alphaS}, respectively, in order to better appreciate the order of magnitude and the location of the structure-preserving corrections.

\subsection{Numerical convergence studies}
The accuracy of the new HGTC schemes is verified on the isentropic vortex problem forwarded in \cite{HuShuVortex1999}.  The computational domain is the square $\Omega=[0;10]^2$ with periodic boundaries, and the generic radial position is $r=\sqrt{(x_1-5)^2+(x_2-5)^2}$. The parameters of the model are such that an ideal inviscid fluid is retrieved, hence we set $c_s=c_h=0$, and the initial condition is prescribed in terms of some perturbations that are superimposed on a background constant state:
\begin{equation}
	\rho(t=0,\xx) = (1+\delta T)^{\frac{1}{\gamma-1}}, \qquad \vv(t=0,\xx) = \mathbf{0}, \qquad p(t=0,\xx) = (1+\delta T)^{\frac{\gamma}{\gamma-1}},
	\label{eq.ConvEul-IC}
\end{equation}
with the perturbations for temperature $\delta T$ given by
\begin{equation}
	\delta T = -\frac{(\gamma-1)\epsilon^2}{8\gamma\pi^2}e^{1-r^2}.
\end{equation}
The simulation is carried out until the final time $t_f=0.25$ on a sequence of successively refined Voronoi meshes, and the errors are measured in $L_2$ norms and reported in Table \ref{tab.convRates}, showing that the formal order of accuracy is retrieved. No numerical dissipation is added to the scheme because the flow does not exhibit any discontinuity, thus $\epsilon^{\ell \err}=0$ in \eqref{eqn.diss}. 

\begin{table}[!htbp]  
	\caption{Numerical convergence results for the isentropic vortex problem using the HGTC scheme. The errors are measured in the $L_2$ norm and refer to the variables $\rho=\rho_0 \, |\A|$ (density), $v_1$ (horizontal velocity) and pressure $p$ at time $t_{f}=0.25$.}  
	\begin{center} 
		\begin{small}
			\renewcommand{\arraystretch}{1.2}
			\begin{tabular}{c|cccccc}
				$h$ & $\left\| \rho_0 \, |\A| \right\|_2$ & $O(\rho_0 \, |\A|)$ & $\left\| v_1 \right\|_2$ & $\mathcal{O}(v_1)$ & $\left\| p \right\|_2$ & $\mathcal{O}(p)$ \\ 
				\hline
				3.20E-01 & 6.2483E-02 & -    & 1.5675E-01 & -    & 7.9011E-02 & -    \\
				1.65E-01 & 3.1941E-02 & 1.01 & 7.9131E-02 & 1.03 & 4.0536E-02 & 1.01 \\
				1.09E-01 & 2.1427E-02 & 0.97 & 5.3292E-02 & 0.96 & 2.7139E-02 & 0.98 \\
				8.57E-02 & 1.6231E-02 & 1.14 & 3.9863E-02 & 1.19 & 2.0620E-02 & 1.13 \\
			\end{tabular}
		\end{small}
	\end{center}
	\label{tab.convRates}
\end{table}

\begin{figure}[!htbp]
	\begin{center}
		\begin{tabular}{cc} 
			\includegraphics[trim=5 0 0 5,clip,width=0.47\textwidth]{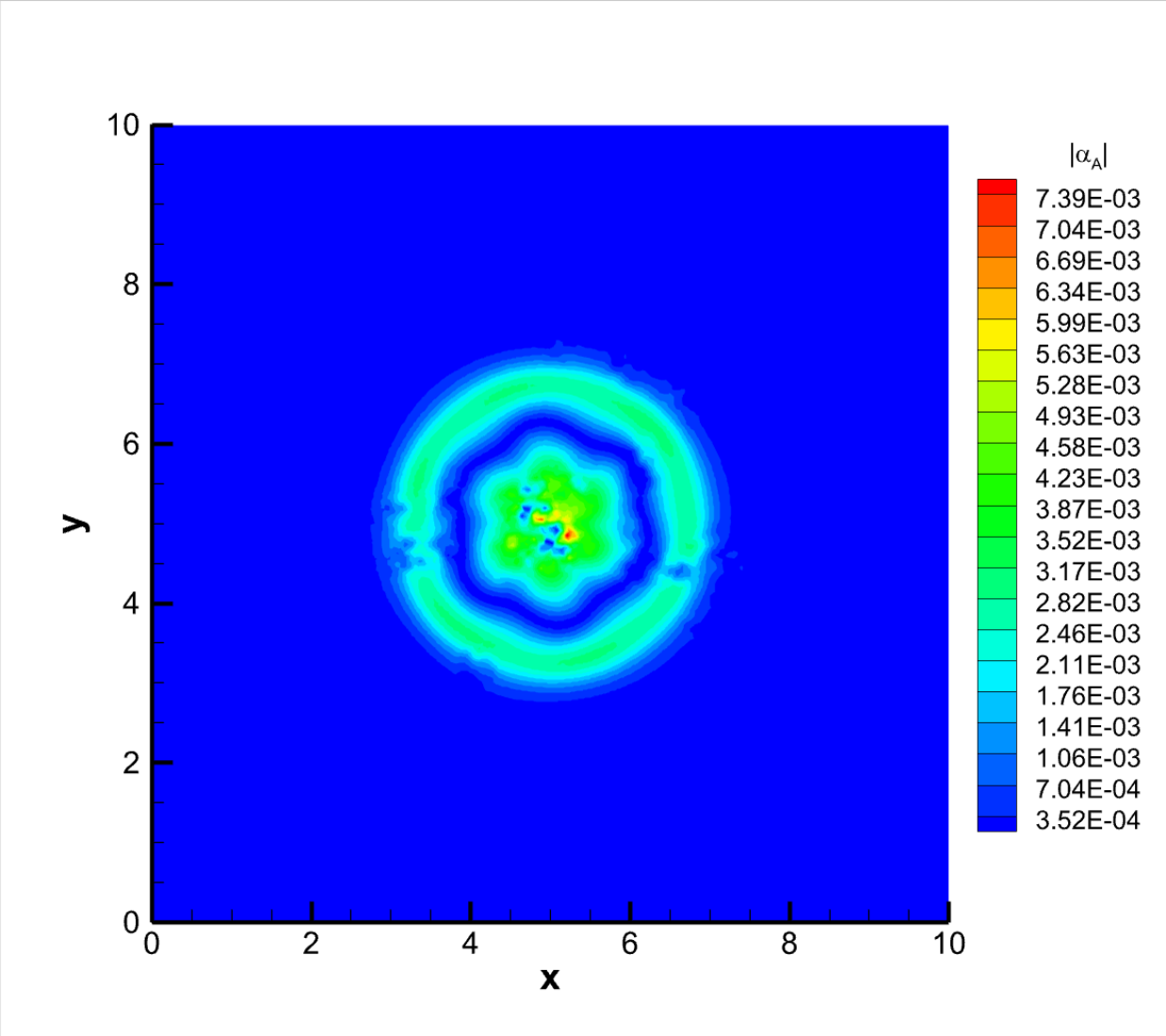} & 
			\includegraphics[trim=5 0 0 5,clip,width=0.47\textwidth]{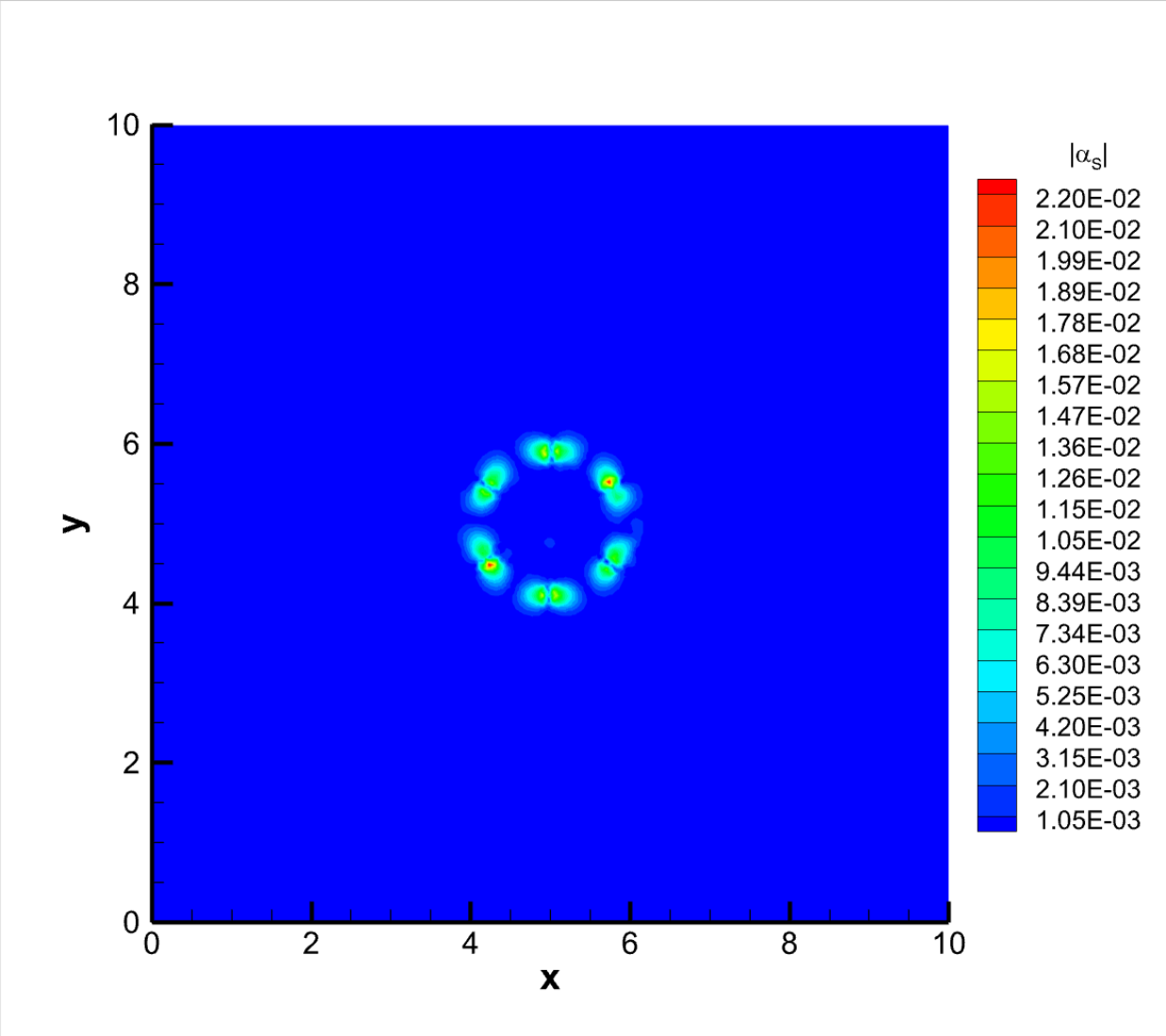} \\
			\includegraphics[trim=0 0 0 5,width=0.47\textwidth]{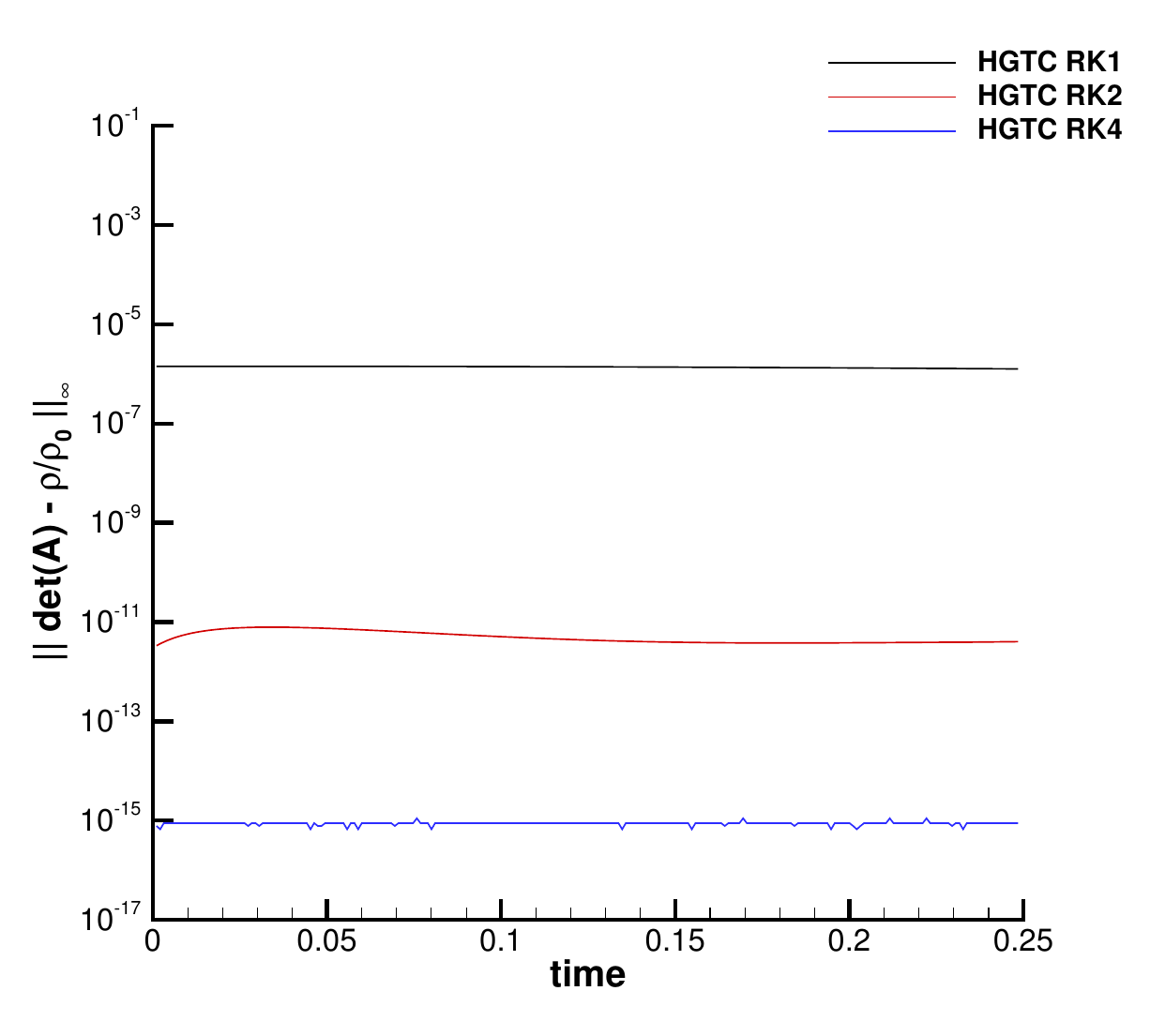} &
			\includegraphics[trim=0 0 0 5,width=0.47\textwidth]{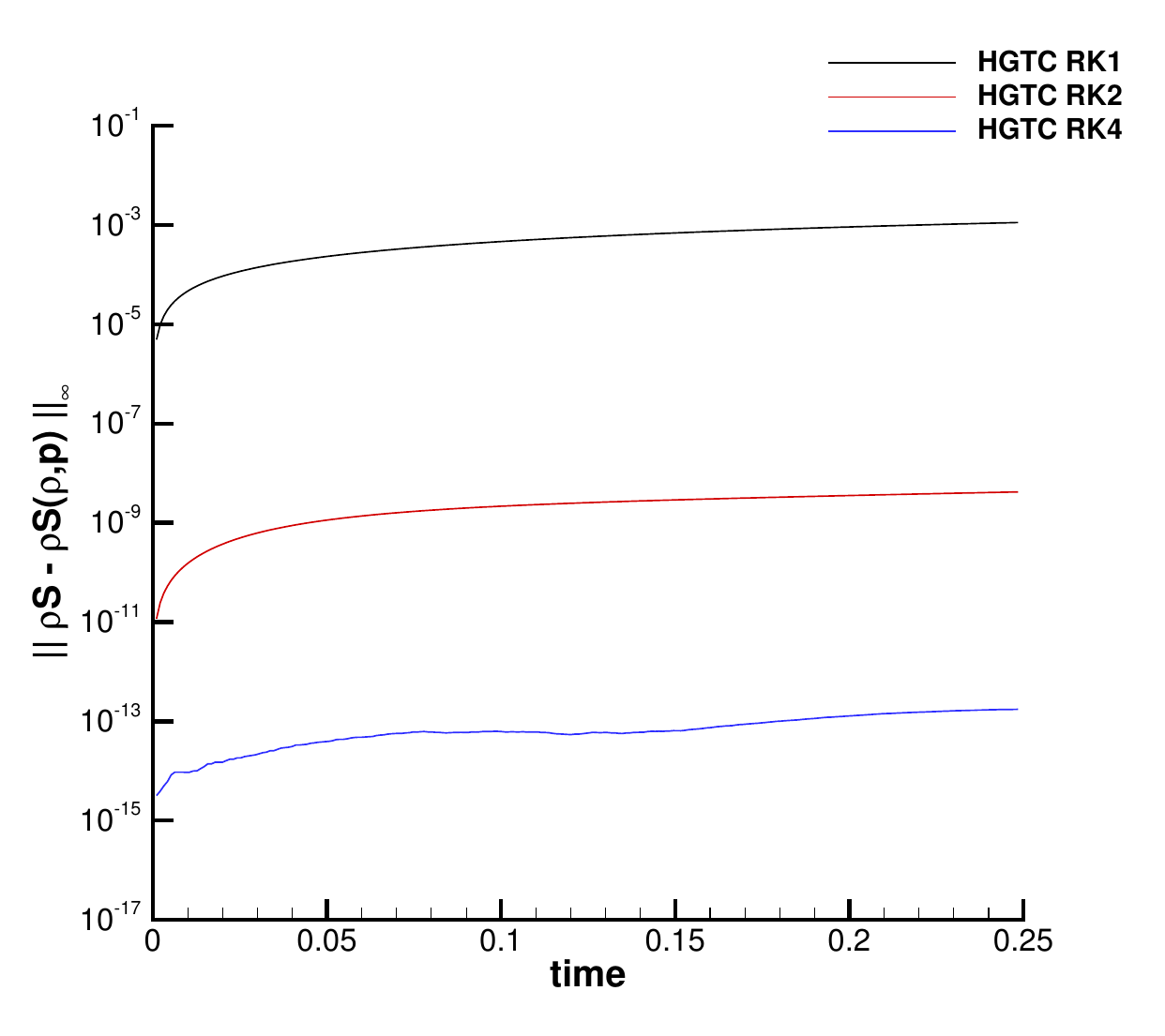} \\
		\end{tabular} 
		\caption{Isentropic vortex problem at time $t_f=0.25$. Top: map of the geometric correction factor $|\alpha_\A|$ (left) and of the thermodynamic correction factor $|\alpha_S|$ (right) with mesh size $h=1/6$. Bottom: time evolution of the mass (left) and entropy (right) conservation errors for Runge-Kutta time integration schemes of order 1 (black line), 2 (red line) and 4 (blue line).}
		\label{fig.shuvortex}
	\end{center}
\end{figure}

We also use this test case to analyze the time convergence which ultimately affects the preservation of the determinant and the entropy compatibility. Therefore we measure the errors of the total mass and entropy conservation according to \eqref{eqn.epsA_S} while running this simulation until the time $t_f=1$ on one single unstructured mesh with characteristic size of $h=1/3$. Three different Runge-Kutta time integrators are used of order $N=\{1,2,4\}$ (see \ref{app.rk}), and the results are collected in Table \ref{tab.RKconv}. We observe that the convergence rates for the entropy conservation exhibits order of accuracy $\mathcal{O}(N+1)$, and convergence of order $\mathcal{O}(N+2)$ is achieved for the total mass conservation. 
\begin{table}[!htbp]  
	\caption{Time convergence study related to total mass and entropy conservation for the isentropic vortex problem at time $t_{f}=1$ with three different Runge-Kutta time integration schemes on a mesh with size $h=1/3$. The errors are measured in the $L_{\infty}$ norm and refer to the geometric thermodynamic errors given by \eqref{eqn.epsA_S}.}  
	\begin{center} 
		\begin{small}
			\renewcommand{\arraystretch}{1.2}
			\begin{tabular}{ccccc} 
				\multicolumn{5}{c}{Runge-Kutta $\mathcal{O}(1)$} \\
				$\dt$  & $ \delta_\A$ &    & $\delta_S$ &  \\
				\hline
				8.00E-03 & 5.1415E-05 & -    & 2.2516E-02 & -    \\
				4.00E-03 & 1.2874E-05 & 2.00 & 1.1238E-02 & 1.00 \\
				2.00E-03 & 3.2211E-06 & 2.00 & 5.6138E-03 & 1.00 \\
				\multicolumn{5}{c}{} \\				
				
				\multicolumn{5}{c}{Runge-Kutta $\mathcal{O}(2)$} \\
				$\dt$  & $ \delta_\A$ &    & $\delta_S$ &   \\
				\hline
				8.00E-03 & 3.3043E-09 & -    & 1.3195E-06 & -    \\
				4.00E-03 & 3.3385E-10 & 3.31 & 2.5736E-07 & 2.36 \\
				2.00E-03 & 3.6793E-11 & 3.18 & 5.5271E-08 & 2.22 \\
				\multicolumn{5}{c}{} \\				
				
				\multicolumn{5}{c}{Runge-Kutta $\mathcal{O}(4)$} \\
				$\dt$ & $ \delta_\A$ &    & $\delta_S$ &   \\
				\hline
				8.00E-03 & 2.9510E-13 & -    & 6.6099E-11 & -    \\
				4.00E-03 & 9.1038E-15 & 5.02 & 4.1208E-12 & 4.00 \\
				2.00E-03 & 2.8818E-16 & 4.98 & 2.5709E-13 & 4.00 \\
			\end{tabular}
		\end{small}
	\end{center}
	\label{tab.RKconv}
\end{table}
The time evolution of the mass and entropy conservation errors are plot in Figure \ref{fig.shuvortex}, where we also show the map of the scalar correction factors $\alpha_\A$ and $\alpha_S$ at the final time. 

\subsection{Riemann problems}
The novel HGTC scheme is here validated against three one-dimensional Riemann problems taken from \cite{ToroBook,HTCAbgrall}. The computational domain is the rectangular box $\Omega=[-0.5;0.5]\times[-0.05;0.05]$ with periodic boundaries in the $y-$direction and transmissive boundaries along the $x-$direction. The computational mesh is unstructured made of polygons and it has a characteristic size of $h=1/4096$ and all the simulations are run in 2D, thus the properties of symmetry preservation of the numerical solution are verified as well. Indeed, despite the one-dimensional setting of the Riemann problems, these test cases become fully multidimensional in the case of unstructured Voronoi meshes, where no mesh edges are in principle aligned with the flow. The initial condition is given in terms of a left and a right state separated at position $x=x_d$. Table \ref{tab.RPinit} summarizes the setup of the three Riemann problems considered here.

\begin{table}[!htbp]  
	\caption{Initialization of Riemann problems. Initial states left (L) and right (R) are reported as well as the final time of the simulation $t_f$ and the position of the initial discontinuity $x_d$. }  
	\begin{center} 
		\begin{small}
			\renewcommand{\arraystretch}{1.0}
			\begin{tabular}{l|cc|cccc|cccc} 
				Name & $t_{f}$ & $x_d$ & $\rho_L$ & $v_{1,L}$ & $v_{2,L}$ & $p_L$ & $\rho_R$ & $v_{1,R}$ & $v_{2,R}$ & $p_R$ \\
				\hline
				RP1 & 0.035 & -0.2 & 5.99924 & 19.5975 & 0.0 & 460.894 & 5.99924 & -6.19633 & 0.0 & 46.095 \\
				RP2 & 0.15 & 0.0 & 1.0 & -2.0 & 0.0 & 0.4 & 1.0 & 2.0 & 0.0 & 0.4 \\
				RP3 & 0.20 & 0.0 & 1.0 & 0.0 & -0.2 & 1.0 & 0.5 & 0.0 & 0.2 & 0.5 \\
			\end{tabular}
		\end{small}
	\end{center}
	\label{tab.RPinit}
\end{table}

\begin{figure}[!htbp]
	\begin{center}
		\begin{tabular}{ccc} 
			\includegraphics[width=0.33\textwidth]{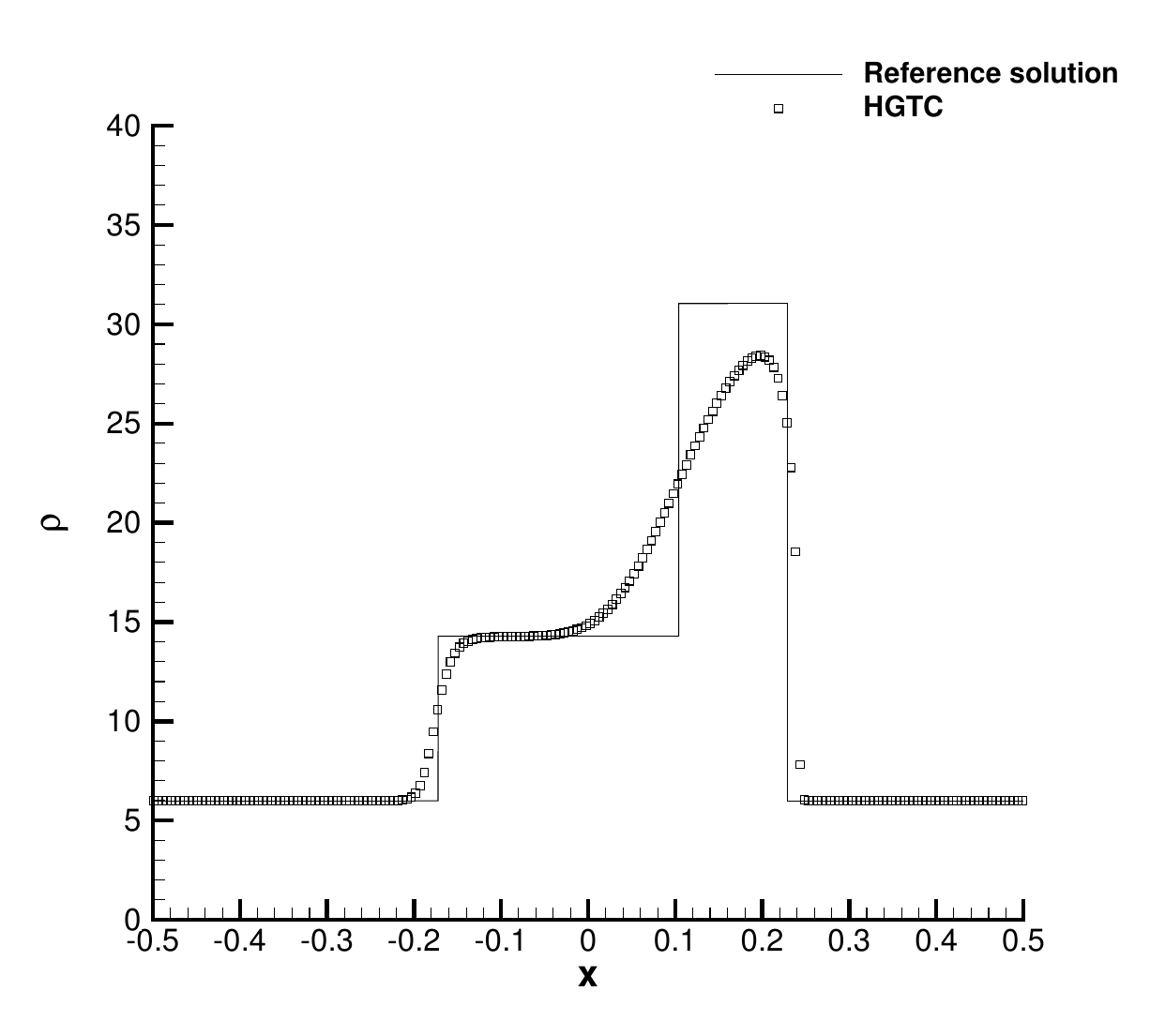}  & 
			\includegraphics[width=0.33\textwidth]{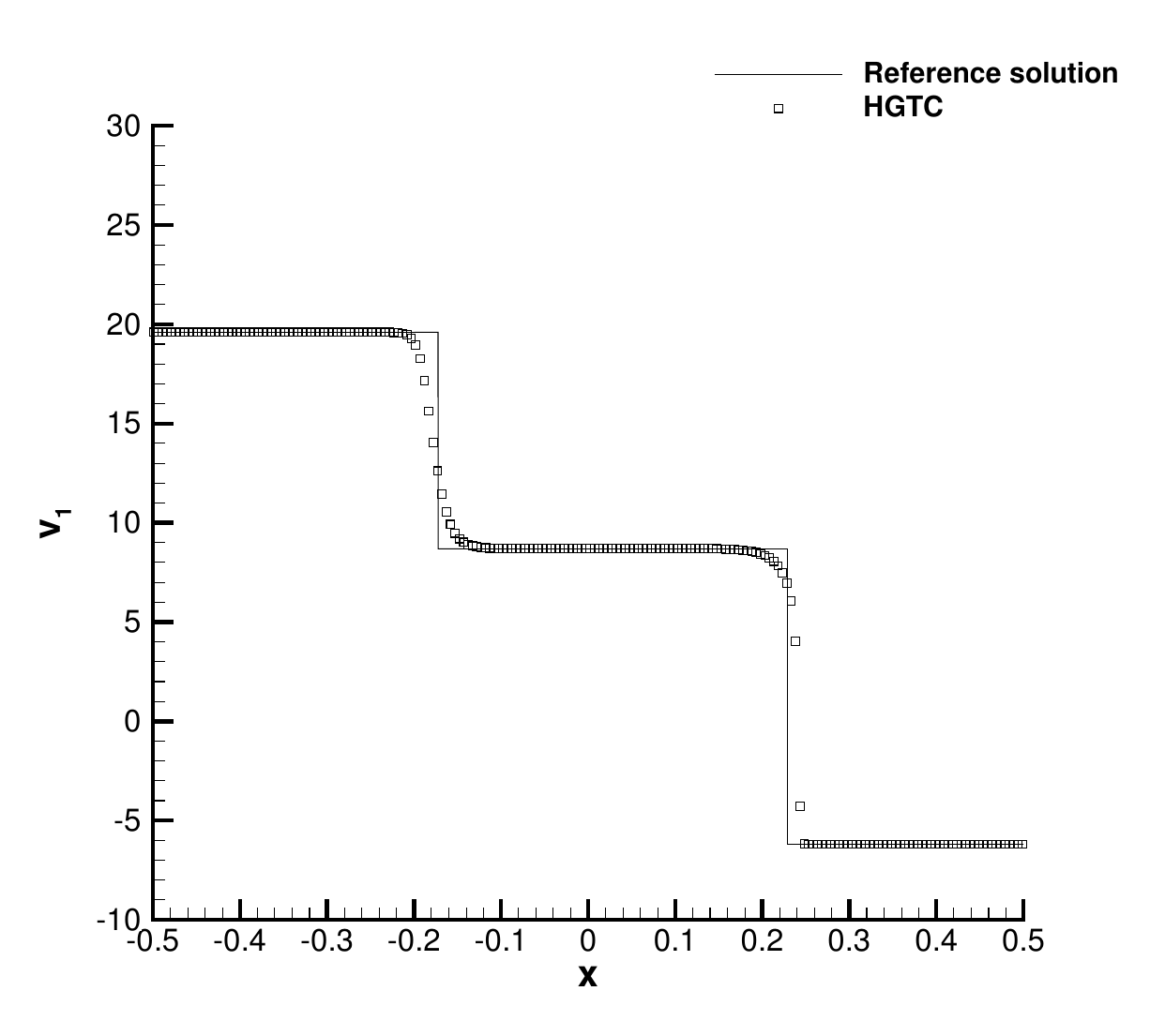} &
    	    \includegraphics[width=0.33\textwidth]{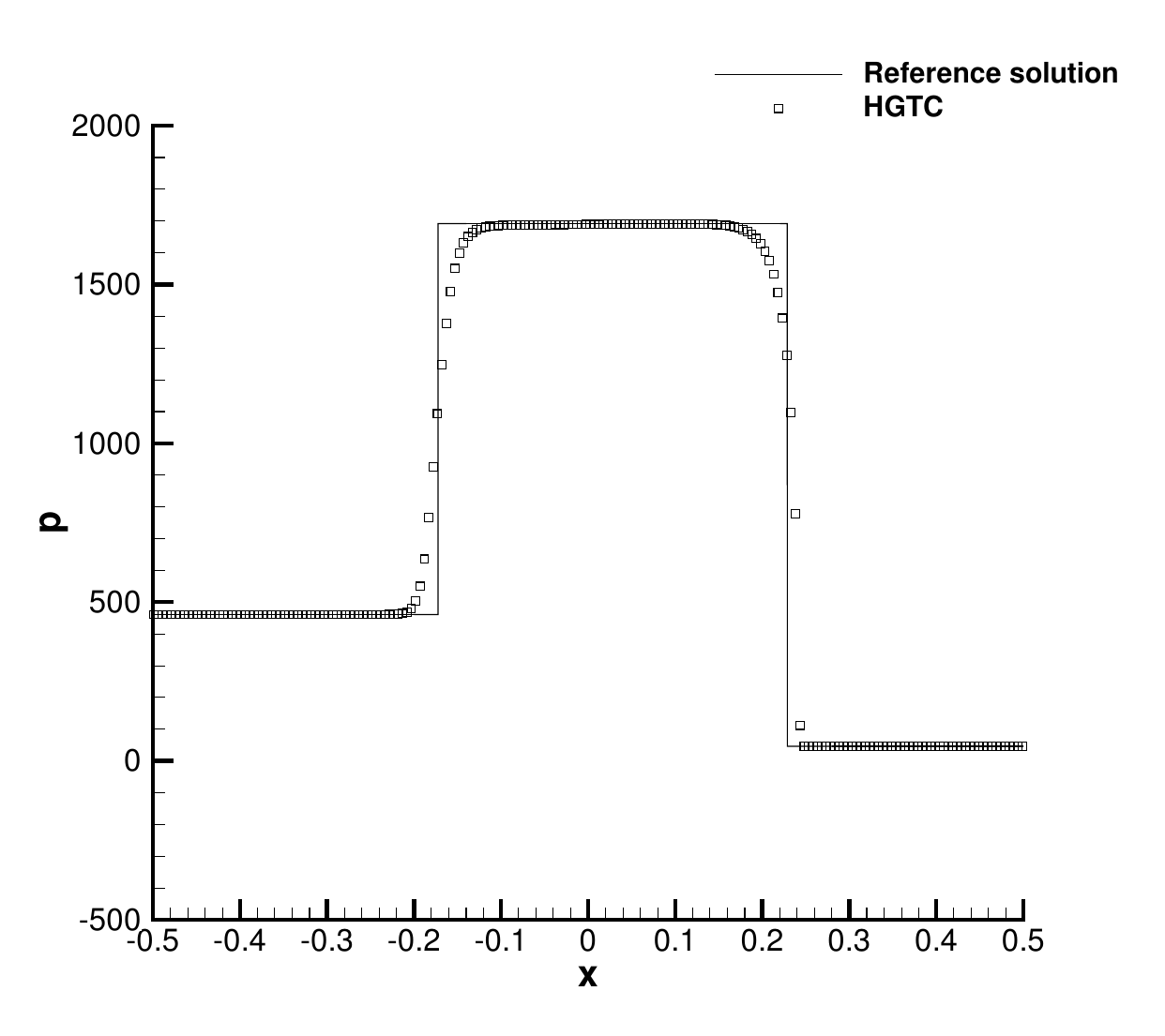}  \\ 
		\end{tabular} 
		\caption{Riemann problem RP1 at final time $t_f=0.035$. Comparison of density, horizontal velocity and pressure against the reference solution extracted with a one-dimensional cut of 200 equidistant points along the $x$-direction at $y=0$.}
		\label{fig.RP2}
	\end{center}
\end{figure}

\begin{figure}[!htbp]
	\begin{center}
		\begin{tabular}{ccc} 
			\includegraphics[width=0.33\textwidth]{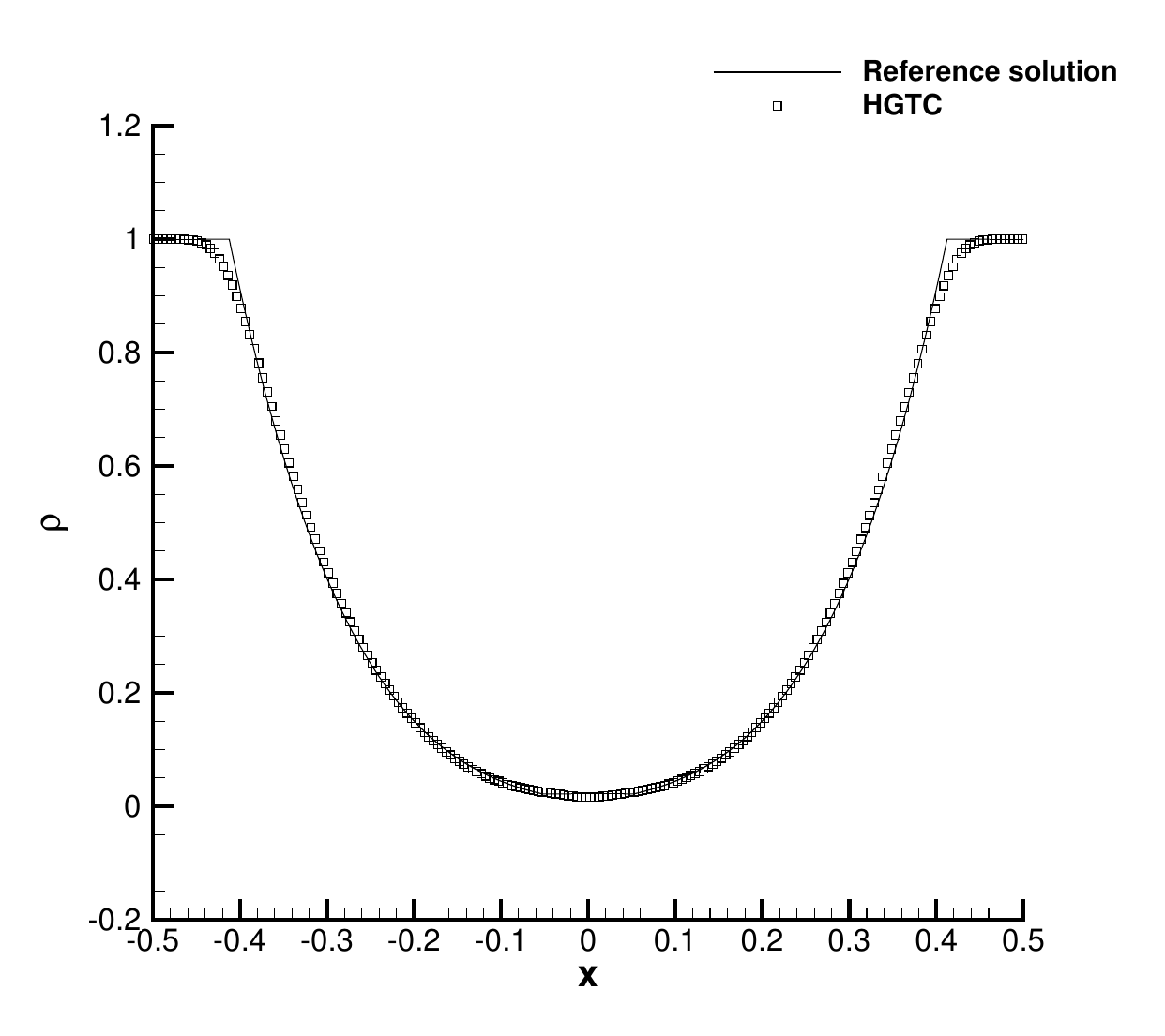}  & 
			\includegraphics[width=0.33\textwidth]{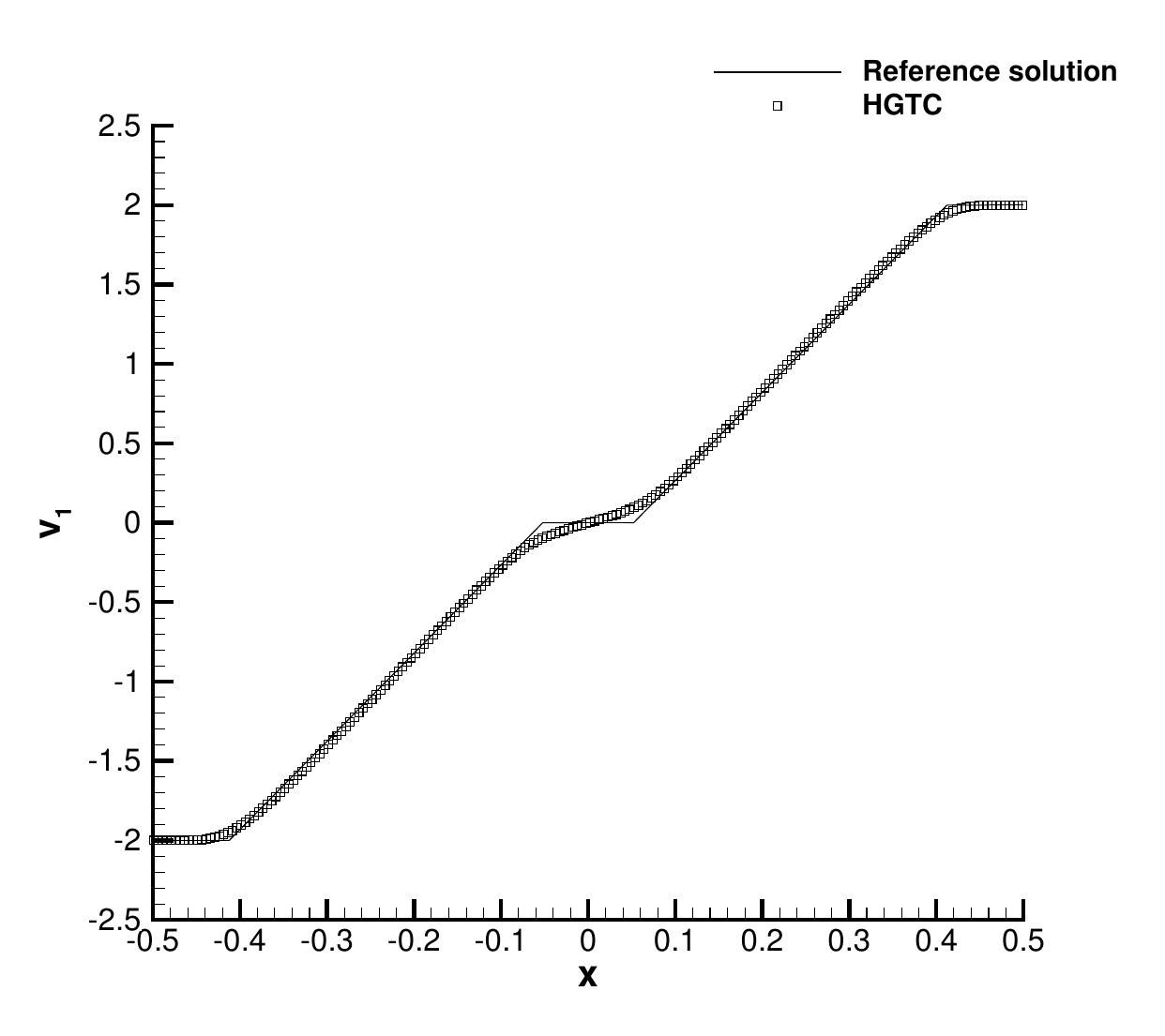} &
			\includegraphics[width=0.33\textwidth]{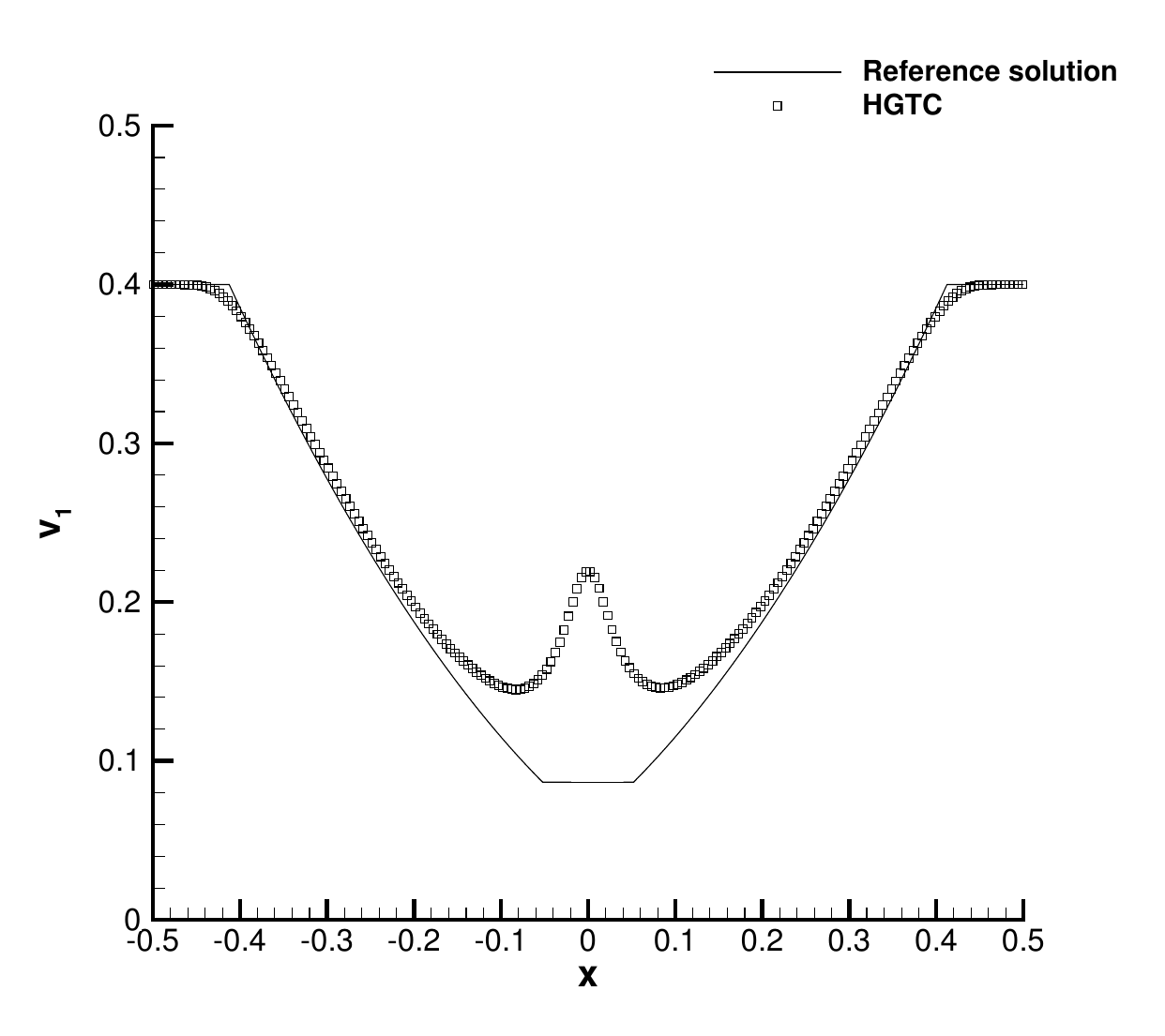}  \\ 
		\end{tabular} 
		\caption{Riemann problem RP2 at final time $t_f=0.15$. Comparison of density, horizontal velocity and temperature against the reference solution extracted with a one-dimensional cut of 200 equidistant points along the $x$-direction at $y=0$.}
		\label{fig.RP123}
	\end{center}
\end{figure}

\begin{figure}[!htbp]
	\begin{center}
		\begin{tabular}{ccc} 
			\includegraphics[width=0.33\textwidth]{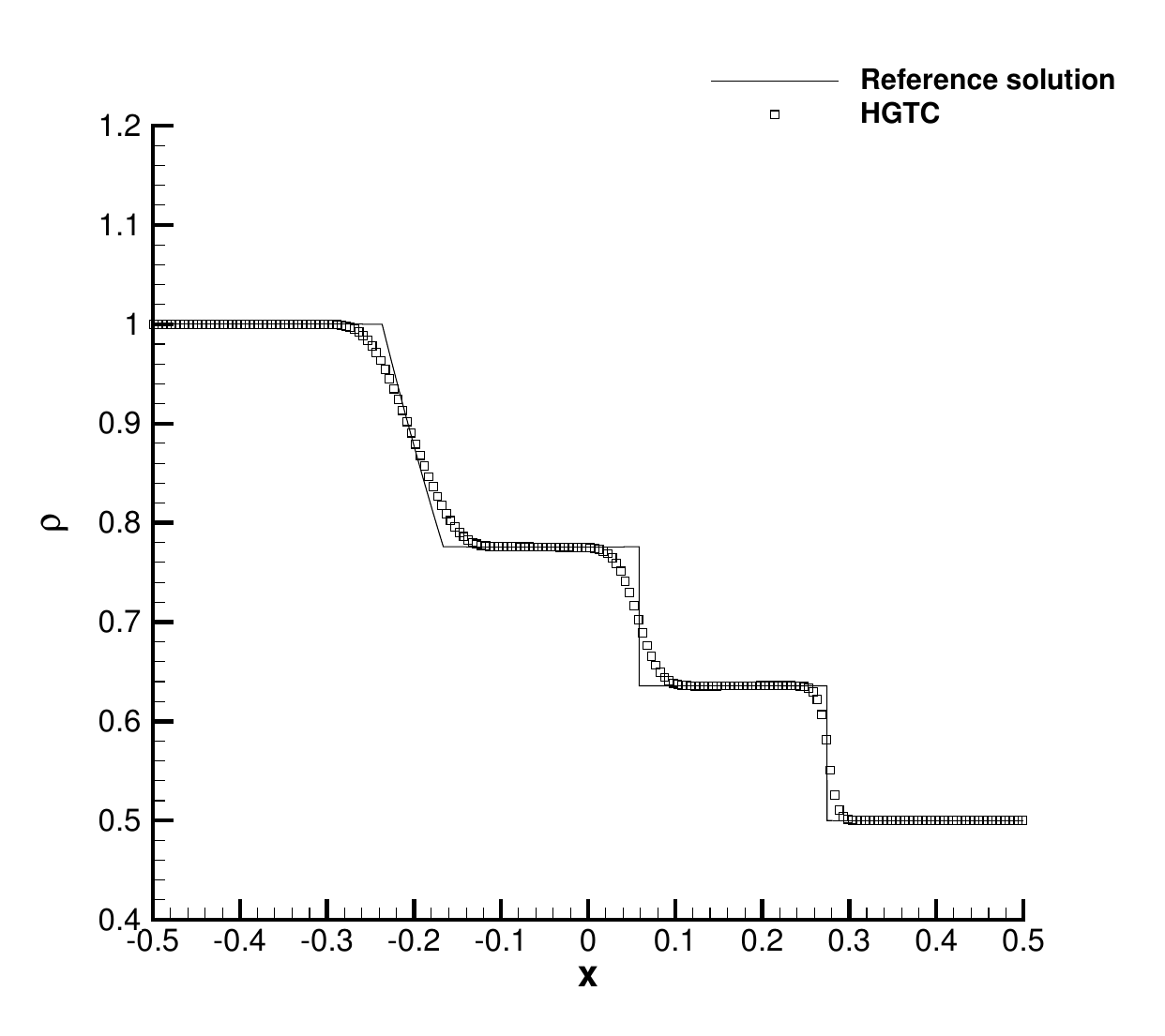}  & 
			\includegraphics[width=0.33\textwidth]{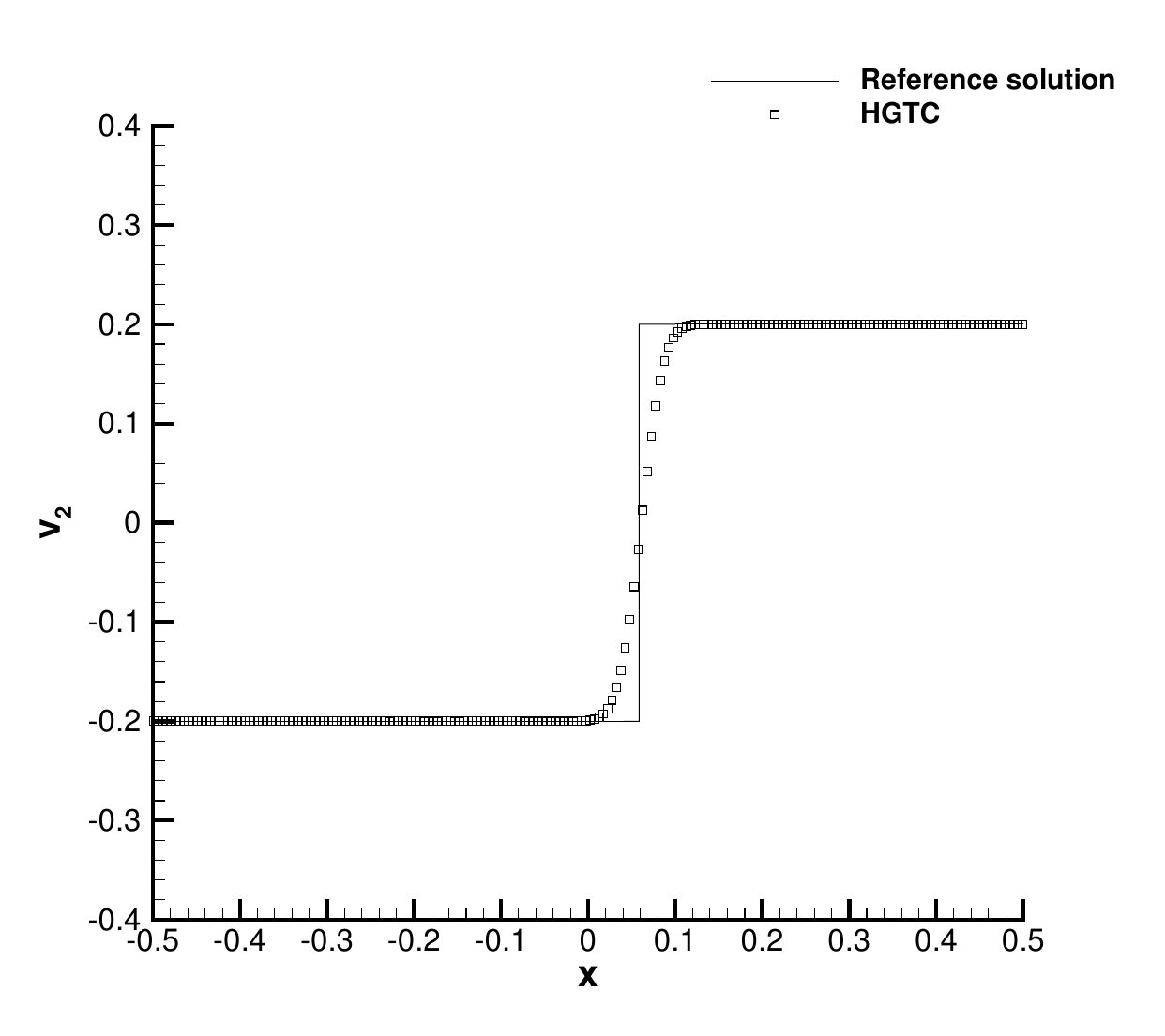} &
			\includegraphics[width=0.33\textwidth]{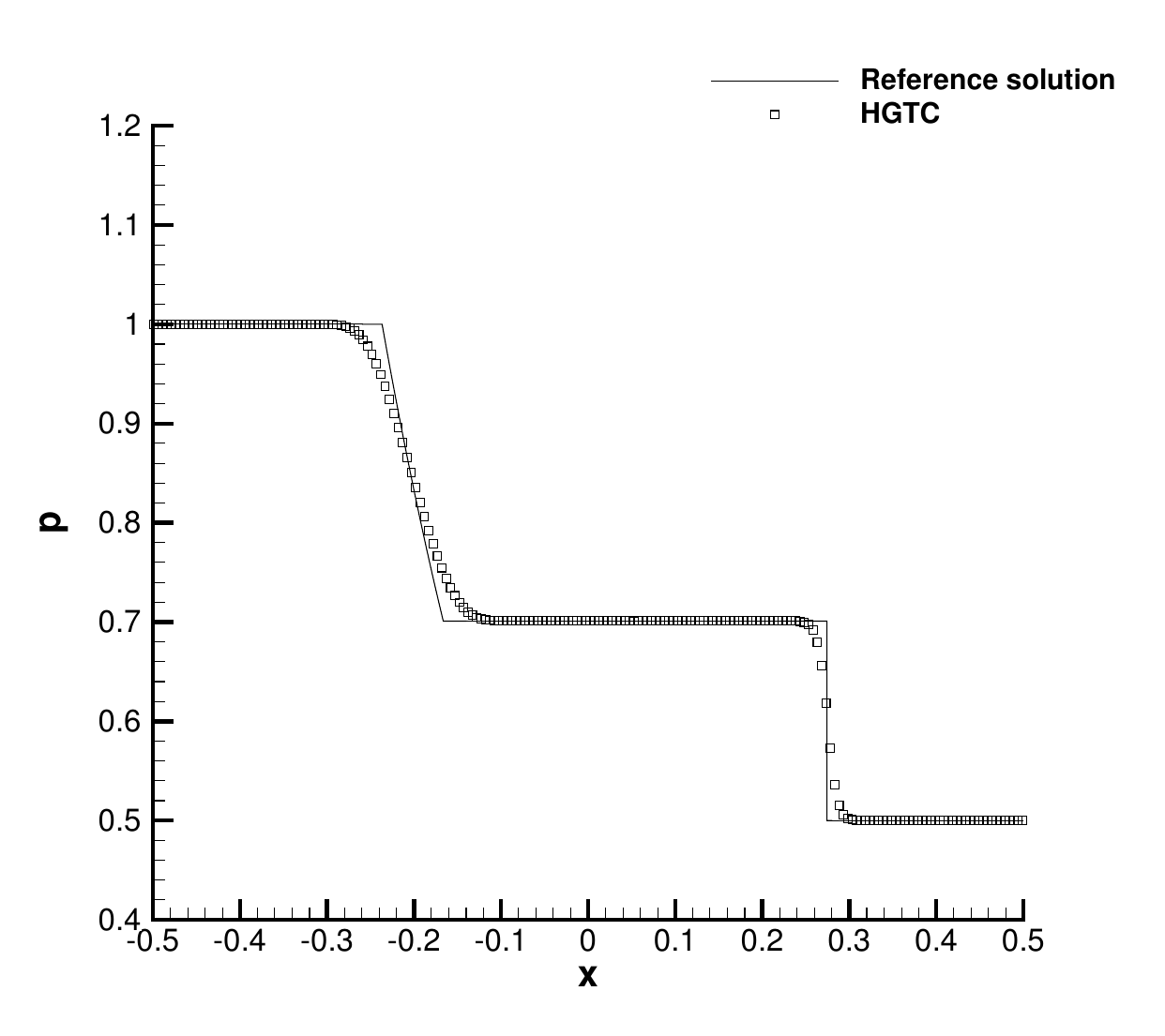}  \\ 
		\end{tabular} 
		\caption{Riemann problem RP3 at final time $t_f=0.2$. Comparison of density, vertical velocity and pressure against the reference solution for the compressible Euler equations extracted with a one-dimensional cut of 200 equidistant points along the $x$-direction at $y=0$.}
		\label{fig.RP3}
	\end{center}
\end{figure}

\begin{figure}[!htbp]
	\begin{center}
		\begin{tabular}{ccc} 
			\includegraphics[width=0.33\textwidth]{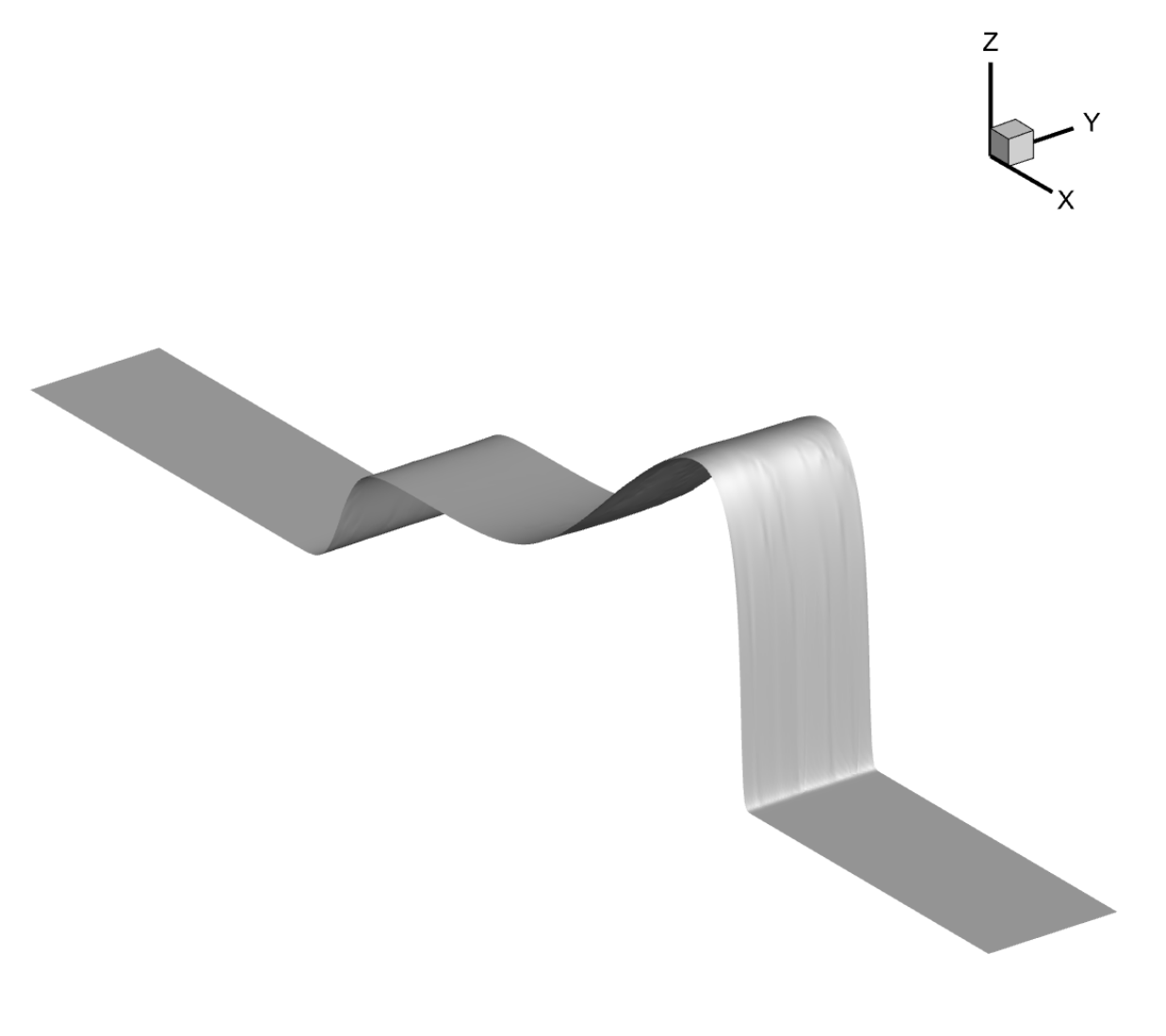}  & 
			\includegraphics[width=0.33\textwidth]{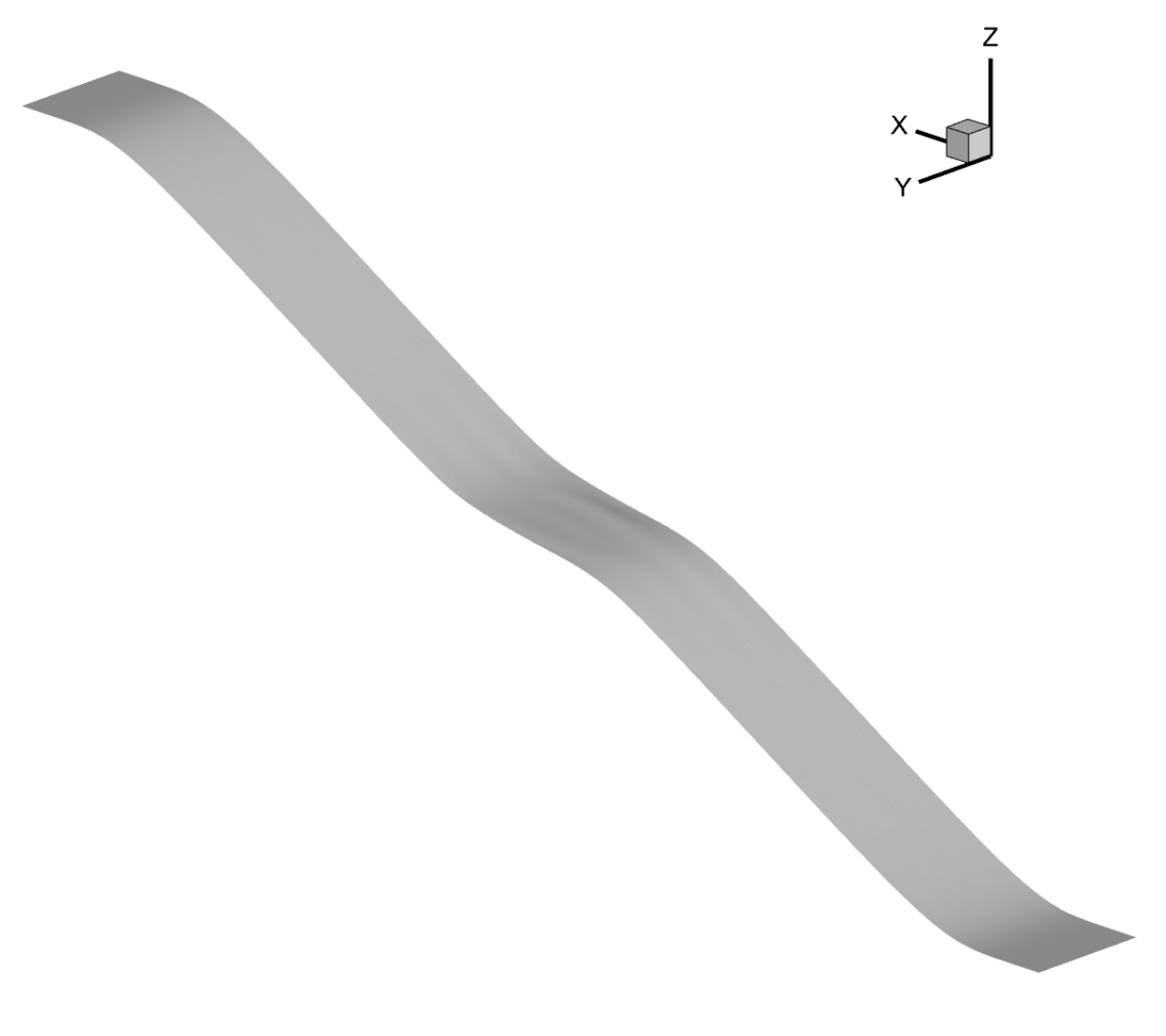} &
			\includegraphics[width=0.33\textwidth]{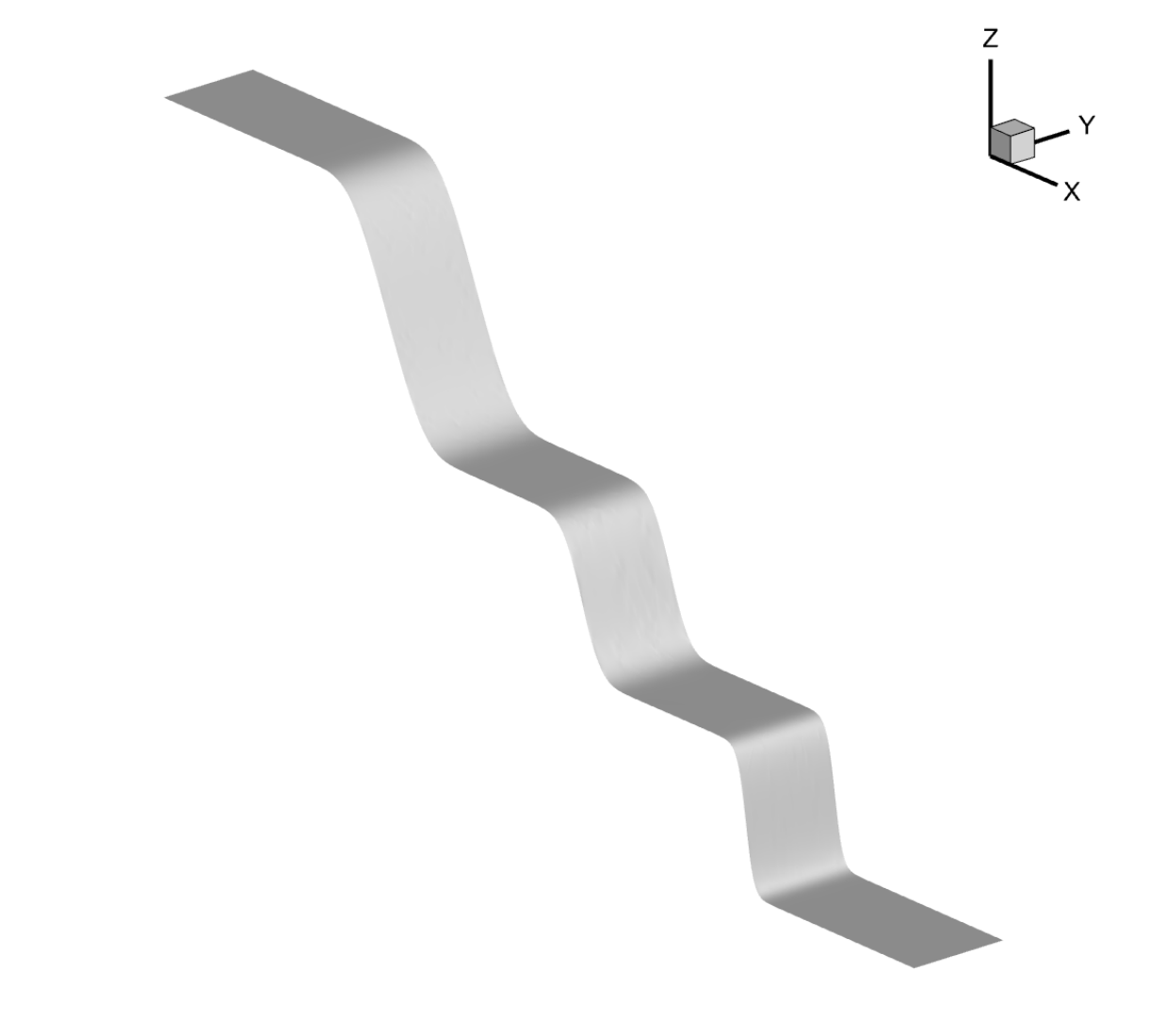}  \\ 
		\end{tabular} 
		\caption{Three-dimensional view of density for RP1 (left), horizontal velocity for RP2 (middle) and density for RP3 (right) at their corresponding final times.}
		\label{fig.RP3D}
	\end{center}
\end{figure}

The first two Riemann problems RP1 and RP2 involve the Euler equations for compressible gas dynamics (i.e. $c_s=c_h=0$), and the reference solution is computed with the exact Riemann solver detailed in \cite{ToroBook}. The last Riemann problem is concerned with the full model \eqref{eqn.GPR} and we set $\mu=\kappa=10^{-5}$, so that the stiff relaxation limit of the model is retrieved and numerically assessed. The reference solution is obtained numerically using a second order TVD finite volume method on a very fine mesh of 100000 control volumes. The results are collected in Figures \ref{fig.RP2}-\ref{fig.RP3}, showing a good agreement with the reference solution in all cases. 
To appreciate that the one-dimensional symmetry of the solution is well preserved, we show in Figure \ref{fig.RP3D} a three-dimensional view of the solution for the three Riemann problems considered here.

\subsection{Circular explosion problem}
\begin{figure}[!htbp]
	\begin{center}
		\begin{tabular}{cc} 
			\includegraphics[width=0.47\textwidth]{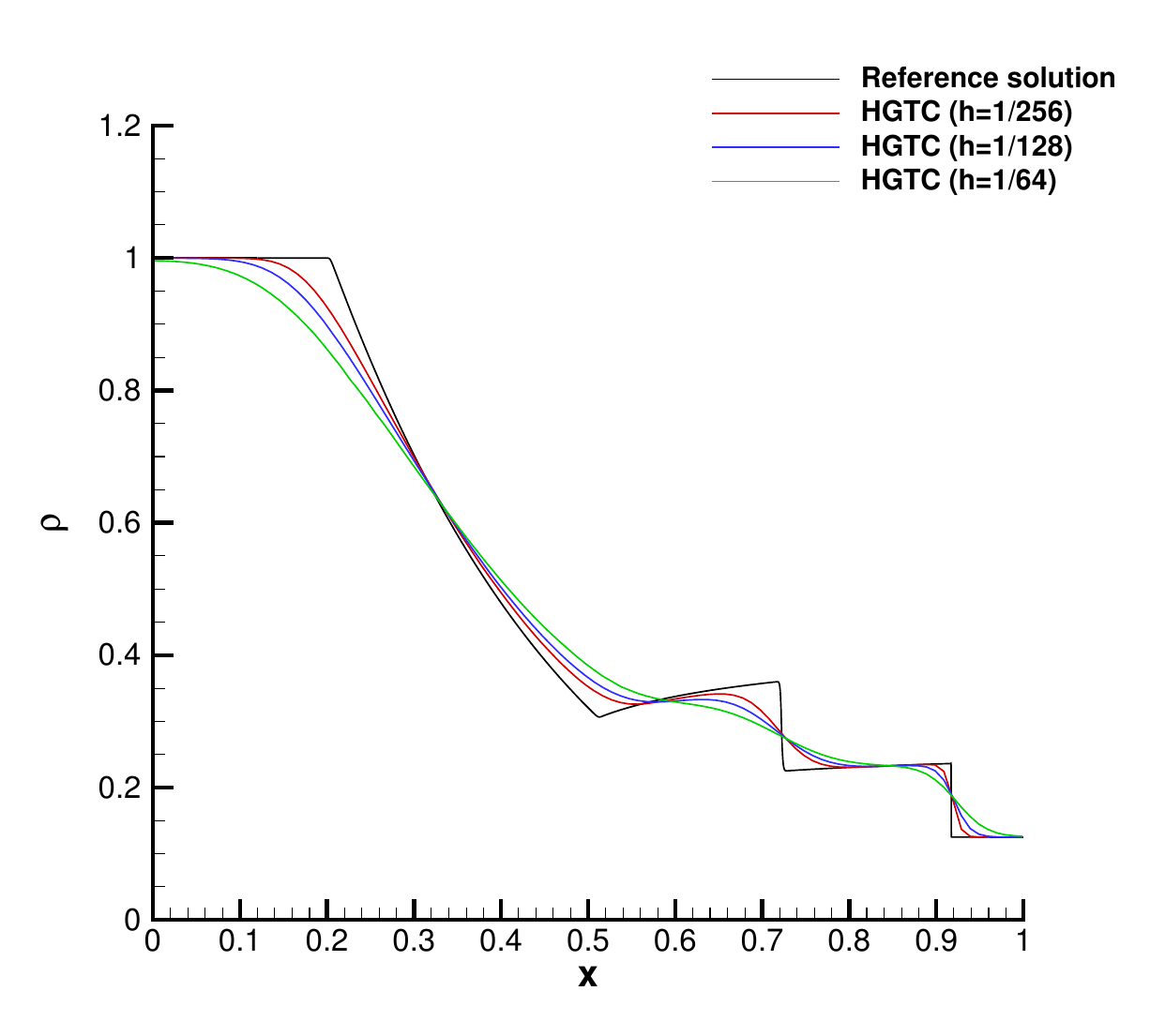} & 
			\includegraphics[width=0.47\textwidth]{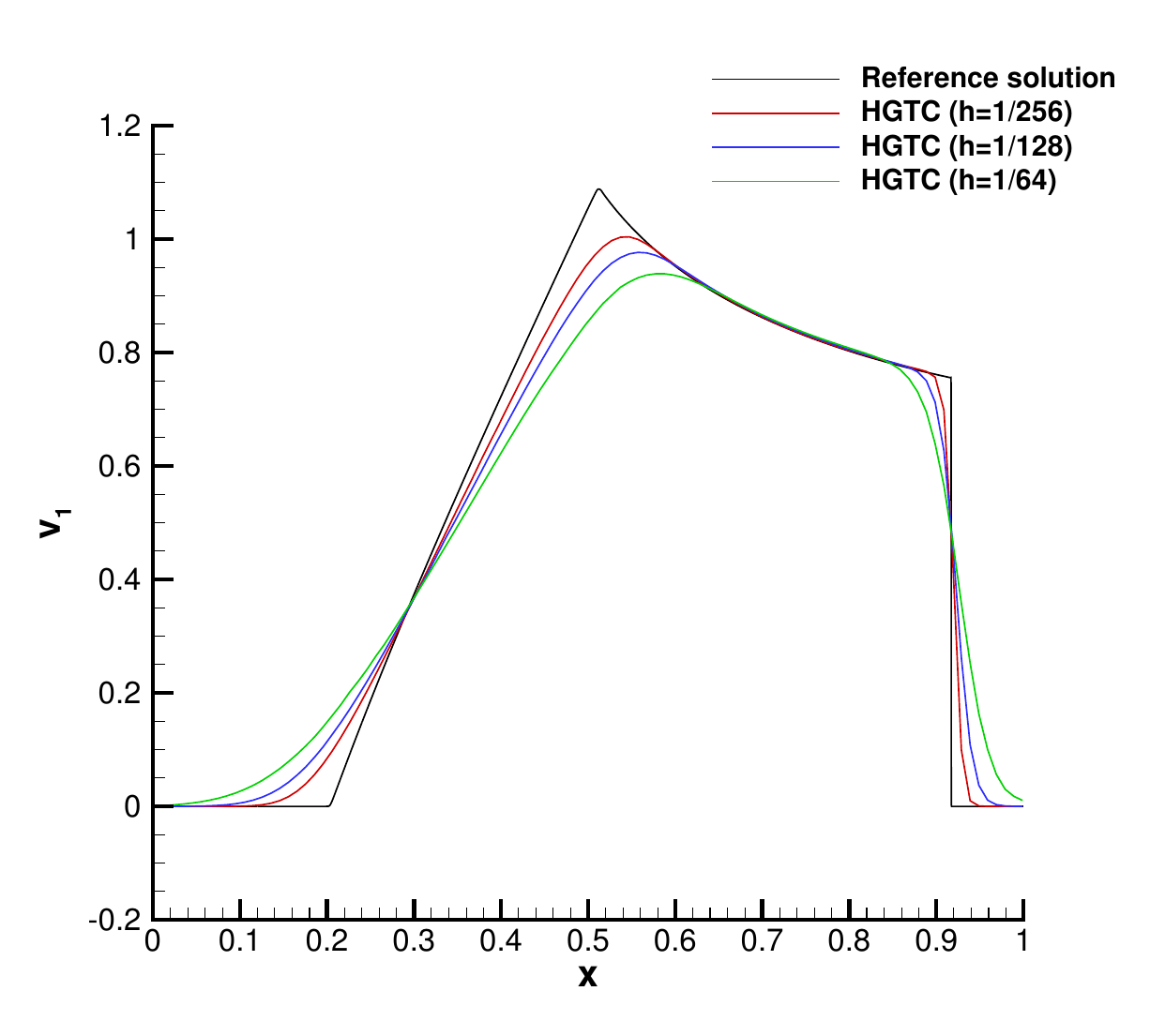} \\   
			\includegraphics[width=0.47\textwidth]{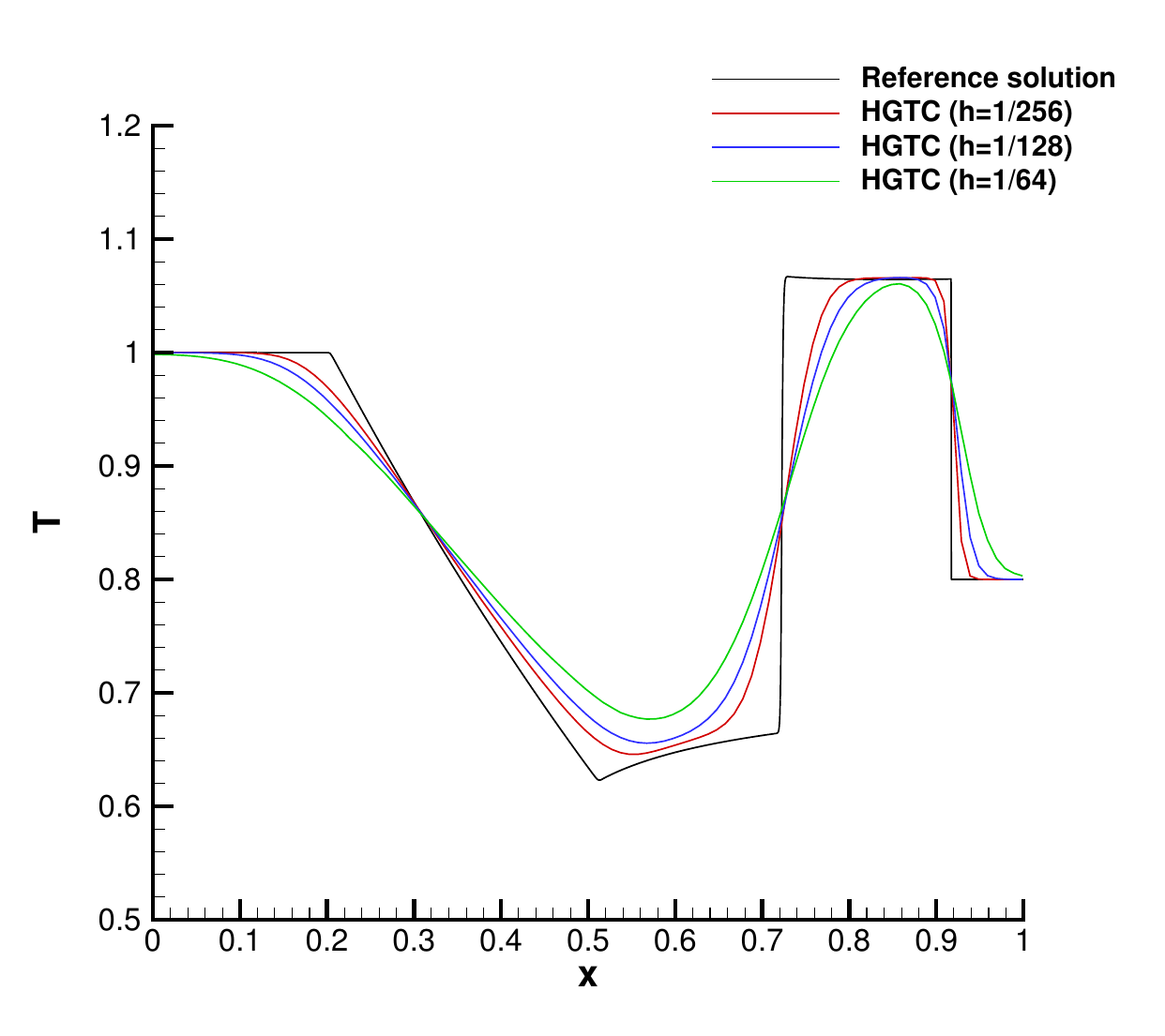} &
			\includegraphics[width=0.47\textwidth]{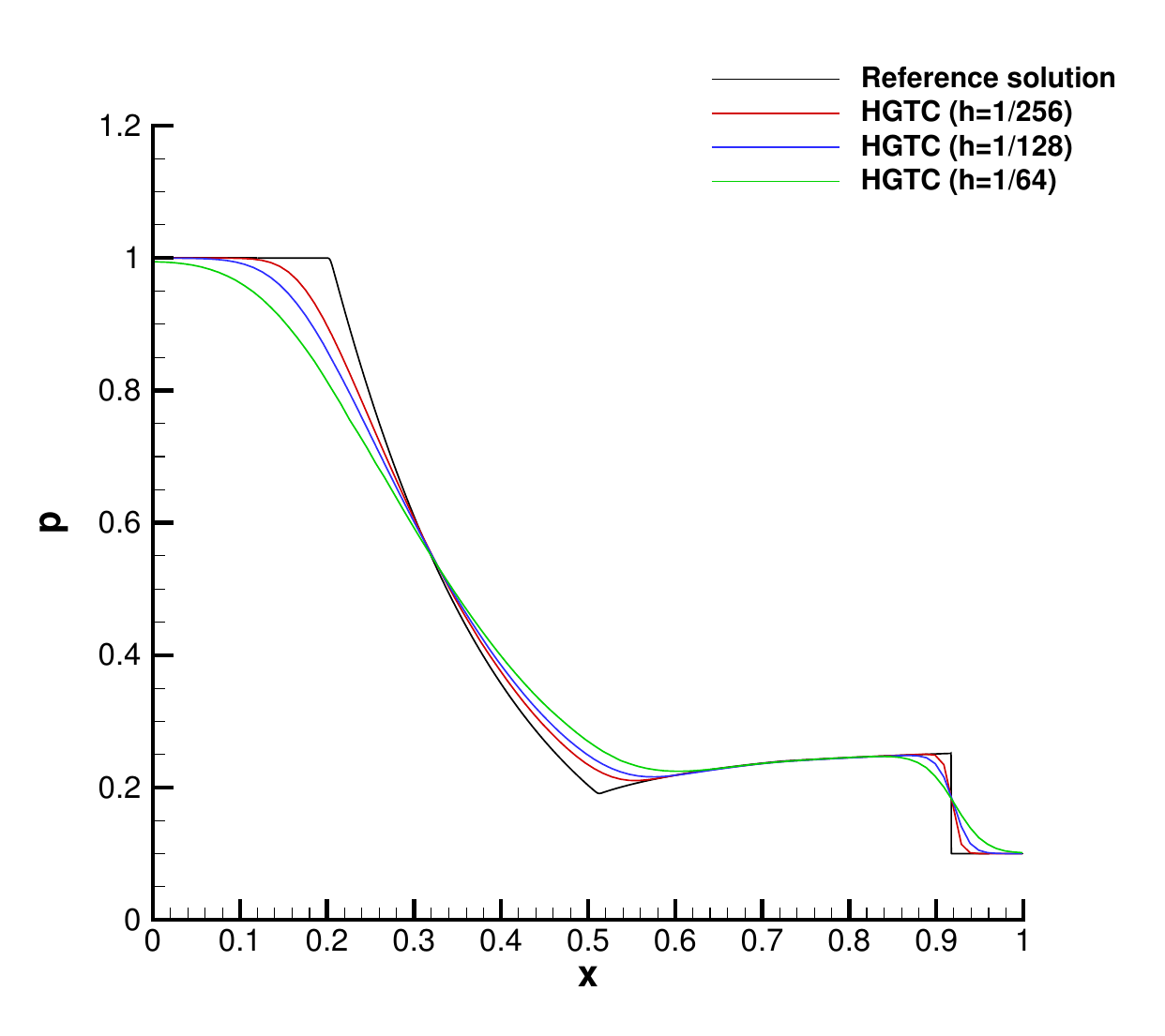} \\
		\end{tabular} 
		\caption{Explosion problem at time $t_f=0.25$. Numerical results for density, horizontal velocity, pressure and temperature (from top left to bottom right) compared against the reference solution extracted with a one-dimensional cut of 200 equidistant points along the $x$-direction at $y=0$. Mesh convergence analysis with characteristic mesh size $h=1/256$ (red line), $h=1/128$ (red line) and $h=1/64$ (green line).}
		\label{fig.ep2d}
	\end{center}
\end{figure}

\begin{figure}[!htbp]
	\begin{center}
		\begin{tabular}{cc} 
			\includegraphics[width=0.47\textwidth]{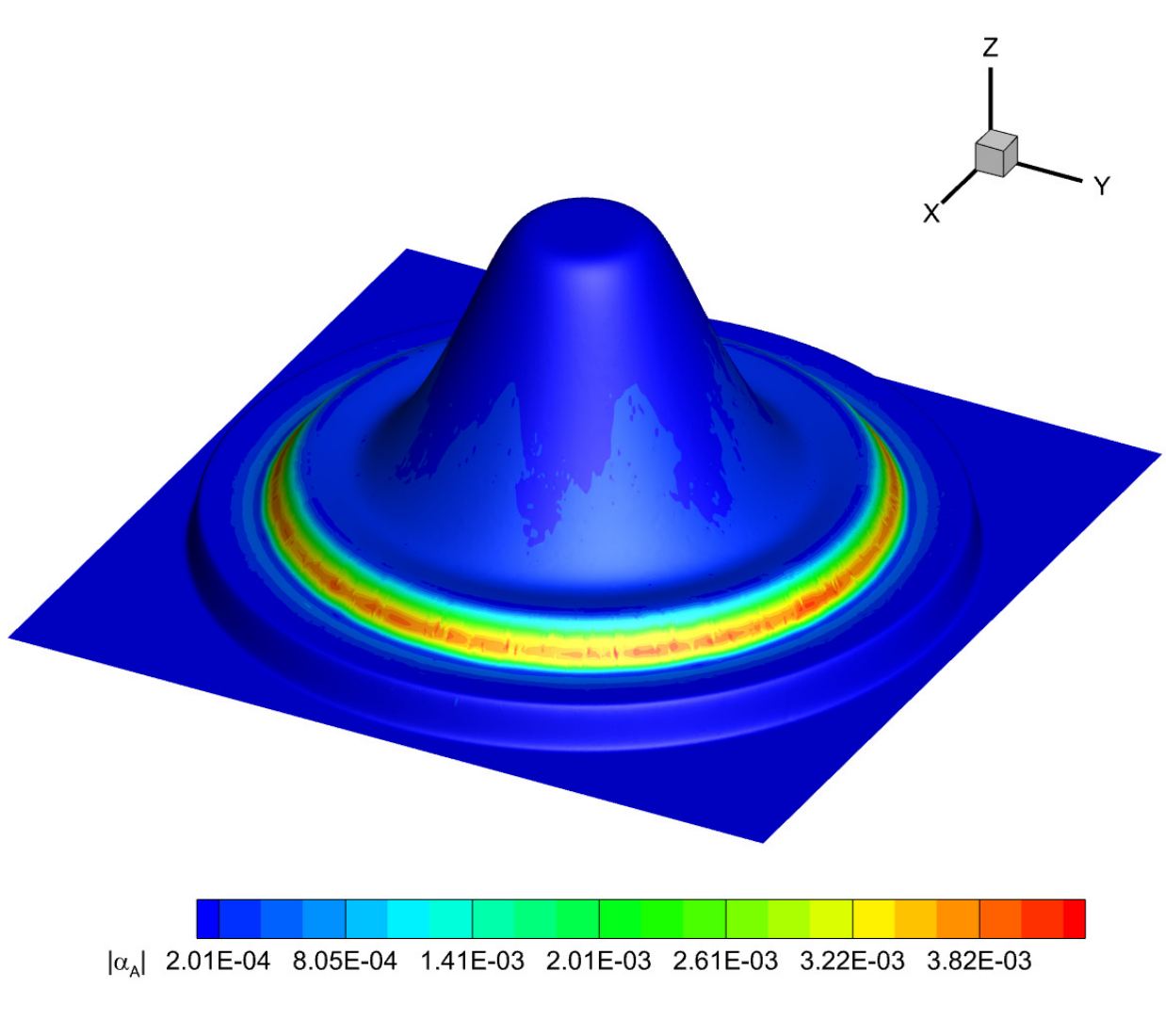} & 
			\includegraphics[width=0.47\textwidth]{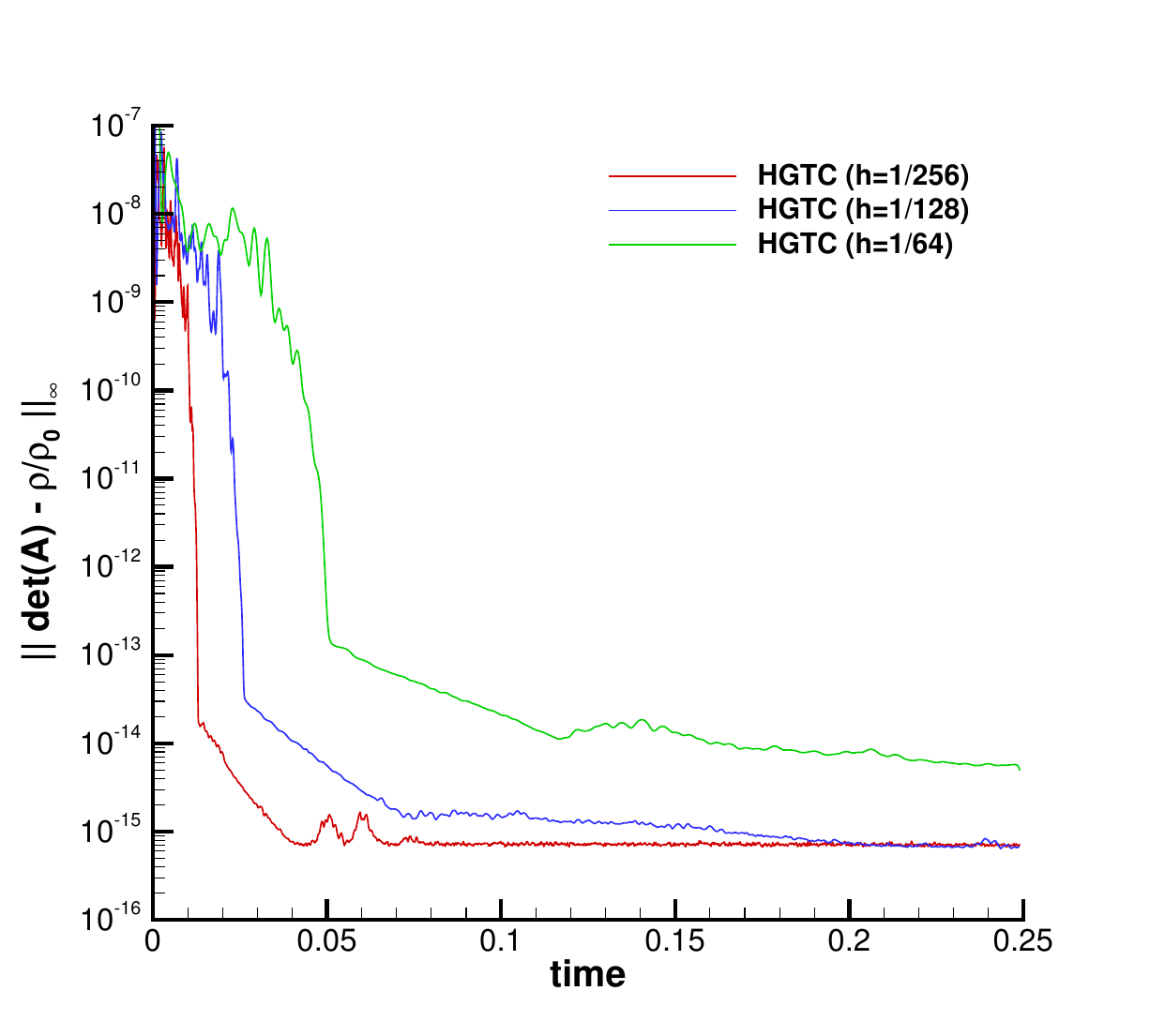} \\
		\end{tabular} 
		\caption{Explosion problem at time $t_f=0.25$. Left: three-dimensional view of the density distribution with the map of the correction factor $\alpha_\A$ for the results obtained with $h=1/256$. Right: time evolution of the mass conservation errors for $h=1/256$ (red line), $h=1/128$ (blue line) and $h=1/64$ (green line).}
		\label{fig.ep2d_detA}
	\end{center}
\end{figure}

We consider a cylindrical explosion problem to test the HGTC schemes with numerical dissipation, which is here activated since the solution exhibits an outward traveling shock wave. The computational domain is given by $\Omega=[-1;1]^2$ with transmissive boundaries, and the fluid is initially assigned as follows:
\begin{equation}
	(\rho,v_1,v_2,v_3,p)=\left\{
	\begin{array}{lll}
		(1, \, 0,\, 0,\, 0,\, 1) & & r<R \\
		(0.125,\, 0,\, 0,\, 0,\, 0.1) & & r\geq R
	\end{array} \right., \qquad t=0, \quad \xx \in \Omega,
\end{equation}
where $R=0.5$ denotes the radius of the initial discontinuity and $r=\sqrt{x_1^2+x_2^2}$ represents the generic radial coordinate. An inviscid fluid is considered by setting $c_s=c_h=0$ and the final time of the simulation is chosen to be $t_f=0.25$. We run this test on three different Voronoi meshes with characteristic mesh size of $h=1/256$, $h=1/128$ and $h=1/64$. The numerical results are compared against the reference solution that has been computed by solving the compressible Euler equations with geometric sources \cite{ToroBook} employing a classical second order TVD finite volume scheme on a very fine mesh composed of 20000 cells. An overall very good agreement can be observed in Figure \ref{fig.ep2d}, that numerically confirms the convergence of the new HGTC schemes as the mesh resolution gets finer. 
Figure \ref{fig.ep2d_detA} depicts a map of the correction factor $\alpha_\A$ as well as the time evolution of the total mass conservation errors for all the simulations. The solution preserves an excellent cylindrical symmetry despite the unstructured nature of the mesh, as shown by the three-dimensional density distribution plot in Figure \ref{fig.ep2d_detA}. Furthermore, we notice that the highest correction for the preservation of the determinant constraint occurs across the contact wave.

\subsection{Viscous shock profile}
\begin{figure}[!htbp]
	\begin{center}
		\begin{tabular}{cc} 
			\includegraphics[width=0.47\textwidth]{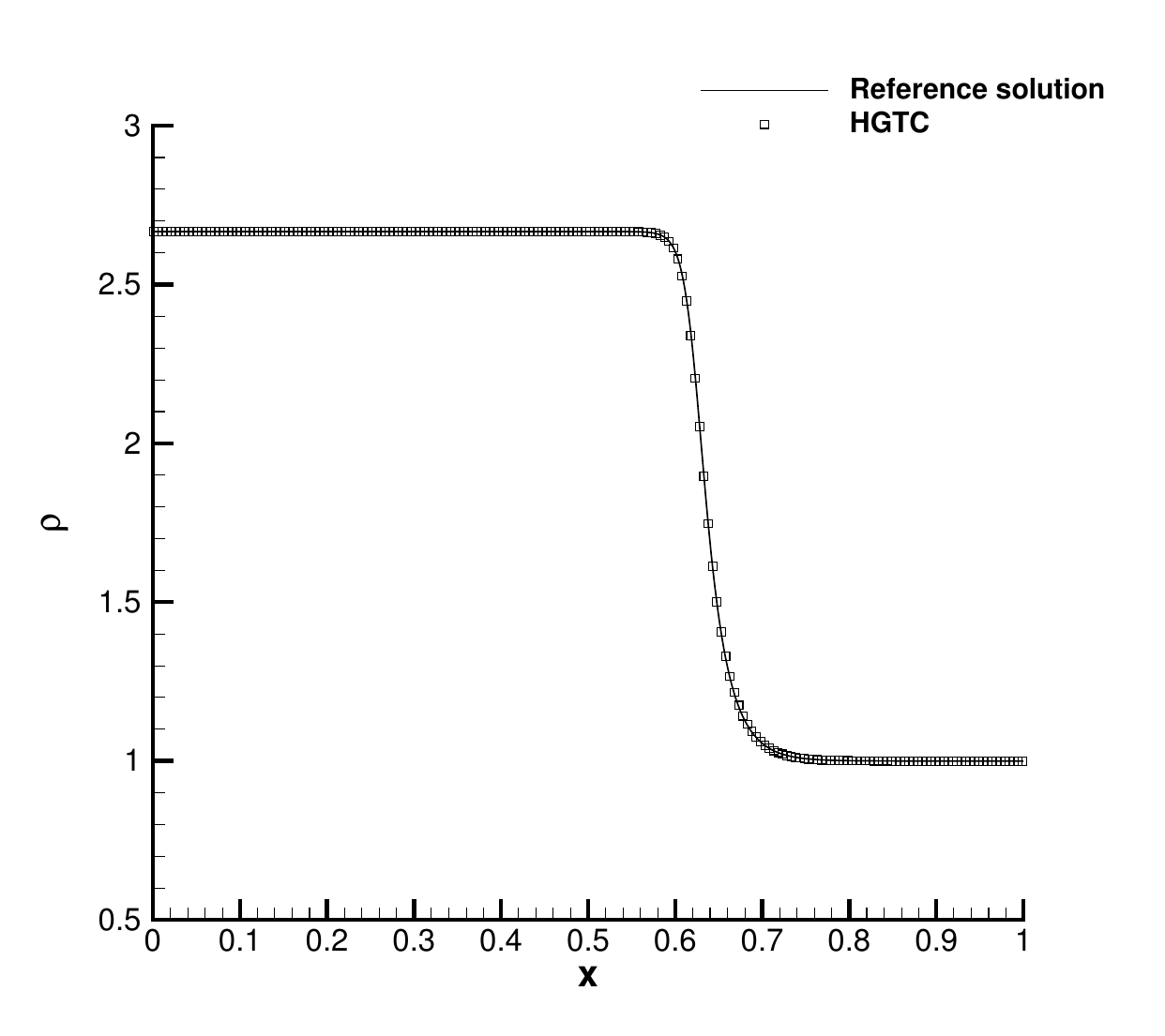} & 
			\includegraphics[width=0.47\textwidth]{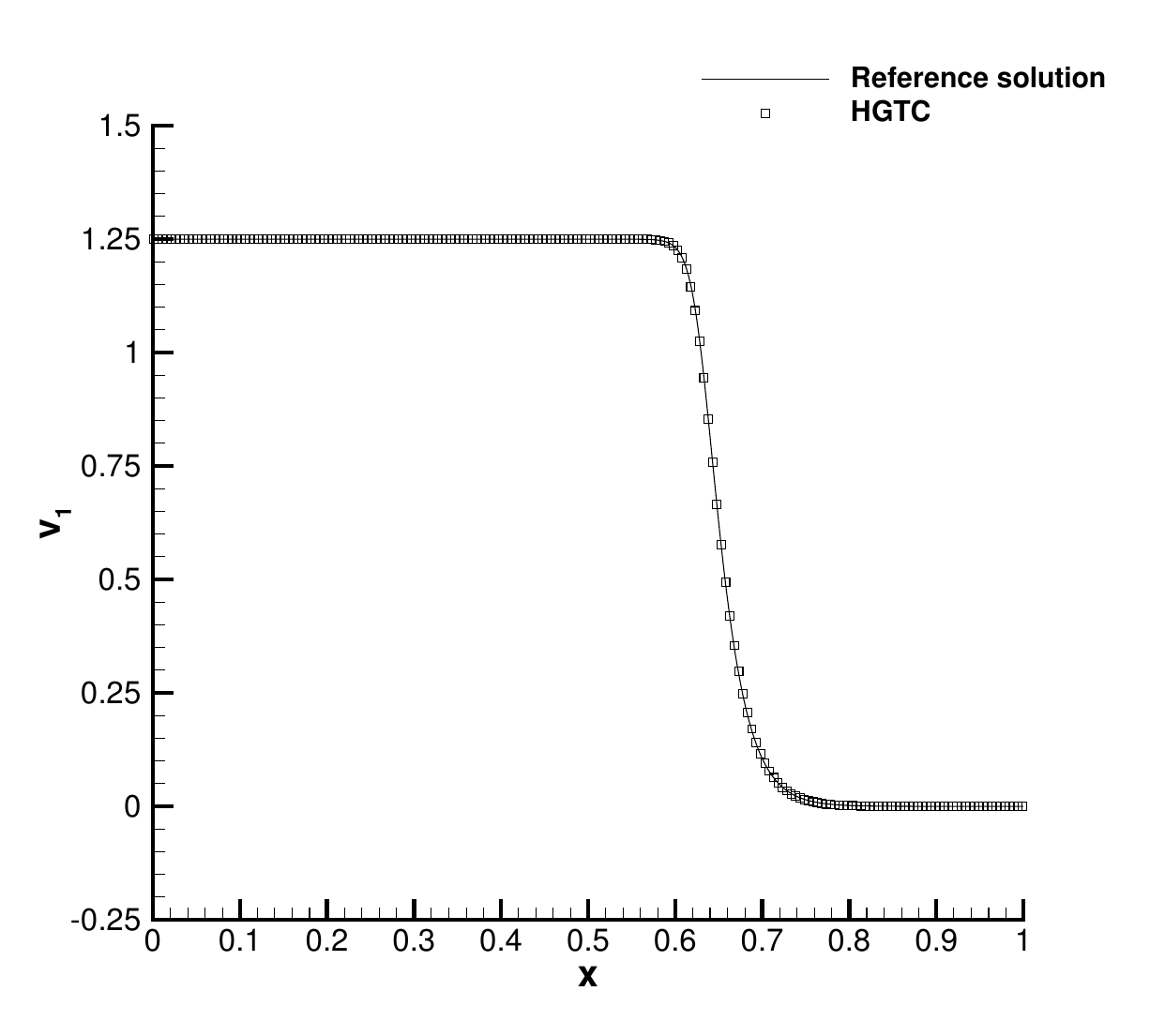} \\   
			\includegraphics[width=0.47\textwidth]{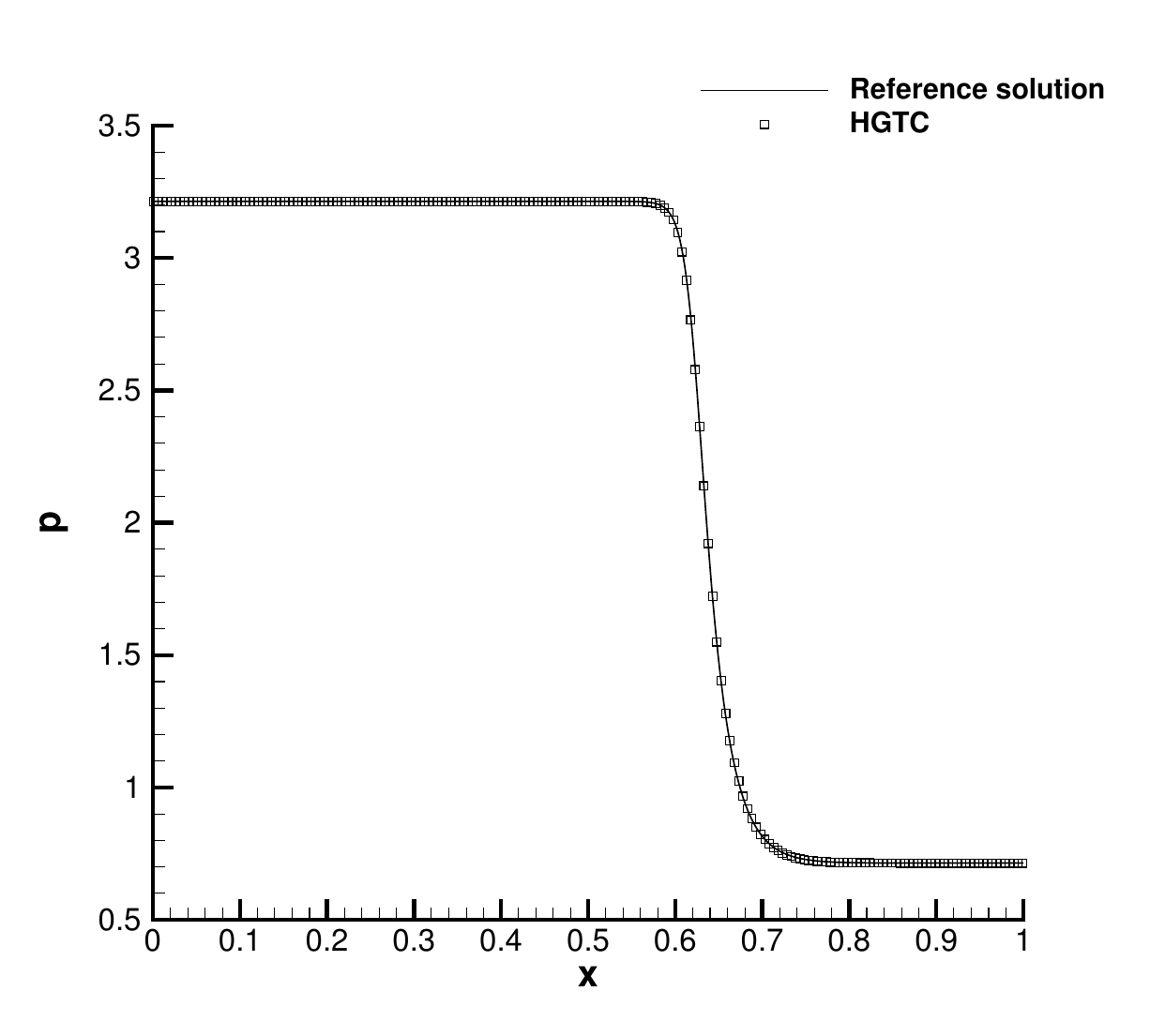} &
			\includegraphics[width=0.47\textwidth]{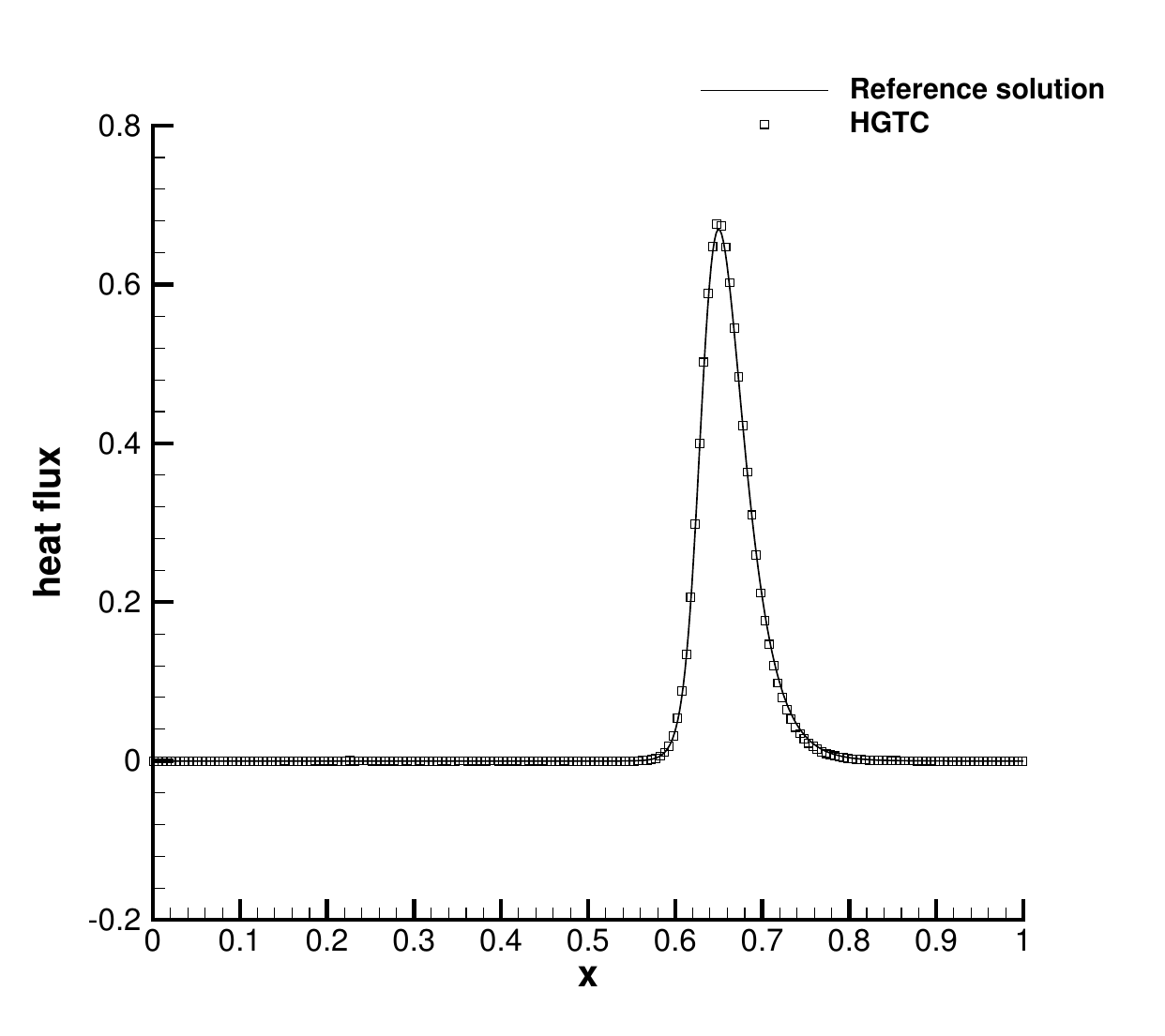} \\ 
		\end{tabular} 
		\caption{Viscous shock problem at time $t_f=0.2$. Numerical results for density, horizontal velocity, pressure and heat flux  extracted with a one-dimensional cut of 200 equidistant points along the $x$-direction at $y=0.1$ compared against the reference solution.}
		\label{fig.viscousshock}
	\end{center}
\end{figure}

Next, we model compressible heat-conducting viscous flows by setting $c_s=c_h=10$, $\mu=2 \cdot 10^{-2}$ and $\kappa=9.3333 \cdot 10^{-2}$. The computational domain is the channel $\Omega=[0;1] \times [0;0.2]$ that is paved with Voronoi polygons of characteristic size of $h=1/1024$. Periodic boundaries are imposed in $y-$direction, while a constant inflow velocity is prescribed for $x=0$ and outflow boundary conditions are set at $x=1$. In \cite{Becker1923}, an exact solution of the one-dimensional compressible Navier-Stokes equations is derived for Prandtl number $\text{Pr}=0.75$ and constant viscosity that involves a stationary viscous shock wave at a shock mach number $M_s$. The Reynolds number is $\text{Re}_s=\rho_0 c_0 M_s L \mu^{-1}$, with the reference length that is assumed to be $L=1$. This is an interesting test case since all the terms characteristics of the one-dimensional compressible Navier-Stokes equations can be verified, including viscous stress and heat conduction. According to \cite{Becker1923}, the exact solution is given in terms of dimensionless density, pressure and velocity. The dimensionless velocity $\bar v = \frac{v}{M_s \, c_0}$ is related to the stationary shock wave, which can be determined as the root of the following equation:
\begin{equation} 
	\label{eqn.alg.u} 
	\frac{|\bar v - 1|}{|\bar v - \lambda^2|^{\lambda^2}} = \left| \frac{1-\lambda^2}{2} \right|^{(1-\lambda^2)} 
	\exp{\left( \frac{3}{4} \textnormal{Re}_s \frac{M_s^2 - 1}{\gamma M_s^2} x \right)},
\end{equation}
with
\begin{equation}
	\lambda^2 = \frac{1+ \frac{\gamma-1}{2}M_s^2}{\frac{\gamma+1}{2}M_s^2}.
\end{equation}
Once the solution of equation \eqref{eqn.alg.u} is computed, the dimensionless velocity $\bar v$ is expressed as a function of $x$. The form of the viscous profile of the dimensionless pressure $\bar p = \frac{p-p_0}{\rho_0 c_0^2 M_s^2}$ is given by the relation 
\begin{equation}
	\label{eqn.alg.p} 
	\bar p = 1 - \bar v +  \frac{1}{2 \gamma}
	\frac{\gamma+1}{\gamma-1} \frac{(\bar v - 1 )}{\bar v} (\bar v - \lambda^2).  
\end{equation}
Finally, the profile of the dimensionless density $\bar \rho = \frac{\rho}{\rho_0}$ is derived from the integrated continuity equation: $\bar \rho \bar v = 1$. Here, we make the simulation unsteady by adding a constant velocity background field $v = M_s c_0$. The initial condition is given by a shock wave centered at $x=0.25$ which is propagating at Mach $M_s=2$ with $\text{Re}_s=100$. The upstream shock state is defined by 
\begin{equation}
	\rho(t=0,\xx) = \rho_0, \qquad \vv(t=0,\xx)=\mathbf{0}, \qquad p_0(t=0,\xx)=1/\gamma,
\end{equation}
with $c_0=1$. The numerical solution obtained with the HGTC schemes without numerical dissipation at the final time $t_f=0.2$ is compared against the reference solution of the one-dimensional compressible Navier-Stokes equations. The numerical solutions for the main primitive variables are plot in Figure \ref{fig.viscousshock}. An excellent agreement is obtained, demonstrating the capability of the HGTC schemes of retrieving the correct physical solution for heat-conducting viscous fluids.

\begin{figure}[!htbp]
	\begin{center}
		\begin{tabular}{cc}   
			\includegraphics[trim=5 0 0 5,clip,width=0.47\textwidth]{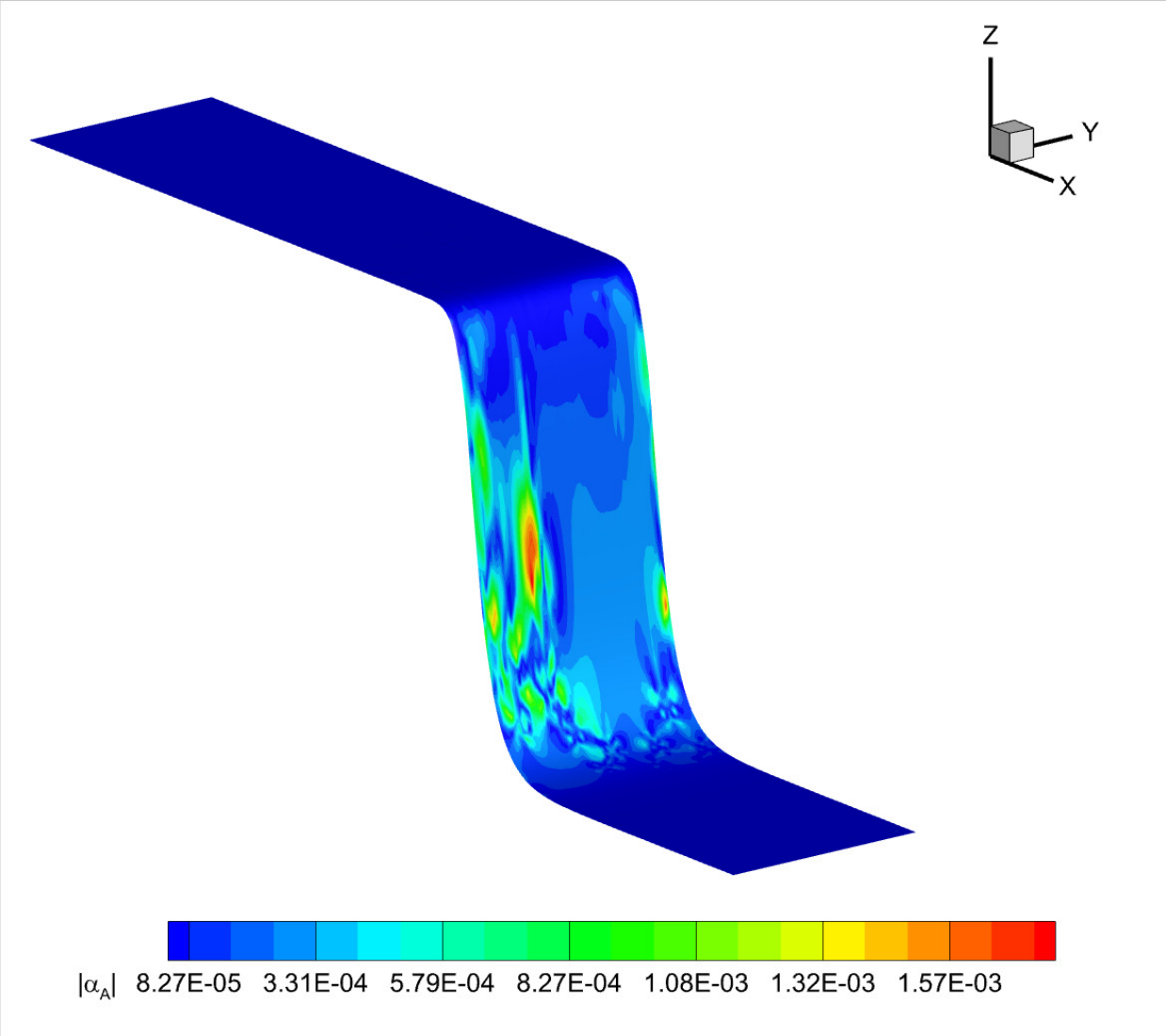} &
			\includegraphics[width=0.47\textwidth]{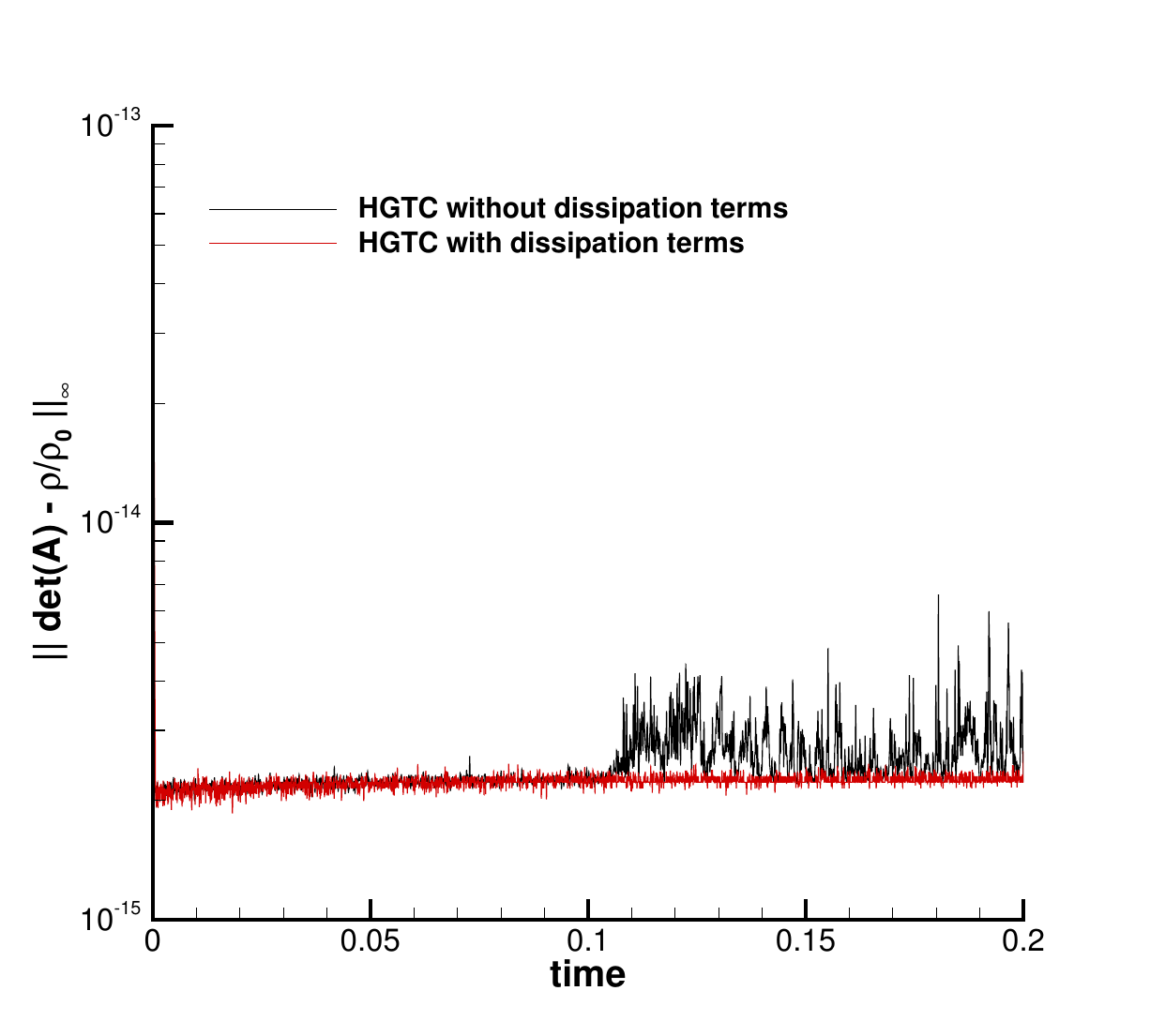} \\
		\end{tabular} 
		\caption{Viscous shock problem at time $t_f=0.2$. Map of the thermodynamic correction factor $|\alpha_\A|$ with a three-dimensional view of the density distribution (left) and time evolution of the mass conservation errors for the HGTC with and without viscous terms (right).}
		\label{fig.viscousshock2}
	\end{center}
\end{figure}

We also plot the correction factor $\alpha_\A$ in Figure \ref{fig.viscousshock2} as well as the time evolution of the total mass conservation errors $\delta_\A$ for the simulations run with and without numerical dissipation. In both cases, the $L_{\infty}$ determinant error remains at machine accuracy.

\subsection{2D Taylor-Green vortex}
\begin{figure}[!htbp]
	\begin{center}
		\begin{tabular}{cc} 
			\includegraphics[trim=2 2 2 2,clip,width=0.47\textwidth]{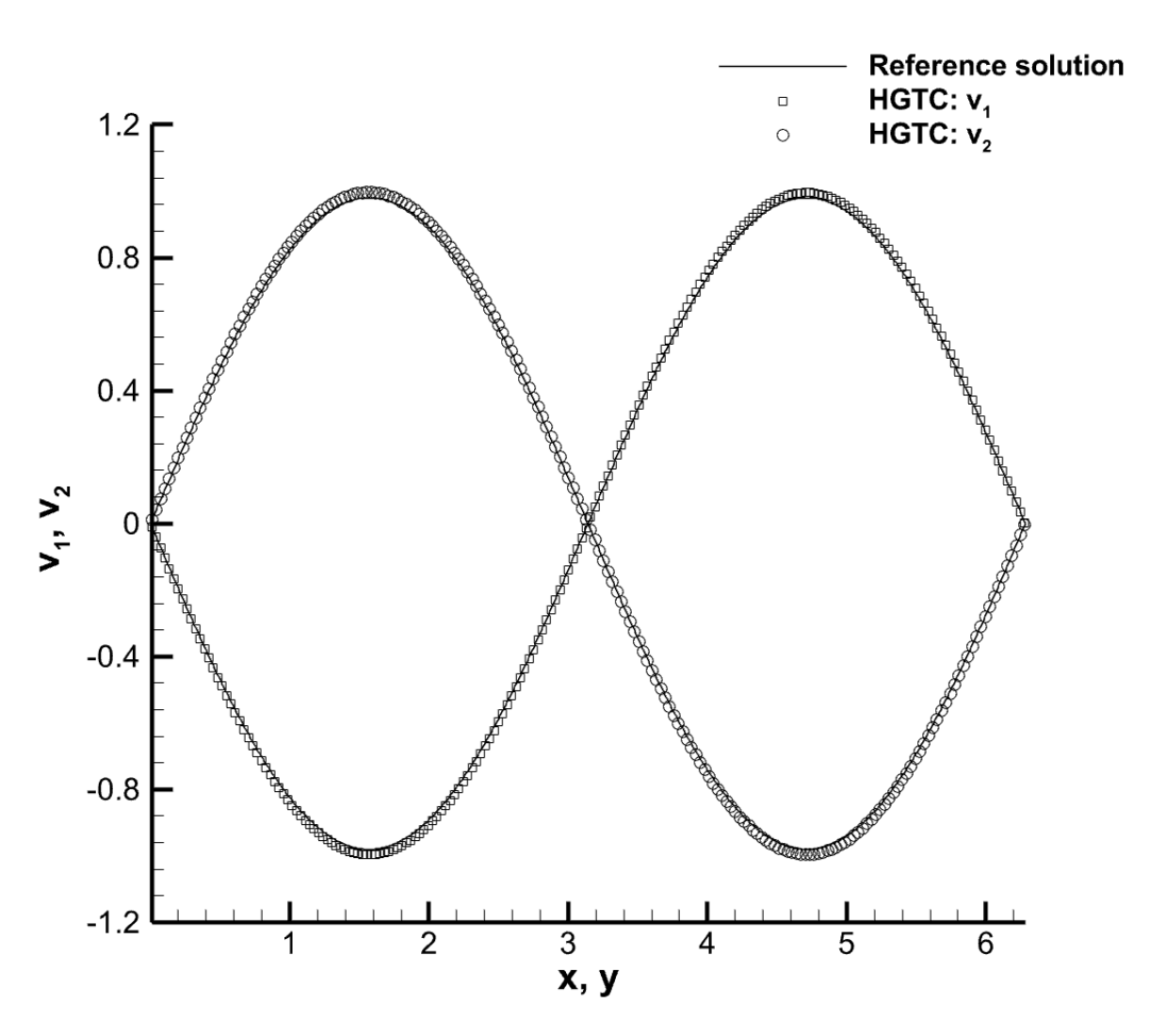} & 
			\includegraphics[trim=2 2 2 2,clip,width=0.47\textwidth]{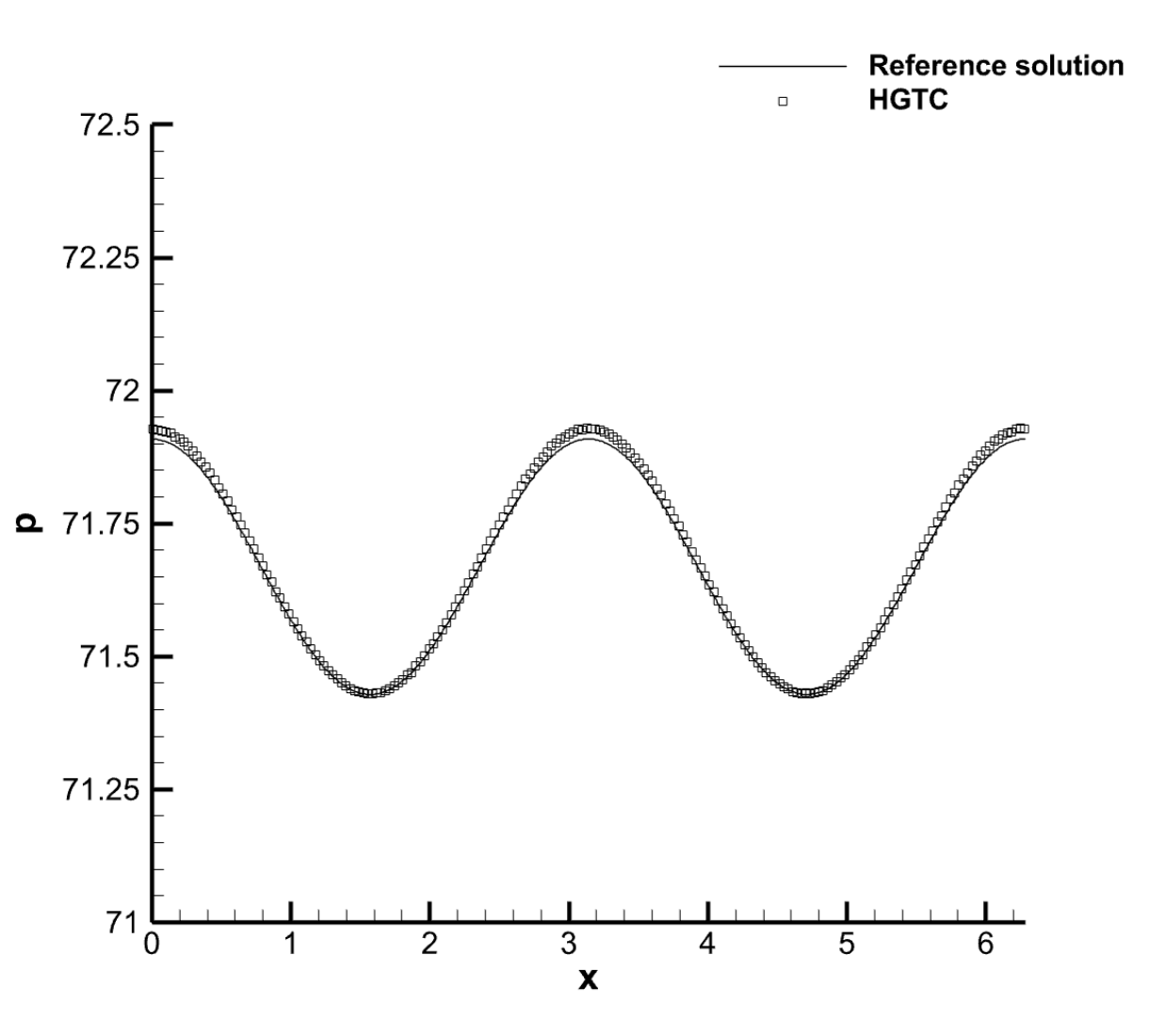} \\
			\includegraphics[trim=2 2 2 2,clip,width=0.47\textwidth]{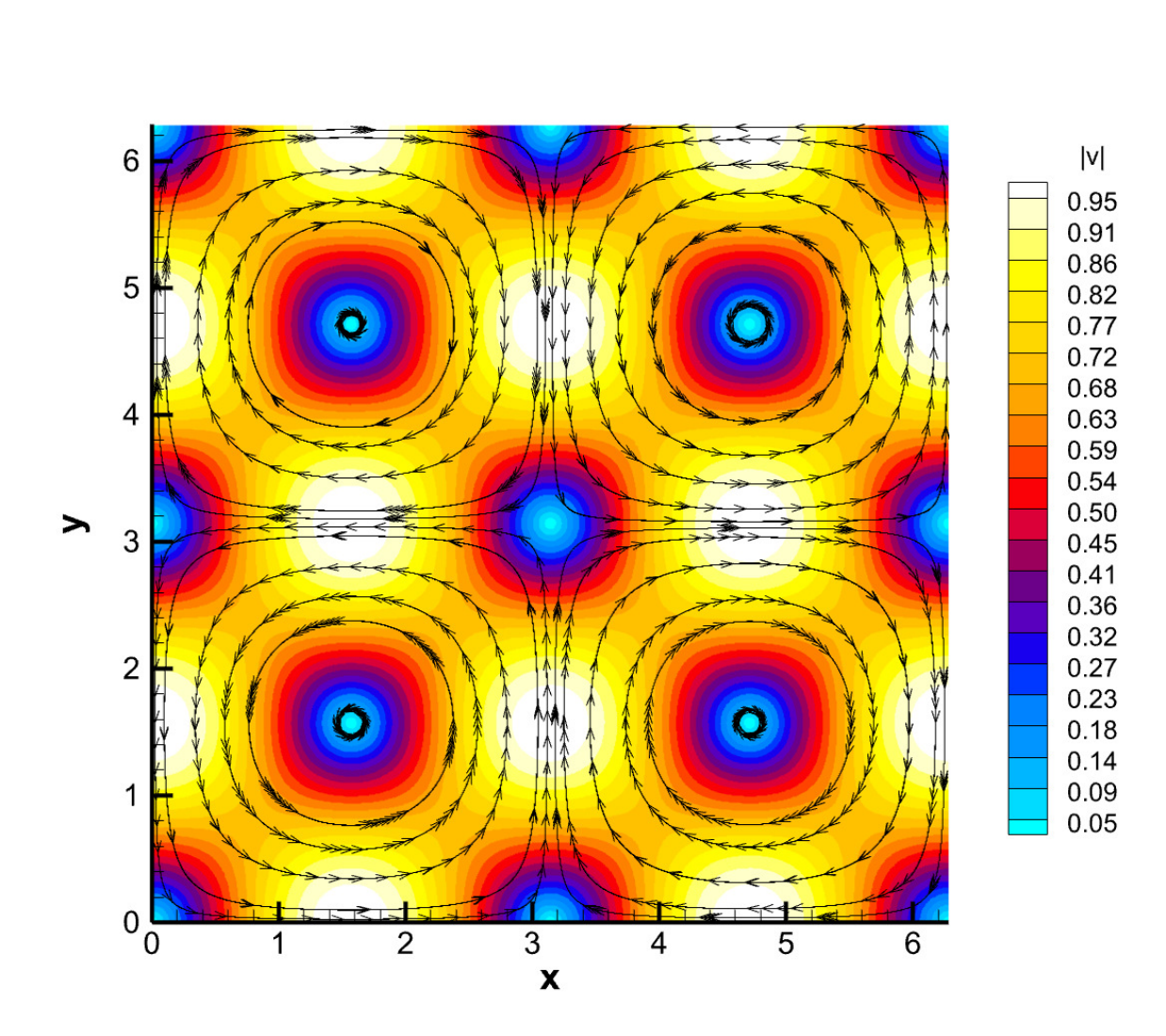} & 
			\includegraphics[trim=2 2 2 2,clip,width=0.47\textwidth]{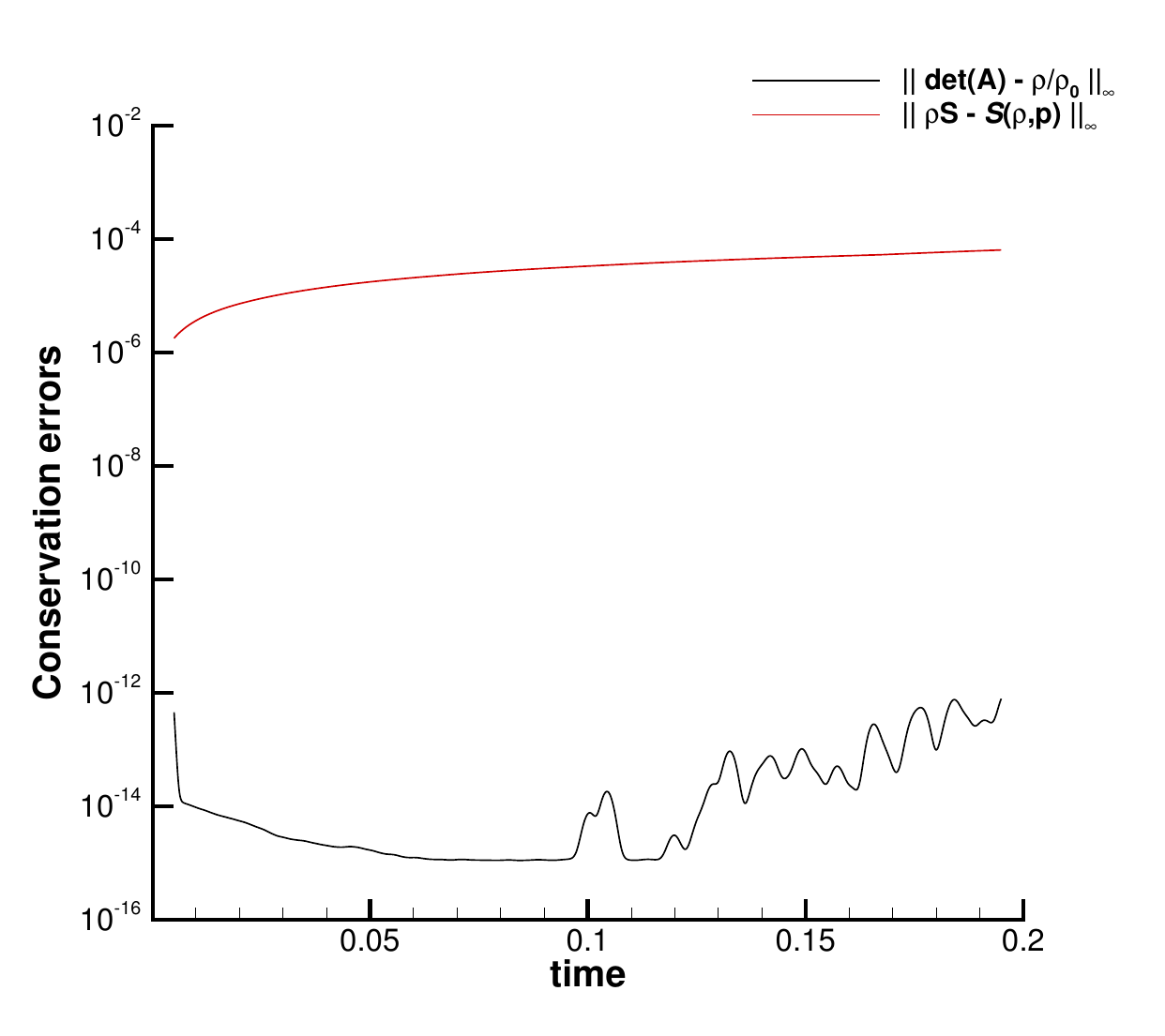} \\
		\end{tabular} 
		\caption{2D Taylor-Green vortex at time $t_f=0.2$ with viscosity $\mu=10^{-2}$. Top: one-dimensional cut of 200 equidistant points along the $x$-axis and the $y-$axis for the velocity components $v_1$ and $v_2$ (left) and for the pressure $p$ (right). Bottom: distribution of the velocity magnitude with stream-traces (left) and time evolution of the mass and entropy conservation errors for the HGTC scheme.}
		\label{fig.TGV2D}
	\end{center}
\end{figure}

The two-dimensional Taylor-Green vortex problem is a well-known test case for the incompressible Navier-Stokes equations. The exact solution writes
\begin{eqnarray}
	u(t,\xx)&=&\phantom{-}\sin(x_1)\cos(x_2) \, e^{-2\nu t},  \nonumber \\
	v(t,\xx)&=&-\cos(x_1)\sin(x_2) \, e^{-2\nu t}, \nonumber \\
	p(t,\xx)&=& C + \frac{1}{4}(\cos(2x_1)+\cos(2x_2)) \, e^{-4\nu t},
	\label{eq:TG_ini}
\end{eqnarray}
where $\nu=\mu/\rho$ denotes the kinematic viscosity of the fluid and the density is $\rho(t,\xx)=1$. To model an incompressible viscous fluid, we set $c_s=10$ and $\mu=10^{-2}$, and the additive constant to the pressure is chosen to be $C=100/\gamma$ so that a maximum Mach number of $0.1$ is retrieved. Heat conduction is neglected, thus we set $c_h=0$. The initial condition is provided by the exact solution \eqref{eq:TG_ini} at time $t=0$. The computational domain is given by $\Omega=[0;2\pi]^2$ with periodic boundary conditions everywhere, and it is paved with a Voronoi grid of characteristic mesh size $h=2\pi/200$. Figure \ref{fig.TGV2D} depicts the numerical results at the final time $t_f=0.2$ that are compared against the reference solution, obtaining an excellent matching. Furthermore, we also show the time evolution of the mass and entropy conservation errors, that remain bounded and preserved thanks to the compatibility corrections of our novel HGTC schemes.

\subsection{Solid rotor problem}
\begin{figure}[!htbp]
	\begin{center}
		\begin{tabular}{cc} 
			\includegraphics[width=0.47\textwidth]{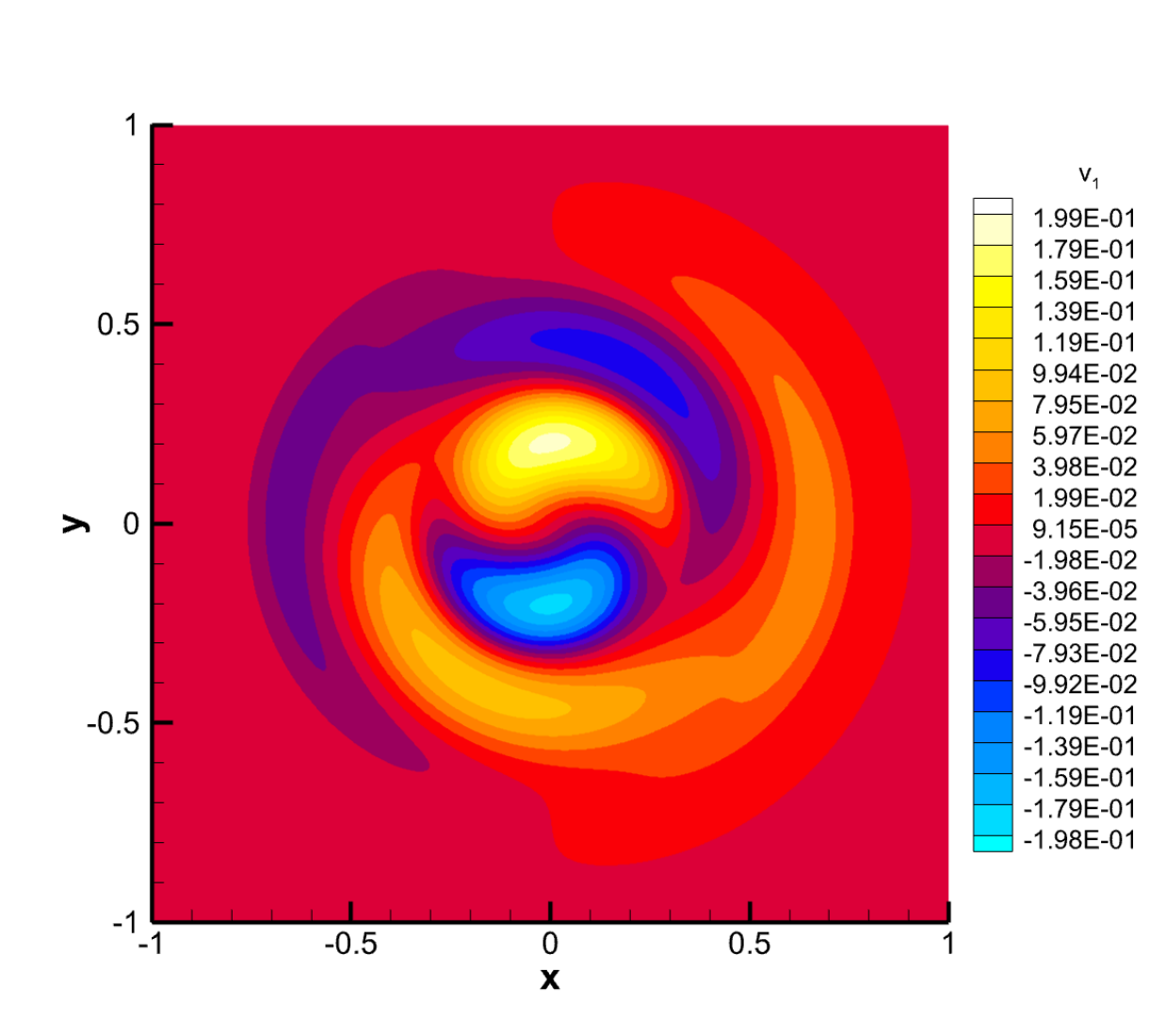} & 
			\includegraphics[width=0.47\textwidth]{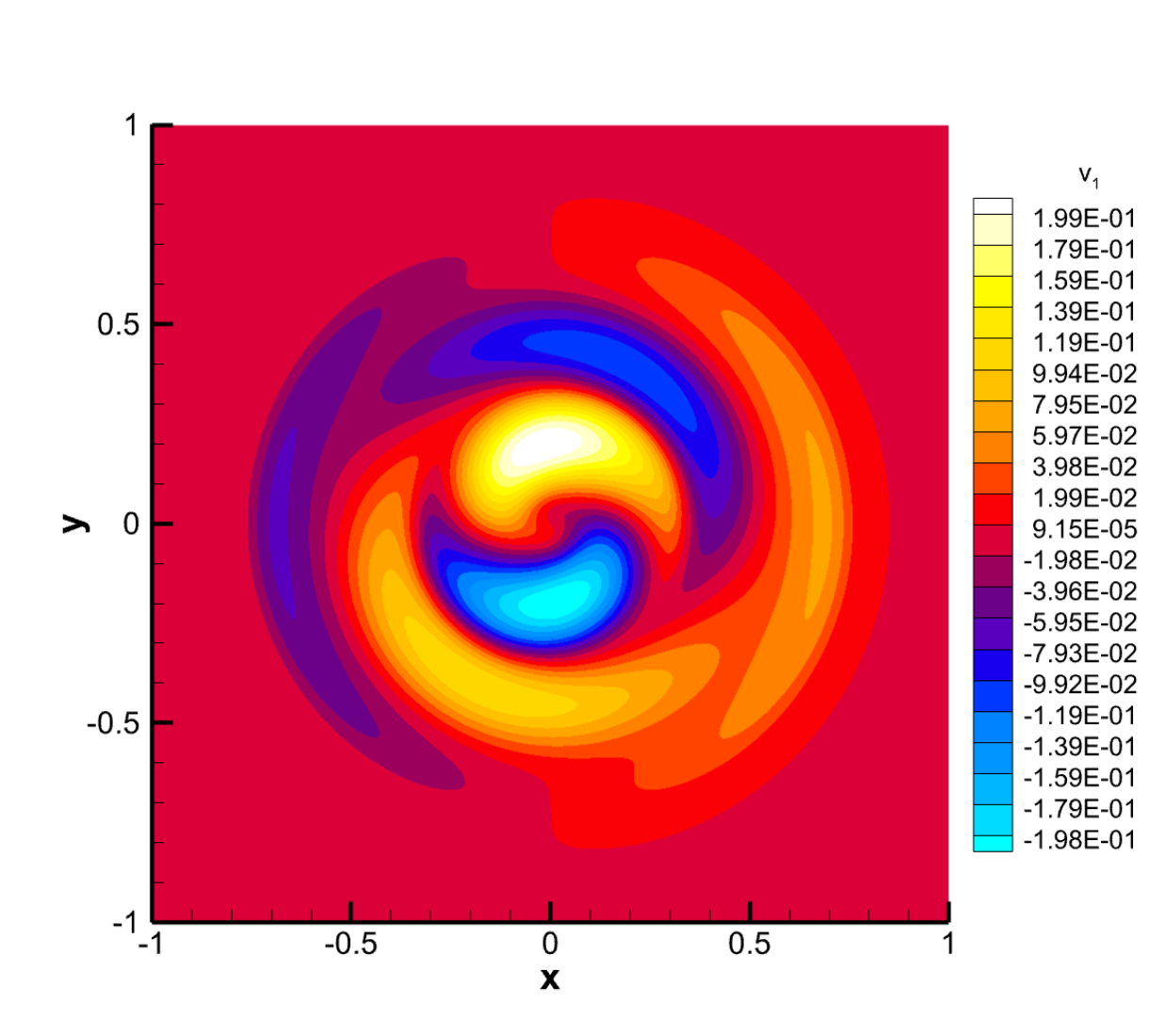} \\ 
			\includegraphics[width=0.47\textwidth]{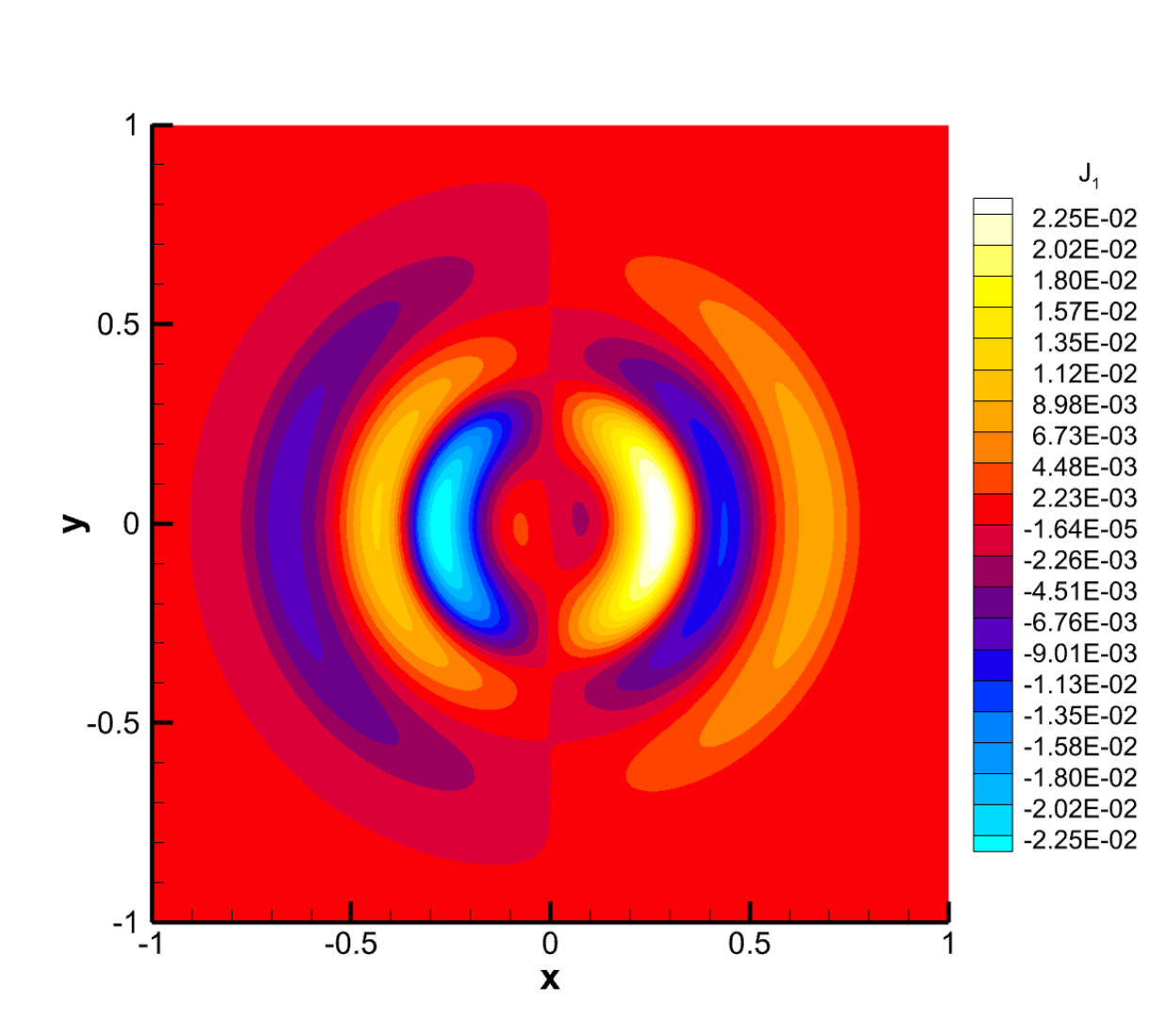} & 
			\includegraphics[width=0.47\textwidth]{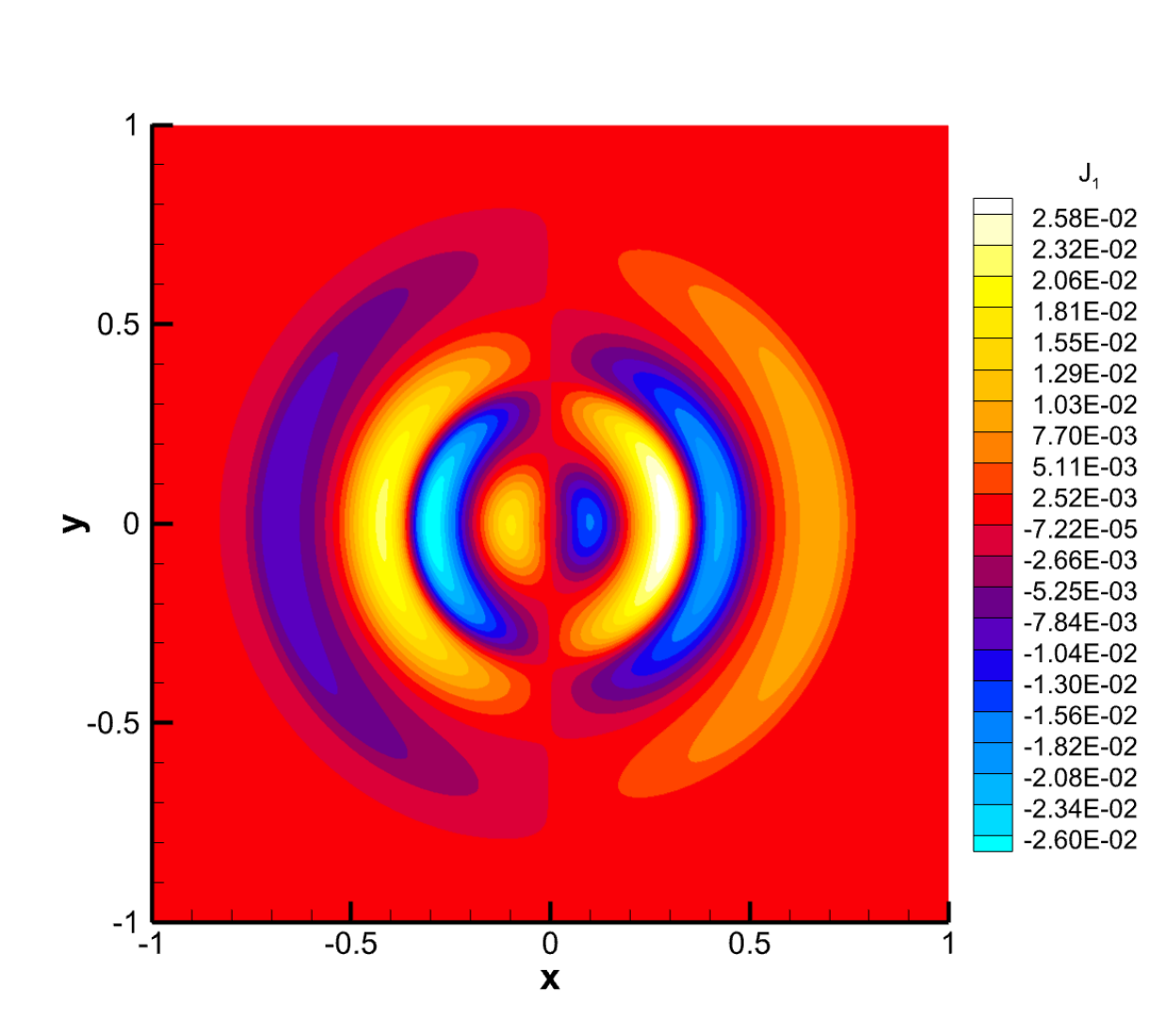} \\  
		\end{tabular} 
		\caption{Solid rotor problem at time $t_f=0.3$. Top: numerical results for the horizontal velocity $v_1$ with $h=1/128$ (left) and $h=1/256$ (right). Bottom: numerical results for the thermal impulse component $J_1$ with $h=1/128$ (left) and $h=1/256$ (right).}
		\label{fig.solidRotor}
	\end{center}
\end{figure}

Finally, a test case for solid mechanics is solved, namely the solid rotor problem introduced in \cite{HyperHypo,SIGPR}. The relaxation times of the mathematical model are set to $\tau_1=\tau_2=10^{20}$, hence nonlinear hyperelastic solids are genuinely modeled by the governing PDE presented in \cite{PeshRom2014}. We fix $c_s=c_h=1$ and the final time of the simulation is $t_f=0.3$. The computational domain is $\Omega=[-1;1]^2$ with periodic boundaries, and the initial condition of the material writes
\begin{equation}
	(\rho,v_1,v_2,v_3,p)=\left\{
	\begin{array}{lll}
		(1, \, -x_2/R, \, x_1/R, \, 0,1) & & r<R \\
		(1,\, 0,\, 0,\, 0, \, 1) & &  r\geq R
	\end{array} \right., \qquad t=0, \quad \xx \in \Omega,
\end{equation}
with the initial discontinuity located at $R=0.2$ and $r=\sqrt{x_1^2+x_2^2}$. To show mesh convergence, we run the solid rotor problem on two different meshes with characteristic size of $h=1/256$ and $h=1/128$. The results are compared with each other in Figure \ref{fig.solidRotor}, where the horizontal velocity distribution is plot. The maps of the scalar correction factors $\alpha_\A$ and $\alpha_S$ at the final time level are also depicted. We observe that both corrections may act at the same spatial locations without negatively interfering between each other.

\begin{figure}[!htbp]
	\begin{center}
		\begin{tabular}{cc} 
			\includegraphics[width=0.47\textwidth]{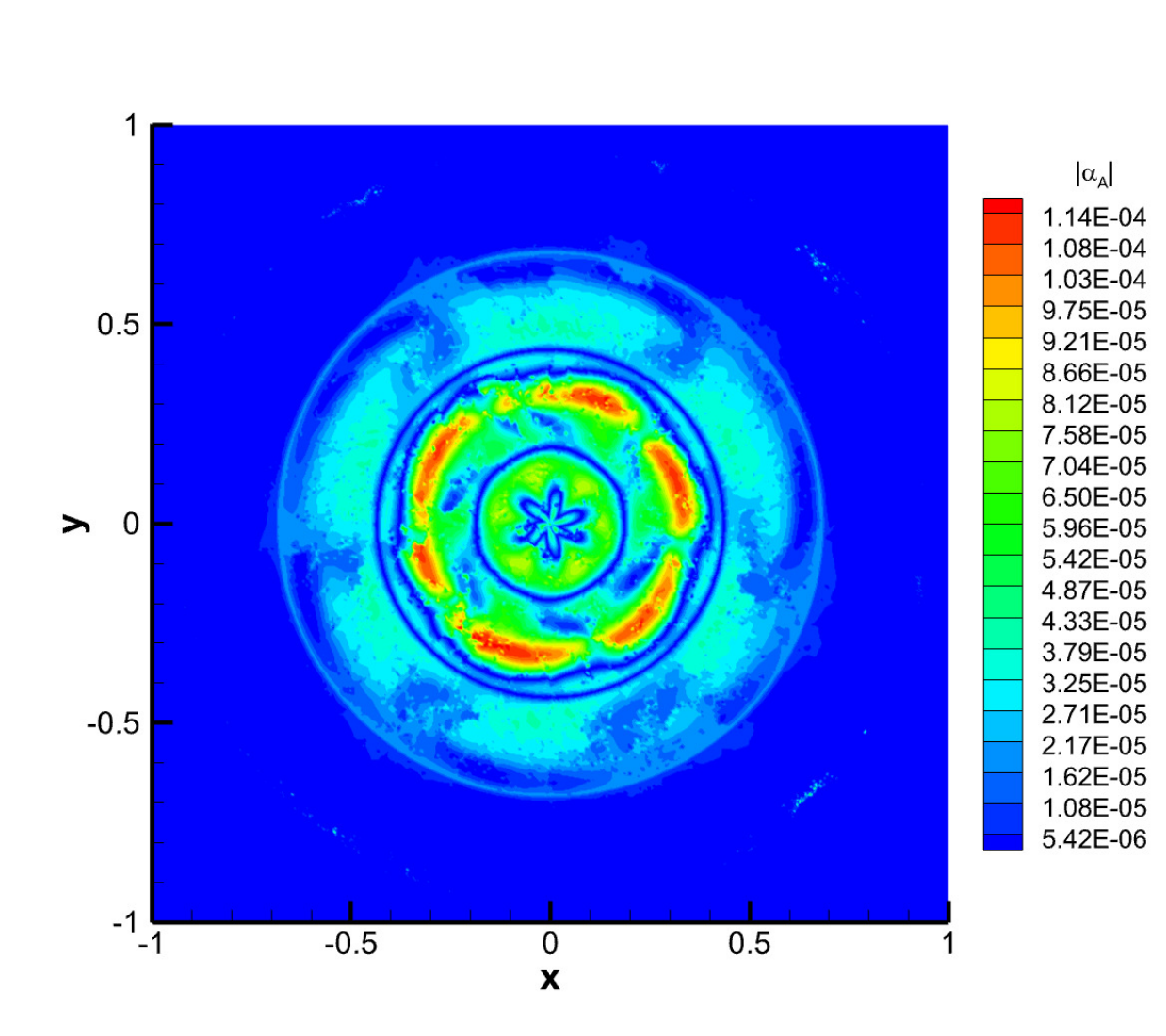} &
			\includegraphics[width=0.47\textwidth]{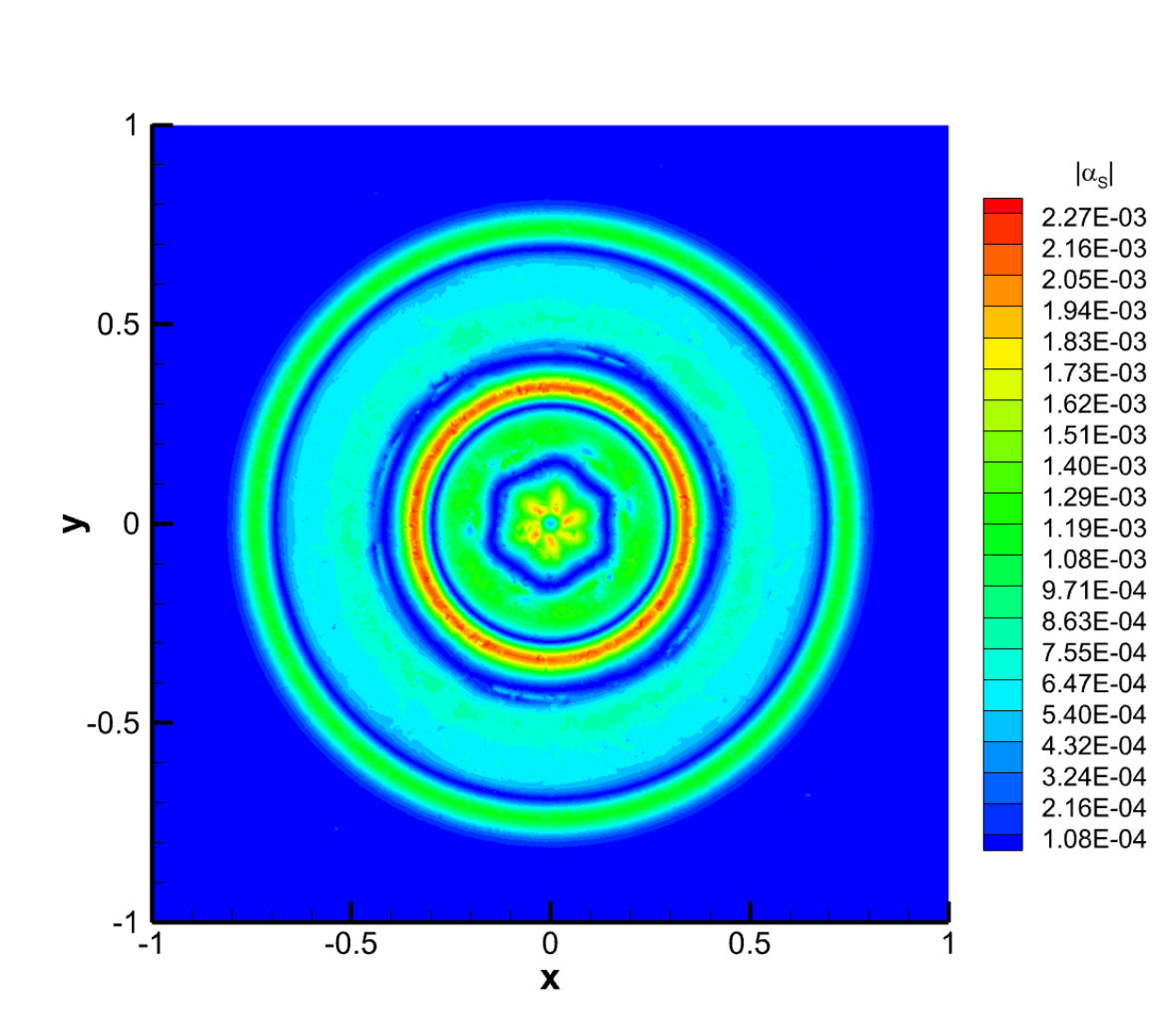} \\ 
		\end{tabular} 
		\caption{Solid rotor problem at time $t_f=0.3$. Map of the geometric correction factor $|\alpha_\A|$ (left) and of the thermodynamic correction factor $|\alpha_S|$ (right) obtained with a characteristic mesh size of $h=1/512$.}
		\label{fig.solidRotor2}
	\end{center}
\end{figure}

\section{Conclusions} \label{sec.concl}
In this paper we have presented a novel finite volume scheme on unstructured Voronoi meshes for the solution of a reduced unified model for continuum mechanics, where the mass conservation equation is discarded. This has been achieved enforcing the compatibility of the new schemes with the Geometric Conservation Law that links the distortion tensor to the density within each control volume. The geometric compatibility is obtained by introducing a new generalized concept of potential, that is assumed to be the determinant of the distortion tensor. Consequently, a set of associated pseudo-dual variables is retrieved, which play the role of the thermodynamic variables for the total energy potential. By means of a conservative correction directly embedded in the numerical fluxes, the novel schemes are proven to be compliant with the GCL at the semi-discrete level. Once the geometric compatibility is achieved, thermodynamic compatibility is also guaranteed using the same strategy that derives from the formalism of symmetric and hyperbolic thermodynamically compatible (SHTC) systems introduced by Godunov in 1961. These two corrections can coexist at the discrete level and they do not interfere with each other, hence making it possible for the first time on unstructured fixed grids to ensure geometric and thermodynamic compatibility at the same time. Two theorems demonstrate that these properties are respected at the semi-discrete level. A two-dimensional first order finite volume scheme with up to fourth order Runge-Kutta time integrators has been implemented and tested. A large suite of test cases is shown to numerically prove the structure preserving properties of the new schemes.

In the future we plan to exploit this strategy to tackle other types of constraints, namely involution-constraints like the solenoidal property of the magnetic field in magnetohydrodynamics or the irrotational behavior of the deformation gradient and the thermal impulse vector in reversible processes in solid mechanics. Finally, the extension of the proposed approach to high order discontinuous Galerkin schemes is also foreseen as well as the development of implicit-explicit \cite{BosPar2021,SIGPR} asymptotic preserving discretizations to make the numerical schemes consistent with the Navier-Stokes-Fourier limit exhibited by the governing equations.

\section*{Acknowledgments}

WB received financial support by Fondazione Cariplo and Fondazione CDP (Italy) under the grant No. 2022-1895 and by the Italian Ministry of University and Research (MUR) with the PRIN Project 2022 No. 2022N9BM3N. WB and RL acknowledge the support of LRC Anabase and the CEA-Cesta.

\bibliographystyle{plain}
\bibliography{biblio}

\appendix

\section{Geometric Conservation Law derived from the determinant potential} \label{app.gcl}
In this appendix we show that the pseudo-\textit{Gibbs relation}$  $ \eqref{eqn.GibbsA} holds true. The starting point is the evolution equation of the distortion tensor $\A=\{A_{ik}\}$, in which the viscous terms are neglected. Let us explicitly compute the dual variables $\w$ of the potential $|\A|$:
\begin{equation}
	\w :=\{w_{ik}\}= \partial_{A_{ik}}|\A| = \left(\begin{array}{ccc}
		\phantom{-}A_{22} A_{33}-A_{23} A_{32} & 
		-A_{21} A_{33}+A_{23} A_{31} &
		\phantom{-}A_{21} A_{32}-A_{22} A_{31} 
		\\ [2mm] 
		-A_{12} A_{33}+A_{13} A_{32} &
		\phantom{-}A_{11} A_{33}-A_{13} A_{31} &
		-A_{11} A_{32}+A_{12} A_{31} 
		\\ [2mm] 
		\phantom{-}A_{12} A_{23}-A_{13} A_{22} &
		-A_{11} A_{23}+A_{13} A_{21} &
		\phantom{-}A_{11} A_{22}-A_{12} A_{21} 
	\end{array}\right).
\end{equation}
The determinant of the distortion tensor $\A$ is explicitly given by
\begin{equation}
	|\A| = A_{11} A_{22} A_{33}-A_{11} A_{23} A_{32}-A_{12} A_{21} A_{33}+A_{12} A_{23} A_{31}+A_{13} A_{21} A_{32}-A_{13} A_{22} A_{31}.
\end{equation}

Firstly, we investigate the compatibility with the source terms which are present in the equations \eqref{eqn.deformation_red}. The source terms $\mathbf{S}_{\A}=\{S_{\A,ik}\}$ are given by
\begin{equation}
	S_{\A,ik} = - \frac{3}{\tau_1} \, |\A|^{\frac{5}{3}} \, A_{im} \, \mathring{G}_{mk}, \qquad \mathring{G}_{mk} = {G}_{mk} - \frac{1}{3} \, G_{jj} \delta_{mk}, \quad G_{mk} = A_{jm} A_{jk}.
\end{equation}
Multiplication of the source terms with the dual variables yields the contributions $\{\tilde{S}_{\A,ik}\} =  \{ w_{ik} \, S_{\A,ik} \}$, which write
\begin{equation}
	\begin{array}{l}
		\tilde{S}_{\A,11} = -\frac{(-2 A_{11}^{3}+(-2 A_{12}^{2}-2 A_{13}^{2}-2 A_{21}^{2}+A_{22}^{2}+A_{23}^{2}-2 A_{31}^{2}+A_{32}^{2}+A_{33}^{2}) A_{11}+(-3 A_{21} A_{22}-3 A_{31} A_{32}) A_{12}-3 A_{13} (A_{21} A_{23}+A_{31} A_{33})) (-A_{22} A_{33}+A_{23} A_{32}) |\A|^{\frac{5}{3}}}{\tau_1}, \\ [2mm]
		\tilde{S}_{\A,12} = \frac{2 (A_{12}^{3}+(A_{32}^{2}-\frac{A_{33}^{2}}{2}+A_{11}^{2}+A_{13}^{2}-\frac{A_{21}^{2}}{2}+A_{22}^{2}-\frac{A_{23}^{2}}{2}-\frac{A_{31}^{2}}{2}) A_{12}+(\frac{3 A_{31} A_{32}}{2}+\frac{3 A_{21} A_{22}}{2}) A_{11}+\frac{3 A_{13} (A_{22} A_{23}+A_{32} A_{33})}{2}) (A_{21} A_{33}-A_{23} A_{31}) |\A|^{\frac{5}{3}}}{\tau_1}, \\ [2mm]
		\tilde{S}_{\A,13} = \frac{(-2 A_{13}^{3}+(-2 A_{11}^{2}-2 A_{12}^{2}+A_{21}^{2}+A_{22}^{2}-2 A_{23}^{2}+A_{31}^{2}+A_{32}^{2}-2 A_{33}^{2}) A_{13}+(-3 A_{21} A_{23}-3 A_{31} A_{33}) A_{11}-3 A_{12} (A_{22} A_{23}+A_{32} A_{33})) (A_{21} A_{32}-A_{22} A_{31}) |\A|^{\frac{5}{3}}}{\tau_1}, \\ [2mm]
		\tilde{S}_{\A,21} = \frac{(-2 A_{21}^{3}+(-2 A_{11}^{2}+A_{12}^{2}+A_{13}^{2}-2 A_{22}^{2}-2 A_{23}^{2}-2 A_{31}^{2}+A_{32}^{2}+A_{33}^{2}) A_{21}+(-3 A_{12} A_{22}-3 A_{13} A_{23}) A_{11}-3 A_{31} (A_{22} A_{32}+A_{23} A_{33})) (-A_{12} A_{33}+A_{13} A_{32}) |\A|^{\frac{5}{3}}}{\tau_1}, \\ [2mm]
		\tilde{S}_{\A,22} = -\frac{2 (A_{11} A_{33}-A_{13} A_{31}) (A_{22}^{3}+(A_{32}^{2}-\frac{A_{33}^{2}}{2}-\frac{A_{11}^{2}}{2}+A_{12}^{2}-\frac{A_{13}^{2}}{2}+A_{21}^{2}+A_{23}^{2}-\frac{A_{31}^{2}}{2}) A_{22}+\frac{3 A_{11} A_{12} A_{21}}{2}+\frac{3 A_{12} A_{13} A_{23}}{2}+\frac{3 A_{32} (A_{21} A_{31}+A_{23} A_{33})}{2}) |\A|^{\frac{5}{3}}}{\tau_1}, \\ [2mm]
		\tilde{S}_{\A,23} = -\frac{(A_{11} A_{32}-A_{12} A_{31}) (-2 A_{23}^{3}+(A_{11}^{2}+A_{12}^{2}-2 A_{13}^{2}-2 A_{21}^{2}-2 A_{22}^{2}+A_{31}^{2}+A_{32}^{2}-2 A_{33}^{2}) A_{23}-3 A_{11} A_{13} A_{21}-3 A_{12} A_{13} A_{22}-3 A_{33} (A_{21} A_{31}+A_{22} A_{32})) |\A|^{\frac{5}{3}}}{\tau_1}, \\ [2mm]
		\tilde{S}_{\A,31} = -\frac{2 (A_{12} A_{23}-A_{13} A_{22}) (A_{31}^{3}+(A_{32}^{2}+A_{33}^{2}+A_{11}^{2}-\frac{A_{12}^{2}}{2}-\frac{A_{13}^{2}}{2}+A_{21}^{2}-\frac{A_{22}^{2}}{2}-\frac{A_{23}^{2}}{2}) A_{31}+(\frac{3 A_{32} A_{12}}{2}+\frac{3 A_{33} A_{13}}{2}) A_{11}+\frac{3 A_{21} (A_{22} A_{32}+A_{23} A_{33})}{2}) |\A|^{\frac{5}{3}}}{\tau_1}, \\ [2mm]
		\tilde{S}_{\A,32} = \frac{2 (A_{32}^{3}+(A_{33}^{2}-\frac{A_{11}^{2}}{2}+A_{12}^{2}-\frac{A_{13}^{2}}{2}-\frac{A_{21}^{2}}{2}+A_{22}^{2}-\frac{A_{23}^{2}}{2}+A_{31}^{2}) A_{32}+\frac{3 A_{11} A_{12} A_{31}}{2}+\frac{3 A_{12} A_{13} A_{33}}{2}+\frac{3 A_{22} (A_{21} A_{31}+A_{23} A_{33})}{2}) (A_{11} A_{23}-A_{13} A_{21}) |\A|^{\frac{5}{3}}}{\tau_1}, \\ [2mm]
		\tilde{S}_{\A,33} = -\frac{2 \left(A_{11} A_{22}-A_{12} A_{21}\right) \left(A_{33}^{3}+\left(A_{32}^{2}-\frac{A_{11}^{2}}{2}-\frac{A_{12}^{2}}{2}+A_{13}^{2}-\frac{A_{21}^{2}}{2}-\frac{A_{22}^{2}}{2}+A_{23}^{2}+A_{31}^{2}\right) A_{33}+\frac{3 A_{11} A_{13} A_{31}}{2}+\frac{3 A_{12} A_{13} A_{32}}{2}+\frac{3 A_{23} \left(A_{21} A_{31}+A_{22} A_{32}\right)}{2}\right) |\A|^{\frac{5}{3}}}{\tau_1}. \\ [2mm]
	\end{array}
\end{equation}
At the aid of a linear algebra software \cite{maple}, by summing up all the above terms, i.e. dot multiplying the source terms with the dual variables, one obtains 
\begin{equation}
	\tilde{S}_{\A,ik} = -w_{ik} \cdot \frac{\alpha_{ik}}{\theta(\tau_1)} = 0.
\end{equation}
Thus, we retrieve no source on the right hand side of the GCL \eqref{eqn.detApde} as expected.

Next, the compatibility with the flux and non-conservative terms on the left hand side of \eqref{eqn.deformation} has to be verified. To that aim, the equation of the distortion tensor as well as the GCL are written in fully non-conservative form as follows:
\begin{align}
	&\frac{\partial A_{i k}}{\partial t}
	+A_{im} \frac{\partial v_m}{\partial x_k} +v_m \frac{\partial A_{ik}}{\partial x_m}	=  -\dfrac{ \alpha_{ik} }{\theta_1(\tau_1)}, \label{eqn.Anc} \\
	&\frac{\partial |\A|}{\partial t} + |\A| \frac{\partial v_k}{\partial x_k} + v_k \cdot \frac{\partial |\A|}{\partial x_k} = 0. \label{eqn.detAnc}
\end{align}
The non-conservative terms in \eqref{eqn.Anc} for each component of the distortion tensor $\A$, i.e. $D_{ik}=A_{im} \frac{\partial v_m}{\partial x_k} +v_m \frac{\partial A_{ik}}{\partial x_m}$, are given by
\begin{equation}
	\mathbf{D} =\{D_{ik}\} = \left(\begin{array}{c}
		A_{11} \left(\frac{\partial v_{1}}{\partial x_{1}}\right)+v_{1} \left(\frac{\partial A_{11}}{\partial x_{1}}\right)+A_{12} \left(\frac{\partial v_{2}}{\partial x_{1}}\right)+v_{2} \left(\frac{\partial A_{11}}{\partial x_{2}}\right)+A_{13} \left(\frac{\partial v_{3}}{\partial x_{1}}\right)+
		v_{3} \left(\frac{\partial A_{11}}{\partial x_{3}}\right) 
		\\ [3mm] 
		A_{11} \left(\frac{\partial v_{1}}{\partial x_{2}}\right)+v_{1} \left(\frac{\partial A_{12}}{\partial x_{1}}\right)+A_{12} \left(\frac{\partial v_{2}}{\partial x_{2}}\right)+v_{2} \left(\frac{\partial A_{12}}{\partial x_{2}}\right)+A_{13} \left(\frac{\partial v_{3}}{\partial x_{2}}\right)+
		v_{3} \left(\frac{\partial A_{12}}{\partial x_{3}}\right) 
		\\ [3mm] 
		A_{11} \left(\frac{\partial v_{1}}{\partial x_{3}}\right)+v_{1} \left(\frac{\partial A_{13}}{\partial x_{1}}\right)+A_{12} \left(\frac{\partial v_{2}}{\partial x_{3}}\right)+v_{2} \left(\frac{\partial A_{13}}{\partial x_{2}}\right)+A_{13} \left(\frac{\partial v_{3}}{\partial x_{3}}\right)+
		v_{3} \left(\frac{\partial A_{13}}{\partial x_{3}}\right) 
		\\ [3mm] 
		A_{21} \left(\frac{\partial v_{1}}{\partial x_{1}}\right)+v_{1} \left(\frac{\partial A_{21}}{\partial x_{1}}\right)+A_{22} \left(\frac{\partial v_{2}}{\partial x_{1}}\right)+v_{2} \left(\frac{\partial A_{21}}{\partial x_{2}}\right)+A_{23} \left(\frac{\partial v_{3}}{\partial x_{1}}\right)+
		v_{3} \left(\frac{\partial A_{21}}{\partial x_{3}}\right) 
		\\ [3mm] 
		A_{21} \left(\frac{\partial v_{1}}{\partial x_{2}}\right)+v_{1} \left(\frac{\partial A_{22}}{\partial x_{1}}\right)+A_{22} \left(\frac{\partial v_{2}}{\partial x_{2}}\right)+v_{2} \left(\frac{\partial A_{22}}{\partial x_{2}}\right)+A_{23} \left(\frac{\partial v_{3}}{\partial x_{2}}\right)+
		v_{3} \left(\frac{\partial A_{22}}{\partial x_{3}}\right) 
		\\ [3mm] 
		A_{21} \left(\frac{\partial v_{1}}{\partial x_{3}}\right)+v_{1} \left(\frac{\partial A_{23}}{\partial x_{1}}\right)+A_{22} \left(\frac{\partial v_{2}}{\partial x_{3}}\right)+v_{2} \left(\frac{\partial A_{23}}{\partial x_{2}}\right)+A_{23} \left(\frac{\partial v_{3}}{\partial x_{3}}\right)+
		v_{3} \left(\frac{\partial A_{23}}{\partial x_{3}}\right) 
		\\ [3mm] 
		A_{31} \left(\frac{\partial v_{1}}{\partial x_{1}}\right)+v_{1} \left(\frac{\partial A_{31}}{\partial x_{1}}\right)+A_{32} \left(\frac{\partial v_{2}}{\partial x_{1}}\right)+v_{2} \left(\frac{\partial A_{31}}{\partial x_{2}}\right)+A_{33} \left(\frac{\partial v_{3}}{\partial x_{1}}\right)+
		v_{3} \left(\frac{\partial A_{31}}{\partial x_{3}}\right) 
		\\ [3mm] 
		A_{31} \left(\frac{\partial v_{1}}{\partial x_{2}}\right)+v_{1} \left(\frac{\partial A_{32}}{\partial x_{1}}\right)+A_{32} \left(\frac{\partial v_{2}}{\partial x_{2}}\right)+v_{2} \left(\frac{\partial A_{32}}{\partial x_{2}}\right)+A_{33} \left(\frac{\partial v_{3}}{\partial x_{2}}\right)+
		v_{3} \left(\frac{\partial A_{32}}{\partial x_{3}}\right) 
		\\ [3mm] 
		A_{31} \left(\frac{\partial v_{1}}{\partial x_{3}}\right)+v_{1} \left(\frac{\partial A_{33}}{\partial x_{1}}\right)+A_{32} \left(\frac{\partial v_{2}}{\partial x_{3}}\right)+v_{2} \left(\frac{\partial A_{33}}{\partial x_{2}}\right)+A_{33} \left(\frac{\partial v_{3}}{\partial x_{3}}\right)+
		v_{3} \left(\frac{\partial A_{33}}{\partial x_{3}}\right) 
	\end{array}\right).
\end{equation}
Then, the product of the above terms with the dual variables leads to
\begin{equation} \label{eqn.wD}
	\left(\begin{array}{c}
		w_{11} \, D_{11} \\ [2mm]
		w_{12} \, D_{12} \\ [2mm]
		w_{13} \, D_{13} \\ [2mm]
		w_{21} \, D_{21} \\ [2mm]
		w_{22} \, D_{22} \\ [2mm]
		w_{23} \, D_{23} \\ [2mm]
		w_{31} \, D_{31} \\ [2mm]
		w_{32} \, D_{32} \\ [2mm]
		w_{33} \, D_{33} 
		\end{array} \right) = 
	\left(\begin{array}{c}
		\left(A_{22} A_{33}-A_{23} A_{32}\right) \left(A_{11} \left(\frac{\partial v_{1}}{\partial x_{1}}\right)+v_{1} \left(\frac{\partial A_{11}}{\partial x_{1}}\right)+A_{12} \left(\frac{\partial v_{2}}{\partial x_{1}}\right)+v_{2} \left(\frac{\partial A_{11}}{\partial x_{2}}\right)+A_{13} \left(\frac{\partial v_{3}}{\partial x_{1}}\right)+v_{3} \left(\frac{\partial A_{11}}{\partial x_{3}}\right)\right) \\ [2mm]
		\left(-A_{21} A_{33}+A_{23} A_{31}\right) \left(A_{11} \left(\frac{\partial v_{1}}{\partial x_{2}}\right)+v_{1} \left(\frac{\partial A_{12}}{\partial x_{1}}\right)+A_{12} \left(\frac{\partial v_{2}}{\partial x_{2}}\right)+v_{2} \left(\frac{\partial A_{12}}{\partial x_{2}}\right)+A_{13} \left(\frac{\partial v_{3}}{\partial x_{2}}\right)+v_{3} \left(\frac{\partial A_{12}}{\partial x_{3}}\right)\right) \\ [2mm]
		\left(A_{21} A_{32}-A_{22} A_{31}\right) \left(A_{11} \left(\frac{\partial v_{1}}{\partial x_{3}}\right)+v_{1} \left(\frac{\partial A_{13}}{\partial x_{1}}\right)+A_{12} \left(\frac{\partial v_{2}}{\partial x_{3}}\right)+v_{2} \left(\frac{\partial A_{13}}{\partial x_{2}}\right)+A_{13} \left(\frac{\partial v_{3}}{\partial x_{3}}\right)+v_{3} \left(\frac{\partial A_{13}}{\partial x_{3}}\right)\right) \\ [2mm]
		\left(-A_{12} A_{33}+A_{13} A_{32}\right) \left(A_{21} \left(\frac{\partial v_{1}}{\partial x_{1}}\right)+v_{1} \left(\frac{\partial A_{21}}{\partial x_{1}}\right)+A_{22} \left(\frac{\partial v_{2}}{\partial x_{1}}\right)+v_{2} \left(\frac{\partial A_{21}}{\partial x_{2}}\right)+A_{23} \left(\frac{\partial v_{3}}{\partial x_{1}}\right)+v_{3} \left(\frac{\partial A_{21}}{\partial x_{3}}\right)\right) \\ [2mm]
		\left(A_{11} A_{33}-A_{13} A_{31}\right) \left(A_{21} \left(\frac{\partial v_{1}}{\partial x_{2}}\right)+v_{1} \left(\frac{\partial A_{22}}{\partial x_{1}}\right)+A_{22} \left(\frac{\partial v_{2}}{\partial x_{2}}\right)+v_{2} \left(\frac{\partial A_{22}}{\partial x_{2}}\right)+A_{23} \left(\frac{\partial v_{3}}{\partial x_{2}}\right)+v_{3} \left(\frac{\partial A_{22}}{\partial x_{3}}\right)\right) \\ [2mm]
		\left(-A_{11} A_{32}+A_{12} A_{31}\right) \left(A_{21} \left(\frac{\partial v_{1}}{\partial x_{3}}\right)+v_{1} \left(\frac{\partial A_{23}}{\partial x_{1}}\right)+A_{22} \left(\frac{\partial v_{2}}{\partial x_{3}}\right)+v_{2} \left(\frac{\partial A_{23}}{\partial x_{2}}\right)+A_{23} \left(\frac{\partial v_{3}}{\partial x_{3}}\right)+v_{3} \left(\frac{\partial A_{23}}{\partial x_{3}}\right)\right) \\ [2mm]
		\left(A_{12} A_{23}-A_{13} A_{22}\right) \left(A_{31} \left(\frac{\partial v_{1}}{\partial x_{1}}\right)+v_{1} \left(\frac{\partial A_{31}}{\partial x_{1}}\right)+A_{32} \left(\frac{\partial v_{2}}{\partial x_{1}}\right)+v_{2} \left(\frac{\partial A_{31}}{\partial x_{2}}\right)+A_{33} \left(\frac{\partial v_{3}}{\partial x_{1}}\right)+v_{3} \left(\frac{\partial A_{31}}{\partial x_{3}}\right)\right) \\ [2mm]
		\left(-A_{11} A_{23}+A_{13} A_{21}\right) \left(A_{31} \left(\frac{\partial v_{1}}{\partial x_{2}}\right)+v_{1} \left(\frac{\partial A_{32}}{\partial x_{1}}\right)+A_{32} \left(\frac{\partial v_{2}}{\partial x_{2}}\right)+v_{2} \left(\frac{\partial A_{32}}{\partial x_{2}}\right)+A_{33} \left(\frac{\partial v_{3}}{\partial x_{2}}\right)+v_{3} \left(\frac{\partial A_{32}}{\partial x_{3}}\right)\right) \\ [2mm]
		\left(A_{11} A_{22}-A_{12} A_{21}\right) \left(A_{31} \left(\frac{\partial v_{1}}{\partial x_{3}}\right)+v_{1} \left(\frac{\partial A_{33}}{\partial x_{1}}\right)+A_{32} \left(\frac{\partial v_{2}}{\partial x_{3}}\right)+v_{2} \left(\frac{\partial A_{33}}{\partial x_{2}}\right)+A_{33} \left(\frac{\partial v_{3}}{\partial x_{3}}\right)+v_{3} \left(\frac{\partial A_{33}}{\partial x_{3}}\right)\right) \\ [2mm]
	\end{array}\right).	
\end{equation}
On the other side, the non-conservative products in \eqref{eqn.detAnc} explicitly write 
\small
\begin{align}
	 |\A| \frac{\partial v_k}{\partial x_k}&= \left(A_{11} A_{22} A_{33}-A_{11} A_{23} A_{32}-A_{12} A_{21} A_{33}+A_{12} A_{23} A_{31}+A_{13} A_{21} A_{32}-A_{13} A_{22} A_{31}\right) \left(\frac{{\partial v_{1}}}{{\partial x_{1}}} + \frac{{\partial v_{2}}}{{\partial x_{2}}} + \frac{{\partial v_{3}}}{{\partial x_{3}}}\right) \label{eqn.Adivv} \\ 
	 v_1 \frac{\partial |\A|}{\partial x_1}&= v_{1} \left( A_{22} A_{33} \left(\frac{\partial A_{11}}{\partial  x_{1}}\right)+A_{11} A_{22} \left(\frac{\partial A_{33}}{\partial x_{1}}\right)+A_{11} A_{33} \left(\frac{\partial A_{22}}{\partial x_{1}}\right)-A_{23} A_{32} \left(\frac{\partial A_{11}}{\partial x_{1}}\right)-A_{11} A_{23} \left(\frac{\partial A_{32}}{\partial x_{1}}\right)-A_{11} A_{32} \left(\frac{\partial A_{23}}{\partial x_{1}}\right) \right. \nn \\ 
	 & \phantom{v_{1} \left( \right. } -A_{21} A_{33} \left(\frac{\partial A_{12}}{\partial x_{1}}\right)-A_{12} A_{21}  \left(\frac{\partial A_{33}}{\partial x_{1}}\right)-A_{12} A_{33} \left(\frac{\partial A_{21}}{\partial x_{1}}\right)+A_{23} A_{31} \left(\frac{\partial A_{12}}{\partial x_{1}}\right)+A_{12} A_{23} \left(\frac{\partial A_{31}}{\partial x_{1}}\right)+A_{12} A_{31} \left(\frac{\partial A_{23}}{\partial x_{1}}\right) \nn \\
	 & \phantom{v_{1} \left( \right. } \left. +A_{21} A_{32} \left(\frac{\partial A_{13}}{\partial x_{1}}\right)+A_{13} A_{21} \left(\frac{\partial A_{32}}{\partial x_{1}}\right)+A_{13} A_{32} \left(\frac{\partial A_{21}}{\partial x_{1}}\right)-A_{22} A_{31} \left(\frac{\partial A_{13}}{\partial x_{1}}\right)-A_{13} A_{22} \left(\frac{\partial A_{31}}{\partial x_{1}}\right)-A_{13} A_{31} \left(\frac{\partial A_{22}}{\partial x_{1}}\right)\right) \label{eqn.v1divA} \\
	 v_2 \frac{\partial |\A|}{\partial x_2}&= v_{2} \left(A_{22} A_{33} \left(\frac{\partial A_{11}}{\partial x_{2}}\right)+A_{11} A_{22} \left(\frac{\partial A_{33}}{\partial x_{2}}\right)+A_{11} A_{33} \left(\frac{\partial A_{22}}{\partial x_{2}}\right)-A_{23} A_{32} \left(\frac{\partial A_{11}}{\partial x_{2}}\right)-A_{11} A_{23} \left(\frac{\partial A_{32}}{\partial x_{2}}\right)-A_{11} A_{32} \left(\frac{\partial A_{23}}{\partial x_{2}}\right) \right. \nn \\
	 & \phantom{v_{1} \left( \right. } -A_{21} A_{33} \left(\frac{\partial A_{12}}{\partial x_{2}}\right)-A_{12} A_{21} \left(\frac{\partial A_{33}}{\partial x_{2}}\right)-A_{12} A_{33} \left(\frac{\partial A_{21}}{\partial x_{2}}\right)+A_{23} A_{31} \left(\frac{\partial A_{12}}{\partial x_{2}}\right)+A_{12} A_{23} \left(\frac{\partial A_{31}}{\partial x_{2}}\right)+A_{12} A_{31} \left(\frac{\partial A_{23}}{\partial x_{2}}\right) \nn \\
	 & \phantom{v_{1} \left( \right. } \left. +A_{21} A_{32} \left(\frac{\partial A_{13}}{\partial x_{2}}\right)+A_{13} A_{21} \left(\frac{\partial A_{32}}{\partial x_{2}}\right)+A_{13} A_{32} \left(\frac{\partial A_{21}}{\partial x_{2}}\right)-A_{22} A_{31} \left(\frac{\partial A_{13}}{\partial x_{2}}\right)-A_{13} A_{22} \left(\frac{\partial A_{31}}{\partial x_{2}}\right)-A_{13} A_{31} \left(\frac{\partial A_{22}}{\partial x_{2}}\right)\right)  \label{eqn.v2divA} \\
	 v_3 \frac{\partial |\A|}{\partial x_3}&= v_{3} \left(A_{22} A_{33} \left(\frac{\partial A_{11}}{\partial x_{3}}\right)+A_{11} A_{22} \left(\frac{\partial A_{33}}{\partial x_{3}}\right)+A_{11} A_{33} \left(\frac{\partial A_{22}}{\partial x_{3}}\right)-A_{23} A_{32} \left(\frac{\partial A_{11}}{\partial x_{3}}\right)-A_{11} A_{23} \left(\frac{\partial A_{32}}{\partial x_{3}}\right)-A_{11} A_{32} \left(\frac{\partial A_{23}}{\partial x_{3}}\right) \right. \nn \\
	 & \phantom{v_{1} \left( \right. } -A_{21} A_{33} \left(\frac{\partial A_{12}}{\partial x_{3}}\right)-A_{12} A_{21} \left(\frac{\partial A_{33}}{\partial x_{3}}\right)-A_{12} A_{33} \left(\frac{\partial A_{21}}{\partial x_{3}}\right)+A_{23} A_{31} \left(\frac{\partial A_{12}}{\partial x_{3}}\right)+A_{12} A_{23} \left(\frac{\partial A_{31}}{\partial x_{3}}\right)+A_{12} A_{31} \left(\frac{\partial A_{23}}{\partial x_{3}}\right) \nn \\
	 & \phantom{v_{1} \left( \right. } \left. +A_{21} A_{32} \left(\frac{\partial A_{13}}{\partial x_{3}}\right)+A_{13} A_{21} \left(\frac{\partial A_{32}}{\partial x_{3}}\right)+A_{13} A_{32} \left(\frac{\partial A_{21}}{\partial x_{3}}\right)-A_{22} A_{31} \left(\frac{\partial A_{13}}{\partial x_{3}}\right)-A_{13} A_{22} \left(\frac{\partial A_{31}}{\partial x_{3}}\right)-A_{13} A_{31} \left(\frac{\partial A_{22}}{\partial x_{3}}\right)\right) \label{eqn.v3divA} 
\end{align}
\normalsize
After some tedious algebraic manipulations, we arrive at the result
\begin{eqnarray}
	w_{ik} \cdot D_{ik} &= & |\A| \frac{\partial v_k}{\partial x_k} + v_k \frac{\partial |\A|}{\partial x_k}, \\
	\eqref{eqn.wD} &= & \eqref{eqn.Adivv} + \eqref{eqn.v1divA} + \eqref{eqn.v2divA} + \eqref{eqn.v3divA} , \nn
\end{eqnarray}
therefore the GCL \eqref{eqn.detApde} is retrieved as the dot product of the dual variables $\w$ with the evolution equations of the distortion tensor $\A$ given by \eqref{eqn.deformation}. 

\section{Runge-Kutta schemes} \label{app.rk}
Runge-Kutta methods represent a quite popular technique to carry out time integration and they are based on the method of lines (MOL) approach. The governing equations can be written in semi-discrete form as
\begin{equation}
	\frac{\text{d}\U}{\text{d}t} = \mathcal{L}_h(\U),
\end{equation}
where $\mathcal{L}_h(\U)$ contains the spatial discretization of the numerical fluxes, non-conservative products and source terms. A generic Runge-Kutta scheme with a total number of $s$ sub-stages is described by a Butcher tableau of the form shown in Table \ref{tab.BTRK}.
\begin{table}[h!]
	\caption{Butcher tableau for Runge-Kutta explicit methods.}
	\begin{center} 
		\begin{tabular}{c|ccccc} 
			0 &  & & & & \\
			$\alpha_2$ & $\beta_{21}$ & & & & \\
			$\alpha_3$ & $\beta_{31}$ & $\beta_{32}$ & & & \\
			$\vdots$  & $\vdots$      & $\vdots$  & $\ddots$ & & \\
			$\alpha_s$ & $\beta_{s1}$ & $\beta_{s2}$ & $...$ & $\beta_{s \left(s-1\right)}$ & \\
			\hline
			& $c_1$ & $c_2$ & $...$ & $c_{s-1}$ & $c_s$
			\label{tab.BTRK}
		\end{tabular}		
	\end{center}	
\end{table}	
The numerical solution is determined at the next time step as
\begin{equation}
	\U^{n+1} = \U^{n} + \dt \ \sum \limits_{i=1}^s c_i \ \kappa_i.
\end{equation}
The generic Runge-Kutta stage $\kappa_i$ is evaluated at the intermediate time level $t^{\left(i\right)}=t^n+\alpha_i \dt$ by
\begin{equation}
	\kappa_i = \mathcal{L}_h\left( \U_h^n + \dt \sum \limits_{j=1}^{i} \beta_{ij} \ \kappa_j \right)
\end{equation}
with $\U_h^n$ denoting the numerical solution at the current time level $t^n$. In this work we consider three different Runge-Kutta schemes:
\begin{itemize}
	\item Euler method with accuracy $\mathcal{O}(1)$
	\begin{center} 
		\begin{tabular}{c|c} 
			0 & 0 \\
			\hline
			& 1
			\label{tab.RK1}
		\end{tabular}		
	\end{center}
	\item Heun method with accuracy $\mathcal{O}(2)$
	\begin{center} 
		\begin{tabular}{c|cc} 
			0 & 0 & 0 \\
			1 & 1 & 0 \\
			\hline
			& 1/2 & 1/2  
			\label{tab.RK2}
		\end{tabular}		
	\end{center}
	\item RK4 method with accuracy $\mathcal{O}(4)$
	\begin{center} 
		\begin{tabular}{c|cccc} 
			0 & 0 & 0 & 0 & 0  \\
			1/2 & 1/2 & 0 & 0 & 0  \\
			1/2 & 0 & 1/2 & 0 & 0  \\
			1 & 0 & 0 & 1 & 0  \\
			\hline
			& 1/6 & 1/3 & 1/3 & 1/6
			\label{tab.RK4}
		\end{tabular}		
	\end{center}	
\end{itemize}

\end{document}